\documentclass{amsart}

\usepackage{amssymb}
\usepackage{amsmath}
\usepackage{amsthm}
\usepackage{amsfonts}

\usepackage{a4wide}

\usepackage{geometry}
\usepackage{longtable}
\usepackage{multirow}
\usepackage{setspace}

\newtheorem{proposition}{Proposition}
\newtheorem{theorem}{Theorem}
\newtheorem{corollary}{Corollary}

\theoremstyle{definition}
\newtheorem{definition}{Definition}

\newtheorem{remark}{Remark}

\newtheorem{example}{Example}

\newtheorem{conjecture}{Conjecture}

\numberwithin{equation}{section}

\begin{document}


\title[Classifying multiplets of totally real cubic fields]{Classifying multiplets of totally real cubic fields}

\author{Daniel C. Mayer}
\address{Naglergasse 53\\8010 Graz\\Austria}
\email{algebraic.number.theory@algebra.at}
\urladdr{http://www.algebra.at}

\thanks{Research supported by the Austrian Science Fund (FWF): projects J0497-PHY and P26008-N25}



\date{Wednesday, February 24, 2021}


\begin{abstract}
The number of non-isomorphic cubic fields \(L\) sharing a common discriminant \(d_L=d\)
is called the multiplicity \(m=m(d)\) of \(d\).
For an assigned value of \(d\), these fields are collected in a multiplet
\(\mathbf{M}_d=(L_1,\ldots,L_m)\).
In this paper,
the information in all existing tables of totally real cubic number fields \(L\)
with positive discriminants \(d_L<10^7\)
is extended by computing the differential principal factorization types
\(\tau(L)\in\lbrace\alpha_1,\alpha_2,\alpha_3,\beta_1,\beta_2,\gamma,\delta_1,\delta_2,\varepsilon\rbrace\)
of the members \(L\) of each multiplet \(\mathbf{M}_d\) of non-cyclic fields,
a new kind of arithmetical invariants
which provide succinct information about
ambiguous principal ideals and capitulation
in the normal closures \(N\) of non-Galois cubic fields \(L\).
The classification is arranged with respect to
increasing \(3\)-class rank \(\varrho\) of the quadratic subfields \(K\)
of the \(S_3\)-fields \(N\),
and to ascending number of prime divisors of the conductor \(f\) of \(N/K\).
The Scholz conjecture
concerning the distinguished index of subfield units \((U_N:U_0)=1\) 
for ramified extensions \(N/K\) with conductor \(f>1\)
is refined and verified.
\end{abstract}

\maketitle


\section{Introduction}
\label{s:Intro}

\noindent
Recently, we have classified all complex and totally real cubic number fields \(L\)
with discriminants in the range \(-20\,000<d_L<100\,000\),
covered in the years between \(1972\) and \(1976\) by Angell
\cite{An1975,An1976}.
The reconstruction
\cite{Ma2021}
was carried out from the viewpoint of \(3\)-ring class fields
with Magma
\cite{MAGMA2020}.

In this paper, we omit simply real cubic fields with a few types only,
and we rather put our focus on \textit{triply real} cubic fields \(L\) with \textit{nine} possible types
\(\tau(L)\in\lbrace\alpha_1,\alpha_2,\alpha_3,\beta_1,\beta_2,\gamma,\delta_1,\delta_2,\varepsilon\rbrace\),
which refine the coarse classification with five types \(\alpha,\beta,\gamma,\delta,\varepsilon\) by Moser
\cite{Mo1979}.
According to the historical development
of systematic investigations of triply real cubic fields,
we present our refined classification of \textit{multiplets} of totally real cubic fields \(L\)
in four steps with increasing upper bounds \(100\,000,200\,000,500\,000\), and \(10\,000\,000\) for the discriminants \(d_L\).

In \S\
\ref{s:Angell},
we start by recalling our results in
\cite{Ma2021}
concerning the range \(d_L<100\,000\) of Angell
\cite{An1975,An1976}.
In \S\
\ref{s:GutensteinStreiteben},
we continue with an update of our extension to \(d_L<200\,000\) in
\cite{Ma1991b,Ma1991c},
which was computed in \(1991\) by means of the Voronoi algorithm
\cite{Vo1896}.
Whereas the count of discriminants and fields and the collection of fields into multiplets were correct,
the classification into type \(\delta_2\) instead of \(\beta_2\) was partially erroneous,
because \textit{absolute principal factors} do not necessarily show up
among the lattice minima in the chains of Voronoi's algorithm.
The coronation of this paper will be established in \S\S\
\ref{s:EnnolaTurunen}
and
\ref{s:LlorenteQuer},
where we use Fieker's class field routines
\cite{Fi2001}
in Magma
\cite{MAGMA2020}
to classify the range \(d_L<500\,000\) of Ennola and Turunen
\cite{EnTu1983,EnTu1985}
and the range \(d_L<10\,000\,000\) of Llorente and Quer
\cite{LlQu1988},
the most extensive ranges deposited in files of unpublished mathematical tables (UMT).
Statistical evaluations and theoretical conclusions are given in \S\
\ref{s:Theoretical}.

We mention that Belabas
\cite{Be1997}
has given a fast method
for simply \textit{counting} cubic fields even in bigger ranges,
without computation of \textit{arithmetical invariants},
like fundamental systems of units and class group structures,
and without classification into \textit{differential principal factorization types}. 

The latter are introduced in \S\
\ref{s:DPFTypes},
where we prove that the \(\mathbb{F}_p\)-vector space \(\mathcal{P}_{N/K}/\mathcal{P}_K\)
of primitive ambiguous principal ideals of a number field extension \(N/K\)
with odd prime degree \(p\)
can be endowed with a natural trichotomic direct product structure.

Finally,
as an application of the notions of multiplets and DPF types,
the Scholz conjecture
\cite{So1933}
concerning the normal closure \(N\) of \(L\)
is stated, refined, and proved completely in \S\
\ref{s:ScholzConjecture}.


We clarify in advance
why we present the results of our classification
for several intervals \(0<d_L<B\) of positive cubic discriminants
instead of focussing on the most extensive of these ranges:
Firstly, we want to add value to the classical tables of
Angell, Ennola and Turunen, Llorente and Quer,
which did not illuminate the constitution of multiplets \((L_1,\ldots,L_m)\) of totally real cubic fields
as subfields of a common \(3\)-ring class field \(K_f\) modulo a \(3\)-admissible conductor \(f\)
over a quadratic base field \(K\),
let alone the differential principal factorization types \((\tau(L_1),\ldots,\tau(L_m))\). 
Secondly, it is our intention to show the increasing \textit{wealth of arithmetical structure}
in three successive extensions of the upper bound \(B\)
from \(100\,000\) to \(200\,000\), \(500\,000\), and finally \(10\,000\,000\).
There arise conductors \(f\) divisible by an ascending number of primes,
new types \(\tau(L)\),
and \textit{heterogeneous} multiplets
\(\mathbf{M}(K_f)=\lbrack\mathbf{M}_{c^2d}\rbrack_{c\mid f}=\lbrack(L_{c,1},\ldots,L_{c,m(c)})\rbrack_{c\mid f}\)
with increasing complexity.


\section{Construction as subfields of a ring class field}
\label{s:RingClassField}

\subsection{Structure and multiplicity of cubic discriminants}
\label{ss:Discriminants}

\noindent
Let \(K\) be a quadratic number field
with fundamental discriminant \(d=d_K\)
(square free, except possibly for the \(2\)-contribution).

\begin{definition}
\label{dfn:Admissible}
A positive integer \(f\) is called a \(3\)-\textit{admissible conductor} for \(K\),
if it has the shape \(f=3^e\cdot q_1\cdots q_t\)
with an integer exponent \(e\in\lbrace 0,1,2\rbrace\),
\(t\ge 0\),
and pairwise distinct prime numbers \(q_1,\ldots,q_t\in\mathbb{P}\setminus\lbrace 3\rbrace\),
such that the following conditions are satisfied:
\[
\text{Kronecker symbol } \left(\frac{d}{q_i}\right)\equiv q_i\,(\mathrm{mod}\,3), \text{ for all } 1\le i\le t, 
\]
and
\[
e\in
\begin{cases}
\lbrace 0,2\rbrace   & \text{ if } d\equiv\pm 1\,(\mathrm{mod}\,3), \\
\lbrace 0,1\rbrace   & \text{ if } d\equiv +3\,(\mathrm{mod}\,9), \\
\lbrace 0,1,2\rbrace & \text{ if } d\equiv -3\,(\mathrm{mod}\,9).
\end{cases}
\]
\end{definition}

\noindent
So, a \(3\)-admissible conductor \(f\) for \(K\) is essentially square free,
except possibly for the \textit{critical} contribution by the prime \(3\).
The condition involving the Kronecker symbol means that
a non-critical prime divisor \(q\ne 3\) of \(f\)
must remain inert in \(K\), if \(q\equiv -1\,(\mathrm{mod}\,3)\),
and must split in \(K\), if \(q\equiv +1\,(\mathrm{mod}\,3)\).
The critical prime divisor \(3\) of \(f\) is \(3\)-admissible
if and only if it ramifies in \(K\), that is, if \(d\equiv\pm 3\,(\mathrm{mod}\,9)\).
Otherwise only the critical prime power divisor \(9\) of \(f\) is \(3\)-admissible.
(Recall that \(3\) remains inert in \(K\) if \(d\equiv -1\,(\mathrm{mod}\,3)\)
and \(3\) splits in \(K\) if \(d\equiv +1\,(\mathrm{mod}\,3)\).)
So far, all contributions to \(f\) are regular.
There is, however, the possibility of an \textit{irregular}
\(3\)-admissible critical prime power divisor \(9\) of \(f\),
when \(d\equiv -3\,(\mathrm{mod}\,9)\).


\begin{definition}
\label{dfn:Formal}
An integer \(D=f^2\cdot d\) is called a \textit{formal cubic discriminant}
if \(f\) is a \(3\)-admissible conductor for the quadratic field \(K\)
with fundamental discriminant \(d\).
(Since the square \(f^2\) and the fundamental discriminant \(d\)
are congruent to \(0\) or \(1\) modulo \(4\),
this is also the case for a formal cubic discriminant \(D\).
We shall see that \(D\) is not necessarily discriminant \(d_L\) of a cubic field \(L\).)
\end{definition}

\noindent
Note that this definition does not include discriminants \(d_L\)
of cyclic cubic fields \(L\)
which are perfect squares \(f^2\) of conductors \(f\)
exactly divisible by primes congruent to \(1\) modulo \(3\),
and possibly also by the prime power \(3^2\).
In order to determine the multiplicity of \(d_L\), we need further definitions.


\begin{definition}
\label{dfn:SelmerSpace}
An algebraic number \(\alpha\ne 0\) in the quadratic field \(K\)
is called a \(3\)-\textit{virtual unit},
if its principal ideal \(\alpha\mathcal{O}_K\) is the cube \(\mathfrak{j}^3\)
of an ideal \(\mathfrak{j}\) of \(K\).
Obviously all units \(\eta\) in \(U_K\) and all third powers \(\alpha^3\ne 0\) in \((K^\times)^3\)
are \(3\)-virtual units of \(K\).
Let \(I\) denote the group of all \(3\)-virtual units of \(K\),
and let \(K^\times=K\setminus\lbrace 0\rbrace\) denote the multiplication group of \(K\).
The \(\mathbb{F}_3\)-vector space \(V:=I/(K^\times)^3\)
is called the \(3\)-\textit{Selmer space} of \(K\).
\end{definition}

\noindent
For any positive integer \(n\)
and a set \(X\) of algebraic numbers,
let \(X(n)\) be the subset of \(X\)
consisting of elements coprime to \(n\).
The \(3\)-Selmer space \(V\) of \(K\)
is isomorphic to the direct product
of the \(3\)-elementary class group
\(\mathrm{Cl}_K/\mathrm{Cl}_K^3\)
and the \(3\)-elementary unit group
\(U_K/U_K^3\)
of \(K\)
\cite[p. 2212]{Ma2014}.
Since any ideal class contains an ideal coprime
to an assigned positive integer \(n\), it follows that
\(\mathrm{Cl}_K=\mathcal{I}_K/\mathcal{P}_K\simeq\mathcal{I}_K(n)/\mathcal{P}_K(n)\),
and since trivially \(U_K=U_K(n)\),
we have \(V\simeq I(n)/K(n)^3\)
\cite[Dfn. 2.2, p. 2211]{Ma2014}.

\begin{definition}
\label{dfn:RingSpace}
Let \(f\) be a positive integer.
Denote by \(S_f:=\lbrace\alpha\in K\mid\alpha\equiv 1\,(\mathrm{mod}\,f)\rbrace\)
the \textit{ray} modulo \(f\) of \(K\),
and by \(R_f:=\mathbb{Q}(f)\cdot S_f\)
the \textit{ring} modulo \(f\) of \(K\).
The subspace
\(V(f):=(I(f)\cap R_f\cdot K(f)^3)/K(f)^3\)
of the \(3\)-Selmer space \(V\) is called the \(3\)-\textit{ring space} modulo \(f\) of \(K\).
Its codimension \(\delta(f):=\mathrm{codim}(V(f))=\dim_{\mathbb{F}_3}(V/V(f))\)
is called the \(3\)-\textit{defect} of \(f\) with respect to \(K\).
\end{definition}

\noindent
In order to enable comparison and binary operations (in particular, intersection)
of two different ring spaces \(V(f)\) and \(V(f^\prime)\),
we need the concept of a \textit{modulus of declaration},
that is a positive integer \(n\) which is a common multiple of \(f\) and \(f^\prime\).
Then
\(V(f)\simeq (I(n)\cap R_f\cdot K(n)^3)/K(n)^3\)
and \(V(f^\prime)\simeq (I(n)\cap R_{f^\prime}\cdot K(n)^3)/K(n)^3\),
whence it makes sense to speak about inclusion and meet.

\begin{remark}
\label{rmk:RingSpace}
It is possible to avoid the requirement of a modulus of declaration,
if the theory of ring spaces is based on the approach via \textit{id\`ele groups}.
This has been done by Satg\'e
\cite{Sa1981}
for prime conductors \(f=q\) and will be expanded further by ourselves
for any \(f\) in a future paper.
\end{remark}


\subsection{Homogeneous and heterogeneous multiplets}
\label{ss:HomoAndHetero}

\noindent
If \(f\) is a \(3\)-admissible conductor
with \(3\)-defect \(\delta(f)\)
for a quadratic field \(K\)
with fundamental discriminant \(d\)
and \(3\)-class rank \(\varrho\),
then the sum of all multiplicities \(m(D)\) of formal cubic discriminants \(D=c^2d\)
with \(c\) running over all divisors of \(f\)
is given by
\[
\sum_{c\mid f}\,m(c^2d)=\frac{1}{2}(3^{\varrho_f}-1)
\]
in terms of the \(3\)-\textit{ring class rank} modulo \(f\) of \(K\)
\cite[Thm. 2.1, p. 2213]{Ma2014},
\[
\varrho_f=\varrho+t+w-\delta(f),
\]
where \(t:=\#\lbrace q\in\mathbb{P}\setminus\lbrace 3\rbrace\mid v_q(f)=1\rbrace\),
and \(w\) is defined in terms of the \(3\)-valuation \(v_3(f)\) of \(f\),
\[
w:=
\begin{cases}
0 & \text{ if } v_3(f)=0, \\
1 & \text{ if } v_3(f)=1 \text{ or } \lbrack v_3(f)=2 \text{ and } d\equiv\pm 1\,(\mathrm{mod}\,3)\rbrack, \\
2 & \text{ if } v_3(f)=2 \text{ and } d\equiv 6\,(\mathrm{mod}\,9).
\end{cases}
\]

\begin{definition}
\label{dfn:Multiplets}
Let \(f\) be a \(3\)-admissible conductor for a quadratic field \(K\).
\begin{enumerate}
\item
For each divisor \(c\) of \(f\) which is also a \(3\)-admissible conductor for \(K\),
the multiplet
\[
\mathbf{M}_{c^2d}:=(L_{c,1},\ldots,L_{c,m}) \text{ with } m=m(c^2d)
\]
is called the \textit{homogeneous} multiplet 
of cubic fields \(L_{c,i}\) with discriminant \(c^2d\).
\item
The multiplet
\(\mathbf{M}(K_f):=\lbrack\mathbf{M}_{c^2d}\rbrack_{c\mid f}\)
is called the \textit{heterogeneous} multiplet
of the \(3\)-ring class field \(K_f\) modulo \(f\) of \(K\).
(The normal closures of all cubic fields \(L_{c,i}\)
with \(c\mid f\) and \(1\le i\le m(c^2d)\)
are subfields of the ring class field \(K_f\).)
\item
The family \(\mathrm{sgn}(\mathbf{M}(K_f)):=\lbrack m(c^2d)\rbrack_{c\mid f}\) of all partial multiplicities associated with \(f\)
is called the \textit{signature} of the heterogeneous multiplet \(\mathbf{M}(K_f)\).
\end{enumerate}
\end{definition}

\noindent
\(D=c^2d\) is only a \textit{formal} but not an actual cubic discriminant
if and only if the multiplicity \(m(c^2d)=0\) vanishes,
that is, if \(\mathbf{M}_{c^2d}=\emptyset\) is a \textit{nilet}
(denoted by the empty set symbol).

\begin{definition}
\label{dfn:Type}
By the \textit{type of the \(3\)-ring class field} \(K_f\) modulo \(f\) of \(K\)
we understand the following pair \((\mathrm{Obj}(K_f),\mathrm{Inv}(K_f))\) of heterogeneous multiplets
\begin{equation}
\label{eqn:Type}
\begin{aligned}
\mathrm{Obj}(K_f) &:= \mathbf{M}(K_f)=\lbrack(L_{c,i})_{1\le i\le m(c^2d)}\rbrack_{c\mid f} \\
\mathrm{Inv}(K_f) &:= \tau(\mathbf{M}(K_f))=\lbrack(\tau(L_{c,i}))_{1\le i\le m(c^2d)}\rbrack_{c\mid f}
\end{aligned}
\end{equation}
consisting of all non-cyclic cubic fields \(L_{c,i}\) with discriminants \(c^2d\) dividing \(f^2d\) as \textit{objects}
and their differential principal factorization types \(\tau(L_{c,i})\) as \textit{invariants}. \\
(See
\cite{Ma2019a,Ma2019b}
and the next section \S\
\ref{s:DPFTypes}.)
\end{definition}


\subsection{Algorithmic process of construction}
\label{ss:Construction}

\noindent
The computational technique which will be employed
for the construction of totally real cubic fields in the sections \S\S\
\ref{s:Angell},
\ref{s:GutensteinStreiteben},
\ref{s:EnnolaTurunen}, and
\ref{s:LlorenteQuer}
consists of two steps.
For an assigned real quadratic field \(K\) with fundamental discriminant \(d\)
and a \(3\)-admissible conductor \(f\) for \(K\),
initially all cyclic cubic extensions \(N/K\) with conductor \(f\)
are constructed as subfields of the \textit{ray class field} modulo \(f\) of \(K\).
Then the members \(N\) of this family are tested for their
absolute automorphism group \(G=\mathrm{Gal}(N/\mathbb{Q})\),
and only those with \(G\simeq S_3\) are permitted to pass the filter.
As a double check, we additionally make sure that
the non-Galois subfields \(L<N\) have the required discriminant \(d_L=f^2\cdot d\),
and thus \(N\) is subfield of \(K_f\), the \textit{ring class field} modulo \(f\) of \(K\),
which is contained in  the ray class field modulo \(f\) of \(K\).
The result is the multiplet \(\mathbf{M}_{f^2d}\),
because the fields \(N\) are certainly not subfields of \(K_c\) for proper divisors \(c\) of \(f\).

Before we apply this algorithm, however, we have to introduce the concept of
differential principal factorizations (DPF) in section \S\
\ref{s:DPFTypes}.


\section{Differential principal factorization types}
\label{s:DPFTypes}

\begin{figure}[ht]
\caption{Dihedral and metacyclic situation}
\label{fig:DihedralMetacyclic}

{\small

\setlength{\unitlength}{1.0cm}
\begin{picture}(12,5)(-7,-9)



\put(-6,-9){\circle*{0.2}}
\put(-6,-9.2){\makebox(0,0)[ct]{\(\mathbb{Q}\)}}
\put(-7,-9){\makebox(0,0)[rc]{rational field}}

\put(-6,-9){\line(2,1){2}}
\put(-5,-8.7){\makebox(0,0)[lt]{\(\lbrack K:\mathbb{Q}\rbrack=2\)}}

\put(-4,-8){\circle*{0.2}}
\put(-4,-8.2){\makebox(0,0)[lt]{\(K=\mathbb{Q}(\sqrt{d})\)}}
\put(-3,-8){\makebox(0,0)[lc]{quadratic field}}


\put(-6.2,-7.5){\makebox(0,0)[rc]{\(\lbrack L:\mathbb{Q}\rbrack=p\)}}
\put(-6,-9){\line(0,1){3}}
\put(-4,-8){\line(0,1){3}}
\put(-3.8,-6.5){\makebox(0,0)[lc]{cyclic extension}}



\put(-6,-6){\circle{0.2}}
\put(-6,-5.8){\makebox(0,0)[cb]{\(L\)}}
\put(-6.1,-6.2){\makebox(0,0)[rt]{\(L_{1},\ldots,L_{p-1}\)}}
\put(-7,-6){\makebox(0,0)[rc]{\(p\) conjugates}}

\put(-6,-6){\line(2,1){2}}

\put(-4,-5){\circle*{0.2}}
\put(-4,-4.8){\makebox(0,0)[cb]{\(N=L\cdot K\)}}
\put(-3,-5){\makebox(0,0)[lc]{dihedral field}}
\put(-3,-5.5){\makebox(0,0)[lc]{of degree \(2p\)}}

\put(-0.9,-9.4){\line(0,1){5}}


\put(2,-9){\circle*{0.2}}
\put(2,-9.2){\makebox(0,0)[ct]{\(\mathbb{Q}\)}}
\put(1,-9){\makebox(0,0)[rc]{rational field}}

\put(2,-9){\line(2,1){2}}
\put(3,-8.7){\makebox(0,0)[lt]{\(\lbrack K:\mathbb{Q}\rbrack=p-1\)}}
\put(3,-8.5){\circle*{0.2}}

\put(4,-8){\circle*{0.2}}
\put(4,-8.2){\makebox(0,0)[lt]{\(K=\mathbb{Q}(\zeta_p)\)}}
\put(5,-8){\makebox(0,0)[lc]{cyclotomic field}}


\put(1.8,-7.5){\makebox(0,0)[rc]{\(\lbrack L:\mathbb{Q}\rbrack=p\)}}
\put(2,-9){\line(0,1){3}}
\put(3,-6.5){\vector(0,1){0.8}}
\put(3,-6.8){\makebox(0,0)[cc]{intermediate}}
\put(3,-7.2){\makebox(0,0)[cc]{fields}}
\put(3,-7.5){\vector(0,-1){0.8}}
\put(4,-8){\line(0,1){3}}
\put(4.2,-6.5){\makebox(0,0)[lc]{Kummer extension}}



\put(2,-6){\circle{0.2}}
\put(1.8,-5.7){\makebox(0,0)[rb]{\(L=\mathbb{Q}(\sqrt[p]{D})\)}}
\put(1.9,-6.2){\makebox(0,0)[rt]{\(L_{1},\ldots,L_{p-1}\)}}
\put(1,-6){\makebox(0,0)[rc]{\(p\) conjugates}}

\put(2,-6){\line(2,1){2}}
\put(3,-5.5){\circle{0.2}}

\put(4,-5){\circle*{0.2}}
\put(4,-4.8){\makebox(0,0)[cb]{\(N=L\cdot K\)}}
\put(5,-5){\makebox(0,0)[lc]{metacyclic field}}
\put(5,-5.5){\makebox(0,0)[lc]{of degree \((p-1)p\)}}


\end{picture}

}

\end{figure}
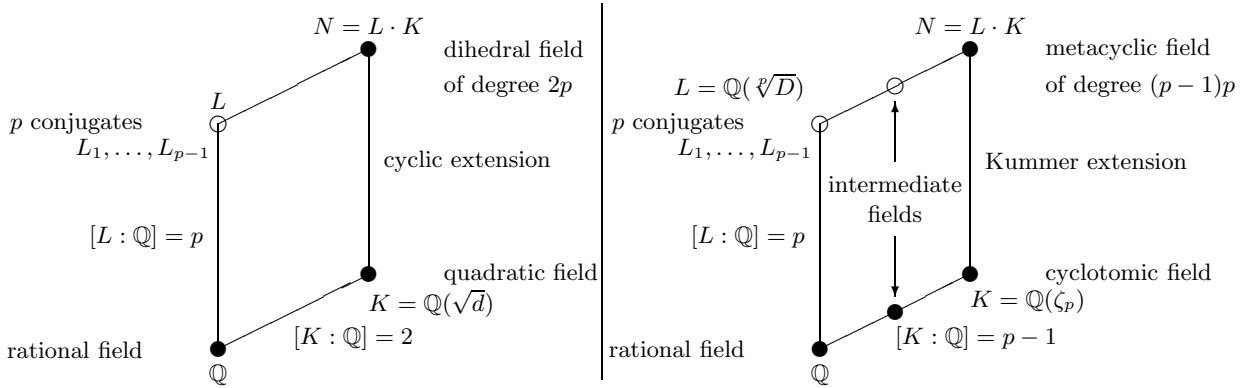

\noindent
Our intention in this section is
to prepare sound foundations for the concept of
\textit{differential principal factorization} (DPF) \textit{types}
and to establish
a common theoretical framework for
the classification
\begin{itemize}
\item
of \textit{dihedral} fields \(N/\mathbb{Q}\) of degree \(2p\) with an odd prime \(p\),
viewed as subfields of suitable \(p\)-ring class fields over a quadratic field \(K\)
(see the left part of Figure
\ref{fig:DihedralMetacyclic}),
and
\item
of \textit{pure metacyclic} fields \(N=K(\sqrt[p]{D})\) of degree \((p-1)\cdot p\) with an odd prime \(p\),
viewed as Kummer extensions of a cyclotomic field \(K=\mathbb{Q}(\zeta_p)\)
(see the right part of Figure
\ref{fig:DihedralMetacyclic}),
\end{itemize}
by the following arithmetical invariants:
\begin{enumerate}
\item
the \(\mathbb{F}_p\)-dimensions of subspaces
of the space \(\mathcal{P}_{N/K}/\mathcal{P}_K\) of primitive ambiguous principal ideals,
which are also called \textit{differential principal factors}, of \(N/K\),
\item
the \textit{capitulation kernel} \(\ker(T_{N/K})\)
of the transfer homomorphism
\(T_{N/K}:\,\mathrm{Cl}_p(K)\to\mathrm{Cl}_p(N)\),
\(\mathfrak{a}\mathcal{P}_K\mapsto(\mathfrak{a}\mathcal{O}_N)\mathcal{P}_N\),
of \(p\)-classes from \(K\) to \(N\), and
\item
the \textit{Galois cohomology} \(\hat{\mathrm{H}}^0(G,U_N)\), \(\mathrm{H}^1(G,U_N)\) of the unit group \(U_N\)
as a module over the cyclic automorphism group \(G=\mathrm{Gal}(N/K)\simeq C_p\).
\end{enumerate}

\subsection{Primitive ambiguous ideals}
\label{ss:PrmAmbIdl}
\noindent
Let \(p\ge 2\) be a prime number,
and \(N/K\) be a relative extension of number fields with degree \(p\)
(\textit{not} necessarily Galois).

\begin{definition}
\label{dfn:Ambiguous}
The group \(\mathcal{I}_N\) of fractional ideals of \(N\)
contains the \textit{subgroup of ambiguous ideals} of \(N/K\),
denoted by the symbol
\(\mathcal{I}_{N/K}:=\lbrace\mathfrak{A}\in\mathcal{I}_N\mid\mathfrak{A}^p\in\mathcal{I}_K\rbrace\).
The quotient \(\mathcal{I}_{N/K}/\mathcal{I}_K\) is called
the \(\mathbb{F}_p\)-\textit{vector space of primitive ambiguous ideals} of \(N/K\).
(Cfr.
\cite[Dfn. 3.1, p. 1991]{Ma2019a}.)
\end{definition}

\begin{proposition}
\label{prp:PrimitiveAmbiguous}
Let \(\mathfrak{L}_1,\ldots,\mathfrak{L}_t\) be the totally ramified prime ideals of \(N/K\),
then a basis and the dimension of the space \(\mathcal{I}_{N/K}/\mathcal{I}_K\) over \(\mathbb{F}_p\) are finite and given by
\begin{equation}
\label{eqn:PrimitiveAmbiguous}
\mathcal{I}_{N/K}/\mathcal{I}_K
\simeq\prod_{i=1}^t\,(\langle\mathfrak{L}_i\rangle/\langle\mathfrak{L}_i^p\rangle)
\simeq\mathbb{F}_p^t, \quad
\dim_{\mathbb{F}_p}(\mathcal{I}_{N/K}/\mathcal{I}_K)=t,
\end{equation}
whereas \(\mathcal{I}_{N/K}\) is an \textit{infinite} abelian group containing \(\mathcal{I}_K\).
\end{proposition}

\begin{proof}
According to the definition of \(\mathcal{I}_{N/K}\),
the quotient \(\mathcal{I}_{N/K}/\mathcal{I}_K\) is an \textit{elementary} abelian \(p\)-group.
By the decomposition law for prime ideals of \(K\) in \(N\), the space
\(\mathcal{I}_{N/K}/\mathcal{I}_K\)
is generated by the \textit{totally ramified} prime ideals (with ramification index \(e=p\)) of \(N/K\),
that is to say
\(\mathcal{I}_{N/K}=\langle\mathfrak{L}\in\mathbb{P}_N\mid\mathfrak{L}^p\in\mathbb{P}_K\rangle\mathcal{I}_K\).
According to the theorem on prime ideals dividing the discriminant,
the number \(t\) of totally ramified prime ideals \(\mathfrak{L}_1,\ldots,\mathfrak{L}_t\) of \(N/K\)
is \textit{finite}.
\end{proof}

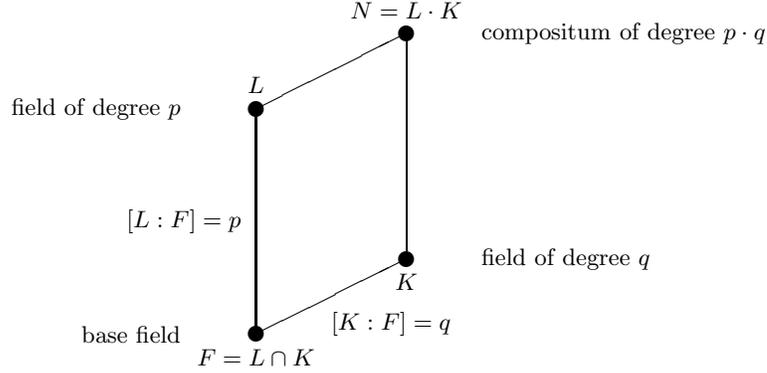
\begin{figure}[ht]
\caption{Hasse subfield diagram of \(N/F\)}
\label{fig:RelativeSubfields}

{\small

\setlength{\unitlength}{1.0cm}
\begin{picture}(5,5)(-7,-9)



\put(-6,-9){\circle*{0.2}}
\put(-6,-9.2){\makebox(0,0)[ct]{\(F=L\cap K\)}}
\put(-7,-9){\makebox(0,0)[rc]{base field}}

\put(-6,-9){\line(2,1){2}}
\put(-5,-8.7){\makebox(0,0)[lt]{\(\lbrack K:F\rbrack=q\)}}

\put(-4,-8){\circle*{0.2}}
\put(-4,-8.2){\makebox(0,0)[ct]{\(K\)}}
\put(-3,-8){\makebox(0,0)[lc]{field of degree \(q\)}}


\put(-6.2,-7.5){\makebox(0,0)[rc]{\(\lbrack L:F\rbrack=p\)}}
\put(-6,-9){\line(0,1){3}}
\put(-4,-8){\line(0,1){3}}



\put(-6,-6){\circle*{0.2}}
\put(-6,-5.8){\makebox(0,0)[cb]{\(L\)}}
\put(-7,-6){\makebox(0,0)[rc]{field of degree \(p\)}}

\put(-6,-6){\line(2,1){2}}

\put(-4,-5){\circle*{0.2}}
\put(-4,-4.8){\makebox(0,0)[cb]{\(N=L\cdot K\)}}
\put(-3,-5){\makebox(0,0)[lc]{compositum of degree \(p\cdot q\)}}


\end{picture}

}

\end{figure}

If \(L\) is another subfield of \(N\)
such that \(N=L\cdot K\) is the compositum of \(L\) and \(K\),
and \(N/L\) is of degree \(q\) \textit{coprime} to \(p\),
then the relative norm homomorphism \(N_{N/L}\) induces an \textit{epimorphism}
\begin{equation}
\label{eqn:InducedNorm}
N_{N/L}:\,\mathcal{I}_{N/K}/\mathcal{I}_K\to\mathcal{I}_{L/F}/\mathcal{I}_F,
\end{equation}
where \(F:=L\cap K\) denotes the intersection of \(L\) and \(K\) in Figure
\ref{fig:RelativeSubfields}.
Thus, by the isomorphism theorem (see also
\cite[Thm. 4.2, pp. 1995--1996]{Ma2019a}),
we have proved:

\begin{theorem}
\label{thm:QualitativeDichotomy}
There are the following two isomorphisms between finite \(\mathbb{F}_p\)-vector spaces:
\begin{equation}
\label{eqn:Dichotomy}
\begin{aligned}
(\mathcal{I}_{N/K}/\mathcal{I}_K)/\ker(N_{N/L}) &\simeq \mathcal{I}_{L/F}/\mathcal{I}_F \quad \text{(quotient)}, \\
\mathcal{I}_{N/K}/\mathcal{I}_K &\simeq (\mathcal{I}_{L/F}/\mathcal{I}_F)\times\ker(N_{N/L}) \quad \text{(direct product)}.
\end{aligned}
\end{equation}
\end{theorem}


\begin{definition}
\label{dfn:Dichotomy}
Since the relative different of \(N/K\) is \textit{essentially} given by
\(\mathfrak{D}_{N/K}=\prod_{i=1}^t\,\mathfrak{L}_i^{p-1}\)
\cite[Thm. 3.2, p. 1993]{Ma2019a},
the space \(\mathcal{I}_{N/K}/\mathcal{I}_K\simeq\prod_{i=1}^t\,(\langle\mathfrak{L}_i\rangle/\langle\mathfrak{L}_i^p\rangle)\)
of primitive ambiguous ideals of \(N/K\) is also called
the space of \textit{differential factors} of \(N/K\).
The two subspaces in the direct product decomposition of
\(\mathcal{I}_{N/K}/\mathcal{I}_K\) in formula
\eqref{eqn:Dichotomy}
are called \\
\(\bullet\) subspace \(\mathcal{I}_{L/F}/\mathcal{I}_F\) of \textit{absolute} differential factors of \(L/F\), and \\
\(\bullet\) subspace \(\ker(N_{N/L})\) of \textit{relative} differential factors of \(N/K\).
\end{definition}


\subsection{Splitting off the norm kernel}
\label{ss:NormKernel}
\noindent
The second isomorphism in formula
\eqref{eqn:Dichotomy}
gives rise to a \textbf{dichotomic decomposition}
of the space \(\mathcal{I}_{N/K}/\mathcal{I}_K\) of primitive ambiguous ideals of \(N/K\)
into two components, whose dimensions can be given
under the following conditions:

\begin{theorem}
\label{thm:QuantitativeDichotomy}
Let \(p\ge 3\) be an odd prime and put \(q=2\).
Among the prime ideals of \(L\) which are totally ramified over \(F\),
denote by \(\mathfrak{p}_1,\ldots,\mathfrak{p}_s\) those which split in \(N\),
\(\mathfrak{p}_i\mathcal{O}_N=\mathfrak{P}_i\mathfrak{P}_i^\prime\) for \(1\le i\le s\),
and by \(\mathfrak{q}_1,\ldots,\mathfrak{q}_n\) those which remain inert in \(N\),
\(\mathfrak{q}_j\mathcal{O}_N=\mathfrak{Q}_j\) for \(1\le j\le n\).
Then the space \(\mathcal{I}_{N/K}/\mathcal{I}_K\) of primitive ambiguous ideals of \(N/K\)
is the direct product of
the subspace \(\mathcal{I}_{L/F}/\mathcal{I}_F\) of \textbf{absolute differential factors} of \(L/F\)
and the subspace \(\ker(N_{N/L})\) of \textbf{relative differential factors} of \(N/K\),
whose bases and dimensions over \(\mathbb{F}_p\) are given by
\begin{equation}
\label{eqn:AbsAndRel}
\begin{aligned}
\mathcal{I}_{L/F}/\mathcal{I}_F
&\simeq \prod_{i=1}^{s}\,(\langle\mathfrak{p}_i\rangle/\langle\mathfrak{p}_i^p\rangle)
\times\prod_{j=1}^{n}\,(\langle\mathfrak{q}_j\rangle/\langle\mathfrak{q}_j^p\rangle)
\simeq\mathbb{F}_p^{s+n}, \quad
\dim_{\mathbb{F}_p}(\mathcal{I}_{L/F}/\mathcal{I}_F)=s+n, \\
\ker(N_{N/L})
&\simeq \prod_{i=1}^{s}\,\Bigl(\langle\mathfrak{P}_i(\mathfrak{P}_i^\prime)^{p-1}\rangle/\langle(\mathfrak{P}_i(\mathfrak{P}_i^\prime)^{p-1})^p\rangle\Bigr)
\simeq\mathbb{F}_p^s, \quad
\dim_{\mathbb{F}_p}(\ker(N_{N/L}))=s.
\end{aligned}
\end{equation}
Consequently, the complete space of differential factors has dimension \(\dim_{\mathbb{F}_p}(\mathcal{I}_{N/K}/\mathcal{I}_K)=n+2s\).
\end{theorem}

\begin{proof}
Whereas the qualitative formula
\eqref{eqn:Dichotomy}
is valid for any prime \(p\ge 2\) and any integer \(q>1\) with \(\gcd(p,q)=1\),
the quantitative description of the norm kernel \(\ker(N_{N/L})\) is only feasible
if we put \(q=2\) and therefore have to select an odd prime \(p\ge 3\).
Replacing \(N\) by \(L\) and \(K\) by \(F\) in formula
\eqref{eqn:PrimitiveAmbiguous},
we get \(t=n+s\) and thus the first isomorphism of formula
\eqref{eqn:AbsAndRel}.
For \(N\) and \(K\), however, we obtain \(t=n+2s\).
We point out that, if \(s=0\),
that is, if none of the totally ramified primes of \(L/F\) splits in \(N\),
then the induced norm mapping \(N_{N/L}\) in formula
\eqref{eqn:InducedNorm}
is an isomorphism.
For the constitution of the norm kernel, see
\cite[Thm. 3.4 and Cor. 3.3(3), p. 1994]{Ma2019a}.
\end{proof}

\newpage

\subsection{Primitive ambiguous principal ideals}
\label{ss:PrmAmbPrcIdl}
\noindent
The preceding result concerned \textit{primitive ambiguous} \textbf{ideals} of \(N/K\),
which can be interpreted as ideal factors of the \textit{relative different} \(\mathfrak{D}_{N/K}\).
Formula
\eqref{eqn:PrimitiveAmbiguous}
and Theorem
\ref{thm:QuantitativeDichotomy}
show that the \(\mathbb{F}_p\)-dimension of the space \(\mathcal{I}_{N/K}/\mathcal{I}_K\)
increases indefinitely with the number \(t\) of totally ramified prime ideals of \(N/K\).

\noindent
Now we restrict our attention to the space \(\mathcal{P}_{N/K}/\mathcal{P}_K\)
of \textit{primitive ambiguous} \textbf{principal ideals} or \textit{differential principal factors} (DPF) of \(N/K\).
We shall see that fundamental constraints from Galois cohomology
prohibit an infinite growth of its dimension over \(\mathbb{F}_p\),
for quadratic base fields \(K\).


\subsection{Splitting off the capitulation kernel}
\label{ss:CapitulationKernel}
\noindent
We have to cope with a difficulty
which arises in the case of a non-trivial class group
\(\mathrm{Cl}_K=\mathcal{I}_K/\mathcal{P}_K>1\),
because then \(\mathcal{P}_{N/K}/\mathcal{P}_K\) cannot be viewed as a subgroup of \(\mathcal{I}_{N/K}/\mathcal{I}_K\).
Therefore we must separate the \textit{capitulation kernel} of \(N/K\),
that is the kernel of the \textit{transfer} homomorphism
\(T_{N/K}:\,\mathrm{Cl}_K\to\mathrm{Cl}_N\), \(\mathfrak{a}\cdot\mathcal{P}_K\mapsto(\mathfrak{a}\mathcal{O}_N)\cdot\mathcal{P}_N\),
which extends classes of \(K\) to classes of \(N\):
\begin{equation}
\label{eqn:Capitulation}
\ker(T_{N/K})
=\lbrace\mathfrak{a}\cdot\mathcal{P}_K\mid(\exists\,A\in N)\,\mathfrak{a}\mathcal{O}_N=A\mathcal{O}_N\rbrace
=(\mathcal{I}_K\cap\mathcal{P}_N)/\mathcal{P}_K.
\end{equation}
On the one hand,
\(\ker(T_{N/K})=(\mathcal{I}_K\cap\mathcal{P}_N)/\mathcal{P}_K\) is a subgroup of \(\mathcal{I}_K/\mathcal{P}_K=\mathrm{Cl}_K\),
consisting of capitulating ideal classes of \(K\).
On the other hand,
since \(\mathcal{I}_K\le\mathcal{I}_{N/K}\) consists of ambiguous ideals of \(N/K\),
\(\ker(T_{N/K})=(\mathcal{I}_K\cap\mathcal{P}_N)/\mathcal{P}_K\) is a subgroup of \(\mathcal{P}_{N/K}/\mathcal{P}_K\),
consisting of special primitive ambiguous principal ideals of \(N/K\),
and we can form the quotient
\begin{equation}
\label{eqn:QuotientSeparation}
(\mathcal{P}_{N/K}/\mathcal{P}_K)/\bigl((\mathcal{I}_K\cap\mathcal{P}_N)/\mathcal{P}_K\bigr)
\simeq\mathcal{P}_{N/K}/(\mathcal{I}_K\cap\mathcal{P}_N)=\mathcal{P}_{N/K}/(\mathcal{I}_K\cap\mathcal{P}_{N/K})
\simeq(\mathcal{P}_{N/K}\cdot\mathcal{I}_K)/\mathcal{I}_K.
\end{equation}
This quotient relation of \(\mathbb{F}_p\)-vector spaces is equivalent to a direct product relation
\begin{equation}
\label{eqn:ProductSeparation}
\mathcal{P}_{N/K}/\mathcal{P}_K
\simeq(\mathcal{P}_{N/K}\cdot\mathcal{I}_K)/\mathcal{I}_K\times\ker(T_{N/K}).
\end{equation}
Since \((\mathcal{P}_{N/K}\cdot\mathcal{I}_K)/\mathcal{I}_K\le\mathcal{I}_{N/K}/\mathcal{I}_K\) is an actual inclusion,
the factorization of \(\mathcal{I}_{N/K}/\mathcal{I}_K\) in formula
\eqref{eqn:Dichotomy}
restricts to a factorization
\begin{equation}
\label{eqn:PrincipalDichotomy}
(\mathcal{P}_{N/K}\cdot\mathcal{I}_K)/\mathcal{I}_K
\simeq(\mathcal{P}_{L/F}/\mathcal{P}_F)\times\Bigl(\ker(N_{N/L})\cap\bigl((\mathcal{P}_{N/K}\cdot\mathcal{I}_K)/\mathcal{I}_K\bigr)\Bigr),
\end{equation}
provided that \(F\) is a field with trivial class group \(\mathrm{Cl}_F\),
that is \(\mathcal{I}_F=\mathcal{P}_F\)
and thus \(\mathcal{P}_{L/F}/\mathcal{P}_F\le\mathcal{I}_{L/F}/\mathcal{I}_F\).
Combining the formulas
\eqref{eqn:ProductSeparation}
and
\eqref{eqn:PrincipalDichotomy}
for the rational base field \(F=\mathbb{Q}\) ,
we obtain:


\begin{theorem}
\label{thm:Trichotomy}
There is a \textbf{trichotomic decomposition}
of the space \(\mathcal{P}_{N/K}/\mathcal{P}_K\) of differential principal factors of \(N/K\)
into three components,
\begin{equation}
\label{eqn:Trichotomy}
\mathcal{P}_{N/K}/\mathcal{P}_K\simeq
\mathcal{P}_{L/\mathbb{Q}}/\mathcal{P}_{\mathbb{Q}}
\times\Bigl(\ker(N_{N/L})\cap\bigl((\mathcal{P}_{N/K}\mathcal{I}_K)/\mathcal{I}_K\bigr)\Bigr)
\times\ker(T_{N/K}),
\end{equation}
\(\bullet\) the \textbf{absolute principal factors}, \(\mathcal{P}_{L/\mathbb{Q}}/\mathcal{P}_{\mathbb{Q}}\), of \(L/\mathbb{Q}\), \\
\(\bullet\) the \textbf{relative principal factors}, \(\ker(N_{N/L})\cap\bigl((\mathcal{P}_{N/K}\mathcal{I}_K)/\mathcal{I}_K\bigr)\), of \(N/K\), and \\
\(\bullet\) the \textbf{capitulation kernel}, \(\ker(T_{N/K})\), of \(N/K\).
\end{theorem}


\subsection{Galois cohomology}
\label{ss:GaloisCohomology}
\noindent
In order to establish a quantitative version of the qualitative formula
\eqref{eqn:Trichotomy},
we suppose that \(N/K\) is a cyclic relative extension of odd prime degree \(p\)
and we use the Galois cohomology of the unit group \(U_N\)
as a module over the automorphism group \(G=\mathrm{Gal}(N/K)=\langle\sigma\rangle\simeq C_p\).
In fact, we combine a theorem of Iwasawa
\cite{Iw1956}
on the first cohomology \(\mathrm{H}^1(G,U_N)\)
with a theorem of Hasse
\cite{Ha1927}
on the Herbrand quotient of \(U_N\)
\cite{Hb1932},
and we use Dirichlet's theorem on the torsion-free unit rank of \(K\).
By \(E_{N/K}=U_N\cap\ker(N_{N/K})\) we denote the group of \textit{relative units} of \(N/K\).
\begin{equation}
\label{eqn:Herbrand}
\begin{aligned}
\mathrm{H}^1(G,U_N) &\simeq E_{N/K}/U_N^{\sigma-1}\simeq\mathcal{P}_{N/K}/\mathcal{P}_K\ \text{(Iwasawa)}, \\
\hat{\mathrm{H}}^0(G,U_N) &\simeq U_K/N_{N/K}(U_N), \text{ with } (U_K:N_{N/K}(U_N))=p^U,\ 0\le U\le r_1+r_2-\theta, \\
\frac{\#\mathrm{H}^1(G,U_N)}{\#\hat{\mathrm{H}}^0(G,U_N)} &= \lbrack N:K\rbrack=p \quad \text{(Hasse)},
\end{aligned}
\end{equation}
where \((r_1,r_2)\) is the signature of \(K\), and \(\theta=0\) if \(K\) contains the \(p\)th roots of unity, but \(\theta=1\) else.


\begin{corollary}
\label{cor:Trichotomy}
If \(N/K\) is cyclic of odd prime degree \(p\ge 3\),
then the \(\mathbb{F}_p\)-dimensions of the spaces of differential principal factors in Theorem
\ref{thm:Trichotomy}
are connected by the \textbf{fundamental equation}
\begin{equation}
\label{eqn:Dimensions}
U+1=A+R+C,\quad \text{where}
\end{equation}
\(\bullet\) \(A:=\dim_{\mathbb{F}_p}(\mathcal{P}_{L/\mathbb{Q}}/\mathcal{P}_{\mathbb{Q}})\) is the dimension of \textbf{absolute} principal factors, \\
\(\bullet\) \(R:=\dim_{\mathbb{F}_p}\Bigl(\ker(N_{N/L})\cap\bigl((\mathcal{P}_{N/K}\mathcal{I}_K)/\mathcal{I}_K\bigr)\Bigr)\)
is the dimension of \textbf{relative} principal factors, and \\
\(\bullet\) \(C:=\dim_{\mathbb{F}_p}(\ker(T_{N/K}))\) is the dimension of the \textbf{capitulation} kernel.
\end{corollary}


\begin{corollary}
\label{cor:Estimates}
Under the assumptions \(p\ge 3\), \(q=2\) of Theorem
\ref{thm:QualitativeDichotomy},
in particular for \(N\) dihedral of degree \(2p\),
the dimensions in Corollary
\ref{cor:Trichotomy}
are bounded by the following \textbf{fundamental estimates}
\begin{equation}
\label{eqn:Estimates}
0\le A\le\min(n+s,m), \quad
0\le R\le\min(s,m), \quad
0\le C\le\min(\varrho_p,m),
\end{equation}
where \(\varrho_p:=\mathrm{rank}_p(\mathrm{Cl}_K)\),
and \(m:=1+r_1+r_2-\theta\) denotes the cohomological maximum of \(U+1\).
In particular, we have \\
\(m=2\) for real quadratic \(K\) with \((r_1,r_2)=(2,0)\) or \(K=\mathbb{Q}(\sqrt{-3})\) if \(p=3\), \\
\(m=1\) for imaginary quadratic \(K\) with \((r_1,r_2)=(0,1)\), except for \(K=\mathbb{Q}(\sqrt{-3})\) when \(p=3\).
\end{corollary}


\begin{remark}
\label{rmk:Estimates}
For \(N\) pure metacyclic of degree \((p-1)p\),
the space \(\mathcal{P}_{L/\mathbb{Q}}/\mathcal{P}_{\mathbb{Q}}\) of absolute principal factors
contains the one-dimensional subspace \(\Delta=\langle\sqrt[p]{D}\rangle\) generated by the \textit{radicals}, and thus
\begin{equation}
\label{eqn:EstimatesMet}
1\le A\le\min(t,m), \
0\le R\le m-1, \
0\le C\le\min(\varrho_p,m-1),
\end{equation}
where \(m=\frac{p+1}{2}\) for cyclotomic \(K\) with \((r_1,r_2)=(0,\frac{p-1}{2})\).
In particular, there is no capitulation, \(C=0\), for a \textit{regular} prime \(p\) with \(\varrho_p=0\), for instance \(p<37\).
\end{remark}


\begin{remark}
\label{rmk:Inclusion}
We mentioned that in general
\(\mathcal{P}_{N/K}/\mathcal{P}_K\) cannot be viewed as a subgroup of \(\mathcal{I}_{N/K}/\mathcal{I}_K\).
In fact, for a dihedral field \(N\) which is unramified with conductor \(f=1\) over \(K\),
we have \(n=s=0\),
consequently \(A=R=0\),
and \(\mathcal{I}_{N/K}/\mathcal{I}_K\simeq 0\) is the nullspace,
whereas \(\mathcal{P}_{N/K}/\mathcal{P}_K\simeq\ker(T_{N/K})\) is at least one-dimensional,
according to Hilbert's Theorem 94
\cite{Hi1897},
and at most two-dimensional by the estimate \(C\le\min(\varrho_p,m)\le\min(\varrho_p,2)\le 2\).
\end{remark}


In the next two sections, we apply the results of \S\S\
\ref{ss:PrmAmbIdl}
--
\ref{ss:GaloisCohomology}
to various extensions \(N/K\).


\subsection{Differential principal factorization (DPF) types of complex dihedral fields}
\label{ss:ComplexDihedralTypes}
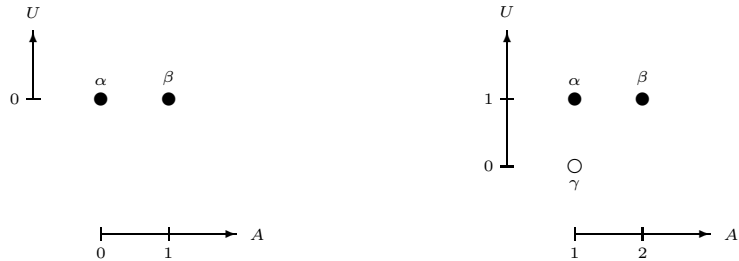
\begin{figure}[ht]
\caption{Classification of simply real dihedral and pure cubic fields}
\label{fig:MoserExtendedCubic}

{\tiny

\setlength{\unitlength}{0.9cm}
\begin{picture}(15,3.5)(-11,-9)




\put(-9,-5.8){\makebox(0,0)[cb]{\(U\)}}
\put(-9,-7){\vector(0,1){1}}
\put(-9,-7){\line(0,1){0}}
\multiput(-9.1,-7)(0,1){1}{\line(1,0){0.2}}

\put(-9.2,-7){\makebox(0,0)[rc]{\(0\)}}

\put(-5.8,-9){\makebox(0,0)[lc]{\(A\)}}
\put(-7,-9){\vector(1,0){1}}
\put(-8,-9){\line(1,0){1}}
\multiput(-8,-9.1)(1,0){2}{\line(0,1){0.2}}

\put(-8,-9.2){\makebox(0,0)[ct]{\(0\)}}
\put(-7,-9.2){\makebox(0,0)[ct]{\(1\)}}


\put(-8,-7){\circle*{0.2}}
\put(-8,-6.8){\makebox(0,0)[cb]{\(\alpha\)}}
\put(-7,-7){\circle*{0.2}}
\put(-7,-6.8){\makebox(0,0)[cb]{\(\beta\)}}




\put(-2,-5.8){\makebox(0,0)[cb]{\(U\)}}
\put(-2,-7){\vector(0,1){1}}
\put(-2,-8){\line(0,1){1}}
\multiput(-2.1,-8)(0,1){2}{\line(1,0){0.2}}

\put(-2.2,-7){\makebox(0,0)[rc]{\(1\)}}
\put(-2.2,-8){\makebox(0,0)[rc]{\(0\)}}

\put(1.2,-9){\makebox(0,0)[lc]{\(A\)}}
\put(0,-9){\vector(1,0){1}}
\put(-1,-9){\line(1,0){1}}
\multiput(-1,-9.1)(1,0){2}{\line(0,1){0.2}}

\put(-1,-9.2){\makebox(0,0)[ct]{\(1\)}}
\put(0,-9.2){\makebox(0,0)[ct]{\(2\)}}


\put(-1,-7){\circle*{0.2}}
\put(-1,-6.8){\makebox(0,0)[cb]{\(\alpha\)}}
\put(0,-7){\circle*{0.2}}
\put(0,-6.8){\makebox(0,0)[cb]{\(\beta\)}}

\put(-1,-8){\circle{0.2}}
\put(-1,-8.2){\makebox(0,0)[ct]{\(\gamma\)}}


\end{picture}
}
\end{figure}

\noindent
Let \(p\) be an odd prime.
We recall the classification theorem
for \textit{pure cubic} fields \(L=\mathbb{Q}(\sqrt[3]{D})\)
and their Galois closure \(N=\mathbb{Q}(\zeta_3,\sqrt[3]{D})\),
that is the metacyclic case \(p=3\).
The \textit{coarse} classification of \(N\)
according to the cohomological invariants \(U\) and \(A\) alone
is closely related to the
classification of \textit{simply real dihedral} fields of degree \(2p\) with any odd prime \(p\)
by Nicole Moser
\cite[Dfn. III.1 and Prop. III.3, p. 61]{Mo1979},
as illustrated in Figure
\ref{fig:MoserExtendedCubic}.
The coarse types \(\alpha\) and \(\beta\)
are completely analogous in both cases.
The additional type \(\gamma\) is required for pure cubic fields,
because there arises the possibility that the primitive cube root of unity \(\zeta_3\)
occurs as relative norm \(N_{N/K}(Z)\) of a unit \(Z\in U_N\).
Due to the existence of radicals in the pure cubic case,
the \(\mathbb{F}_p\)-dimension \(A\) of the vector space of absolute DPF
exceeds the corresponding dimension for simply real dihedral fields by one.


\noindent
The \textit{fine} classification of \(N\)
according to the invariants \(U\), \(A\), \(R\) and \(C\)
in the simply real dihedral situation with \(U+1=A+R+C\)
splits type \(\alpha\) with \(A=0\) further in
type \(\alpha_1\) with \(C=1\) (capitulation) and
type \(\alpha_2\) with \(R=1\) (relative DPF).
In the pure cubic situation, however, no further splitting occurs,
since \(C=0\), and \(R=U+1-A\) is determined uniquely by \(U\) and \(A\) already.
We oppose the two classifications in the following theorems.

\newpage

\begin{theorem}
\label{thm:MainComplex}
Each \textbf{simply real dihedral} field \(N/\mathbb{Q}\)
of absolute degree \(\lbrack N:\mathbb{Q}\rbrack=2p\) with an odd prime \(p\)
belongs to precisely one of the following \(3\) differential principal factorization types,
in dependence on the triplet \((A,R,C)\):


\renewcommand{\arraystretch}{1.1}

\begin{table}[ht]
\label{tbl:ComplexDPFTypes}
\begin{center}
\begin{tabular}{|c||c||c||ccc|}
\hline
 Type        & \(U\) & \(U+1=A+R+C\) & \(A\) & \(R\) & \(C\) \\
\hline
\(\alpha_1\) & \(0\) & \(1\) & \(0\) & \(0\) & \(1\) \\
\(\alpha_2\) & \(0\) & \(1\) & \(0\) & \(1\) & \(0\) \\
\(\beta\)    & \(0\) & \(1\) & \(1\) & \(0\) & \(0\) \\
\hline
\end{tabular}
\end{center}
\end{table}


\end{theorem}

\begin{proof}
Consequence of Cor.
\ref{cor:Trichotomy}
and 
\ref{cor:Estimates}.
See \cite[Dfn. III.1 and Prop. III.3, p. 61]{Mo1979} and \cite{Ma1991b}.
\end{proof}

\begin{theorem}
\label{thm:MainCubic}
Each \textbf{pure metacyclic} field \(N=\mathbb{Q}(\zeta_3,\sqrt[3]{D})\)
of absolute degree \(\lbrack N:\mathbb{Q}\rbrack=6\)
with cube free radicand \(D\in\mathbb{Z}\), \(D\ge 2\),
belongs to precisely one of the following \(3\) differential principal factorization types,
in dependence on the invariant \(U\) and the pair \((A,R)\):


\renewcommand{\arraystretch}{1.1}

\begin{table}[ht]
\label{tbl:CubicDPFTypes}
\begin{center}
\begin{tabular}{|c||c||c||cc|}
\hline
 Type      & \(U\) & \(U+1=A+R\) & \(A\) & \(R\) \\
\hline
\(\alpha\) & \(1\) & \(2\) & \(1\) & \(1\) \\
\(\beta\)  & \(1\) & \(2\) & \(2\) & \(0\) \\
\hline
\(\gamma\) & \(0\) & \(1\) & \(1\) & \(0\) \\
\hline
\end{tabular}
\end{center}
\end{table}


\end{theorem}

\begin{proof}
A part of the proof is due to Barrucand and Cohn
\cite{BaCo1971}
who distinguished \(4\) different types,
\(\mathrm{I}\hat{=}\beta\), \(\mathrm{II}\), \(\mathrm{III}\hat{=}\alpha\), and \(\mathrm{IV}\hat{=}\gamma\).
However, Halter-Koch
\cite{HK1976}
showed the impossibility of one of these types, namely type \(\mathrm{II}\).
Our new proof with different methods is given in
\cite[Thm. 2.1, p. 254]{AMITA2020}.
\end{proof}


\subsection{Differential principal factorization (DPF) types of real dihedral fields}
\label{ss:RealDihedralAndQuinticTypes}
\begin{figure}[ht]
\caption{Classification of totally real dihedral and pure quintic fields}
\label{fig:MoserExtendedQuintic}

{\tiny

\setlength{\unitlength}{1.0cm}
\begin{picture}(15,5)(-11,-10)




\put(-9,-5.8){\makebox(0,0)[cb]{\(U\)}}
\put(-9,-7){\vector(0,1){1}}
\put(-9,-8){\line(0,1){1}}
\multiput(-9.1,-8)(0,1){2}{\line(1,0){0.2}}

\put(-9.2,-7){\makebox(0,0)[rc]{\(1\)}}
\put(-9.2,-8){\makebox(0,0)[rc]{\(0\)}}

\put(-4.8,-10){\makebox(0,0)[lc]{\(A\)}}
\put(-6,-10){\vector(1,0){1}}
\put(-8,-10){\line(1,0){2}}
\multiput(-8,-10.1)(1,0){3}{\line(0,1){0.2}}

\put(-8,-10.2){\makebox(0,0)[ct]{\(0\)}}
\put(-7,-10.2){\makebox(0,0)[ct]{\(1\)}}
\put(-6,-10.2){\makebox(0,0)[ct]{\(2\)}}


\put(-8,-7){\circle*{0.2}}
\put(-8,-6.8){\makebox(0,0)[cb]{\(\alpha\)}}
\put(-7,-7){\circle*{0.2}}
\put(-7,-6.8){\makebox(0,0)[cb]{\(\beta\)}}
\put(-6,-7){\circle*{0.2}}
\put(-6,-6.8){\makebox(0,0)[cb]{\(\gamma\)}}

\put(-8,-8){\circle*{0.2}}
\put(-8,-7.8){\makebox(0,0)[cb]{\(\delta\)}}
\put(-7,-8){\circle*{0.2}}
\put(-7,-7.8){\makebox(0,0)[cb]{\(\varepsilon\)}}




\put(-2,-5.8){\makebox(0,0)[cb]{\(U\)}}
\put(-2,-7){\vector(0,1){1}}
\put(-2,-9){\line(0,1){2}}
\multiput(-2.1,-9)(0,1){3}{\line(1,0){0.2}}

\put(-2.2,-7){\makebox(0,0)[rc]{\(2\)}}
\put(-2.2,-8){\makebox(0,0)[rc]{\(1\)}}
\put(-2.2,-9){\makebox(0,0)[rc]{\(0\)}}

\put(2.2,-10){\makebox(0,0)[lc]{\(A\)}}
\put(1,-10){\vector(1,0){1}}
\put(-1,-10){\line(1,0){2}}
\multiput(-1,-10.1)(1,0){3}{\line(0,1){0.2}}

\put(-1,-10.2){\makebox(0,0)[ct]{\(1\)}}
\put(0,-10.2){\makebox(0,0)[ct]{\(2\)}}
\put(1,-10.2){\makebox(0,0)[ct]{\(3\)}}


\put(-1,-7){\circle*{0.2}}
\put(-1,-6.8){\makebox(0,0)[cb]{\(\alpha\)}}
\put(0,-7){\circle*{0.2}}
\put(0,-6.8){\makebox(0,0)[cb]{\(\beta\)}}
\put(1,-7){\circle*{0.2}}
\put(1,-6.8){\makebox(0,0)[cb]{\(\gamma\)}}

\put(-1,-7.9){\circle*{0.2}}
\put(-1,-7.7){\makebox(0,0)[cb]{\(\delta\)}}
\put(0,-7.9){\circle*{0.2}}
\put(0,-7.7){\makebox(0,0)[cb]{\(\varepsilon\)}}
\put(-1,-8.1){\circle{0.2}}
\put(-1,-8.3){\makebox(0,0)[ct]{\(\zeta\)}}
\put(0,-8.1){\circle{0.2}}
\put(0,-8.3){\makebox(0,0)[ct]{\(\eta\)}}

\put(-1,-9){\circle{0.2}}
\put(-1,-9.2){\makebox(0,0)[ct]{\(\vartheta\)}}


\end{picture}
}
\end{figure}
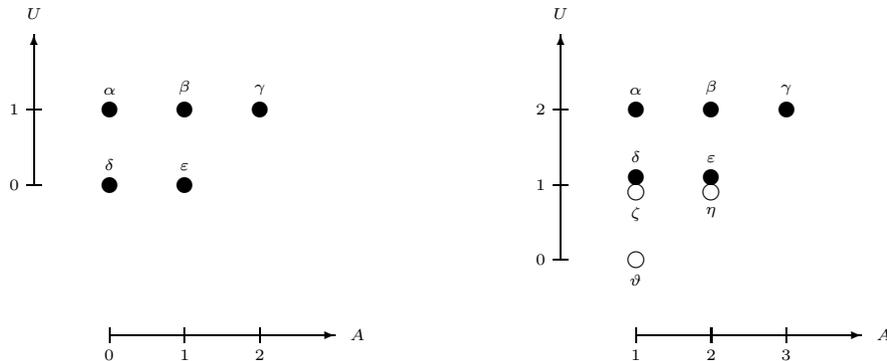

\noindent
Now we state the classification theorem
for \textit{pure quintic} fields \(L=\mathbb{Q}(\sqrt[5]{D})\)
and their Galois closure \(N=\mathbb{Q}(\zeta_5,\sqrt[5]{D})\),
that is the metacyclic case \(p=5\).
The \textit{coarse} classification of \(N\)
according to the invariants \(U\) and \(A\) alone
is closely related to the
classification of \textit{totally real dihedral} fields of degree \(2p\) with any odd prime \(p\)
by Nicole Moser
\cite[Thm. III.5, p. 62]{Mo1979},
as illustrated in Figure
\ref{fig:MoserExtendedQuintic}.
The coarse types \(\alpha\), \(\beta\), \(\gamma\), \(\delta\), \(\varepsilon\)
are completely analogous in both cases.
Additional types \(\zeta\), \(\eta\), \(\vartheta\) are required for pure quintic fields,
because there arises the possibility that the primitive fifth root of unity \(\zeta_5\)
occurs as relative norm \(N_{N/K}(Z)\) of a unit \(Z\in U_N\).
Due to the existence of radicals in the pure quintic case,
the \(\mathbb{F}_p\)-dimension \(A\) of the vector space of absolute DPF
exceeds the corresponding dimension for totally real dihedral fields by one
(see Remark
\ref{rmk:Estimates}).


\noindent
The \textit{fine} classification of \(N\) according to the invariants \(U\), \(A\), \(R\) and \(C\)
in the totally real dihedral situation with \(U+1=A+R+C\)
splits type \(\alpha\) with \(U=1\), \(A=0\) further in
type \(\alpha_1\) with \(C=2\) (double capitulation),
type \(\alpha_2\) with \(C=R=1\) (mixed capitulation and relative DPF),
type \(\alpha_3\) with \(R=2\) (double relative DPF),
type \(\beta\) with \(U=A=1\) in
type \(\beta_1\) with \(C=1\) (capitulation),
type \(\beta_2\) with \(R=1\) (relative DPF), and
type \(\delta\) with \(U=A=0\) in
type \(\delta_1\) with \(C=1\) (capitulation),
type \(\delta_2\) with \(R=1\) (relative DPF).
In the pure quintic situation with \(U+1=A+I+R\)
\cite{Ma2019a},
however, we arrive at the second of the following theorems
where we oppose the two classifications.

\begin{theorem}
\label{thm:MainReal}
Each \textbf{totally real dihedral} field \(N/\mathbb{Q}\)
of absolute degree \(\lbrack N:\mathbb{Q}\rbrack=2p\) with an odd prime \(p\)
belongs to precisely one of the following \(9\) differential principal factorization types,
in dependence on the invariant \(U\) and the triplet \((A,R,C)\).


\renewcommand{\arraystretch}{1.1}

\begin{table}[ht]
\label{tbl:RealDPFTypes}
\begin{center}
\begin{tabular}{|c||c||c||ccc|}
\hline
 Type           & \(U\) & \(U+1=A+R+C\) & \(A\) & \(R\) & \(C\) \\
\hline
\(\alpha_1\)    & \(1\) & \(2\) & \(0\) & \(0\) & \(2\) \\
\(\alpha_2\)    & \(1\) & \(2\) & \(0\) & \(1\) & \(1\) \\
\(\alpha_3\)    & \(1\) & \(2\) & \(0\) & \(2\) & \(0\) \\
\(\beta_1\)     & \(1\) & \(2\) & \(1\) & \(0\) & \(1\) \\
\(\beta_2\)     & \(1\) & \(2\) & \(1\) & \(1\) & \(0\) \\
\(\gamma\)      & \(1\) & \(2\) & \(2\) & \(0\) & \(0\) \\
\hline
\(\delta_1\)    & \(0\) & \(1\) & \(0\) & \(0\) & \(1\) \\
\(\delta_2\)    & \(0\) & \(1\) & \(0\) & \(1\) & \(0\) \\
\(\varepsilon\) & \(0\) & \(1\) & \(1\) & \(0\) & \(0\) \\
\hline
\end{tabular}
\end{center}
\end{table}


\end{theorem}

\begin{proof}
Consequence of the Corollaries
\ref{cor:Trichotomy}
and 
\ref{cor:Estimates}.
See also \cite[Thm. III.5, p. 62]{Mo1979} and \cite{Ma1991b}.
\end{proof}


\begin{theorem}
\label{thm:MainQuintic}
Each \textbf{pure metacyclic} field \(N=\mathbb{Q}(\zeta_5,\sqrt[5]{D})\)
of absolute degree \(\lbrack N:\mathbb{Q}\rbrack=20\)
with \(5\)-th power free radicand \(D\in\mathbb{Z}\), \(D\ge 2\),
belongs to precisely one of the following \(13\) differential principal factorization types,
in dependence on the invariant \(U\) and the triplet \((A,I,R)\).


\renewcommand{\arraystretch}{1.1}

\begin{table}[ht]
\label{tbl:QuinticDPFTypes}
\begin{center}
\begin{tabular}{|c||c||c||ccc|}
\hline
 Type           & \(U\) & \(U+1=A+I+R\) & \(A\) & \(I\) & \(R\) \\
\hline
\(\alpha_1\)    & \(2\) & \(3\) & \(1\) & \(0\) & \(2\) \\
\(\alpha_2\)    & \(2\) & \(3\) & \(1\) & \(1\) & \(1\) \\
\(\alpha_3\)    & \(2\) & \(3\) & \(1\) & \(2\) & \(0\) \\
\(\beta_1\)     & \(2\) & \(3\) & \(2\) & \(0\) & \(1\) \\
\(\beta_2\)     & \(2\) & \(3\) & \(2\) & \(1\) & \(0\) \\
\(\gamma\)      & \(2\) & \(3\) & \(3\) & \(0\) & \(0\) \\
\hline
\(\delta_1\)    & \(1\) & \(2\) & \(1\) & \(0\) & \(1\) \\
\(\delta_2\)    & \(1\) & \(2\) & \(1\) & \(1\) & \(0\) \\
\(\varepsilon\) & \(1\) & \(2\) & \(2\) & \(0\) & \(0\) \\
 \hline
\(\zeta_1\)     & \(1\) & \(2\) & \(1\) & \(0\) & \(1\) \\
\(\zeta_2\)     & \(1\) & \(2\) & \(1\) & \(1\) & \(0\) \\
\(\eta\)        & \(1\) & \(2\) & \(2\) & \(0\) & \(0\) \\
 \hline
\(\vartheta\)   & \(0\) & \(1\) & \(1\) & \(0\) & \(0\) \\
\hline
\end{tabular}
\end{center}
\end{table}


The types \(\delta_1\), \(\delta_2\), \(\varepsilon\)
are characterized additionally by \(\zeta_5\not\in N_{N/K}(U_N)\),
and the types \(\zeta_1\), \(\zeta_2\), \(\eta\)
by \(\zeta_5\in N_{N/K}(U_N)\).
\end{theorem}

\begin{proof}
The proof is given in
\cite[Thm. 6.1]{Ma2019a}.
\end{proof}


\begin{remark}
\label{rmk:RealDPFTypes}
Our classification of totally real dihedral fields in Theorem
\ref{thm:MainReal}
refines the classification by Moser
\cite{Mo1979}
who uses the results on integral representations
of the dihedral group \(D_p\) by Lee
\cite{Le1964}.
She denotes by \(b=(U_K:N_{N/K}(U_N))\) the unit norm index
and obtains \(U_N=U_K\cdot E_{N/K}\) as a \textit{split extension} (direct product)
of \(U_K\) by \(E_{N/K}\) for \(b=p\) (types \(\alpha,\beta,\gamma\)),
and \((U_N:U_K\cdot E_{N/K})=p\) as a \textit{non-split extension}
(of modules over \(\mathbb{Z}\lbrack D_p\rbrack\))
for \(b=1\) (types \(\delta,\varepsilon\)),
due to a non-trivial relation \(N_{N/K}(H)=H^{1+\sigma+\ldots+\sigma^{p-1}}=\eta\)
for the fundamental unit \(\eta\) of \(K\) and a unit \(H\in U_N\setminus(U_K\cdot E_{N/K})\).

For \(p=3\), a geometric interpretation of the
unit lattice in logarithmic space,
i.e., the Dirichlet-Minkowski image of \(U_N\),
has been given by Hasse
\cite{Ha1948,Ha1950}.
\end{remark}

\newpage

\section{Classifying Angell's range \(0<d_L<100\,000\)}
\label{s:Angell}

\noindent
In Table
\ref{tbl:Angell}
and all the following tables,
we present the results of our classification
of totally real cubic fields \(L\)
and their normal closures \(N\)
into \textit{differential principal factorization types} \(\tau(L)=\tau(N)\).
The rows correspond to the numerous steps
where we applied our algorithm (\S\
\ref{ss:Construction})
to various configurations of
\(3\)-class rank \(\varrho_3\) of the real quadratic subfield \(K\) of \(N\)
and \(3\)-admissible conductors \(f\) of \(N/K\).
Here, \(d\) denotes the fundamental discriminant of \(K\),
and \(q,q_1,q_2\), resp. \(\ell,\ell_1,\ell_2\),
denote prime numbers congruent to \(2\), resp. \(1\) modulo \(3\).
In Table
\ref{tbl:Angell},
the types \(\alpha_2\) and \(\alpha_3\) do not yet occur.


\renewcommand{\arraystretch}{1.1}

\begin{table}[ht]
\caption{Totally real cubic discriminants \(d_L=f^2\cdot d\) in the range \(0<d_L<10^5\)}
\label{tbl:Angell}
\begin{center}
\begin{tabular}{|rl||r|rrrr||rrrr|rrr|r|}
\hline
       &           & \multicolumn{5}{c||}{Multiplicity}    & \multicolumn{7}{c|}{Differential Principal Factorization} &       \\
 \(f\) & Condition & \(0\) & \(1\) & \(2\) & \(3\) & \(4\) & \(\alpha_1\) & \(\beta_1\) & \(\beta_2\) & \(\gamma\) & \(\delta_1\) & \(\delta_2\) & \(\varepsilon\) & Total \\
\hline
 \(1\)      &                   \(\varrho_3=0\) &\(27089\) &    \(\) &  \(\) &   \(\) & \(\) & \(\) &  \(\) &   \(\) &   \(\) &    \(\)&   \(\) &    \(\) &    \(0\) \\
\hline
 \(q\)      &   \(\equiv 2\,(\mathrm{mod}\,3)\) & \(2219\) & \(806\) &  \(\) &   \(\) & \(\) & \(\) &  \(\) &   \(\) &   \(\) &   \(\) &   \(\) & \(806\) &  \(806\) \\
 \(3\)      &  \(d\equiv 3\,(\mathrm{mod}\,9)\) &  \(287\) & \(109\) &  \(\) &   \(\) & \(\) & \(\) &  \(\) &   \(\) &   \(\) &   \(\) &   \(\) & \(109\) &  \(109\) \\
 \(3\)      &  \(d\equiv 6\,(\mathrm{mod}\,9)\) &  \(284\) & \(105\) &  \(\) &   \(\) & \(\) & \(\) &  \(\) &   \(\) &   \(\) &   \(\) &   \(\) & \(105\) &  \(105\) \\
 \(9\)      &  \(d\equiv 6\,(\mathrm{mod}\,9)\) &    \(9\) &  \(38\) &  \(\) &  \(1\) & \(\) & \(\) &  \(\) &   \(\) &   \(\) &   \(\) &   \(\) &  \(41\) &   \(41\) \\
 \(9\)      &  \(d\equiv 2\,(\mathrm{mod}\,3)\) &  \(102\) &  \(34\) &  \(\) &   \(\) & \(\) & \(\) &  \(\) &   \(\) &   \(\) &   \(\) &   \(\) &  \(34\) &   \(34\) \\
\hline
 \(9\)      &  \(d\equiv 1\,(\mathrm{mod}\,3)\) &   \(96\) &  \(31\) &  \(\) &   \(\) & \(\) & \(\) &  \(\) &  \(8\) &   \(\) &   \(\) & \(20\) &   \(3\) &   \(31\) \\
 \(\ell\)   &   \(\equiv 1\,(\mathrm{mod}\,3)\) &  \(316\) &  \(86\) &  \(\) &   \(\) & \(\) & \(\) &  \(\) & \(20\) &   \(\) &   \(\) & \(59\) &   \(7\) &   \(86\) \\
\hline
 \(q_1q_2\) &                                   &   \(30\) &  \(38\) & \(2\) &   \(\) & \(\) & \(\) &  \(\) &   \(\) & \(38\) &   \(\) &   \(\) &   \(4\) &   \(42\) \\
 \(3q\)     &  \(d\equiv 3\,(\mathrm{mod}\,9)\) &   \(23\) &  \(23\) &  \(\) &   \(\) & \(\) & \(\) &  \(\) &   \(\) & \(23\) &   \(\) &   \(\) &    \(\) &   \(23\) \\
 \(3q\)     &  \(d\equiv 6\,(\mathrm{mod}\,9)\) &   \(19\) &  \(25\) & \(1\) &   \(\) & \(\) & \(\) &  \(\) &   \(\) & \(25\) &   \(\) &   \(\) &   \(2\) &   \(27\) \\
 \(9q\)     &  \(d\equiv 6\,(\mathrm{mod}\,9)\) &     \(\) &    \(\) & \(4\) &  \(1\) & \(\) & \(\) &  \(\) &   \(\) &  \(9\) &   \(\) &   \(\) &   \(2\) &   \(11\) \\
 \(9q\)     &  \(d\equiv 2\,(\mathrm{mod}\,3)\) &    \(6\) &   \(8\) &  \(\) &   \(\) & \(\) & \(\) &  \(\) &   \(\) &  \(8\) &   \(\) &   \(\) &    \(\) &    \(8\) \\
\hline
 \(9q\)     &  \(d\equiv 1\,(\mathrm{mod}\,3)\) &    \(5\) &  \(10\) &  \(\) &   \(\) & \(\) & \(\) &  \(\) & \(10\) &   \(\) &   \(\) &   \(\) &    \(\) &   \(10\) \\
 \(q\ell\)  &                                   &   \(13\) &  \(29\) & \(1\) &   \(\) & \(\) & \(\) &  \(\) & \(29\) &   \(\) &   \(\) &   \(\) &   \(2\) &   \(31\) \\
 \(3\ell\)  &  \(d\equiv 3\,(\mathrm{mod}\,9)\) &    \(1\) &   \(5\) &  \(\) &   \(\) & \(\) & \(\) &  \(\) &  \(5\) &   \(\) &   \(\) &   \(\) &    \(\) &    \(5\) \\
 \(3\ell\)  &  \(d\equiv 6\,(\mathrm{mod}\,9)\) &    \(2\) &   \(3\) &  \(\) &   \(\) & \(\) & \(\) &  \(\) &  \(3\) &   \(\) &   \(\) &   \(\) &    \(\) &    \(3\) \\
 \(9\ell\)  &  \(d\equiv 2\,(\mathrm{mod}\,3)\) &     \(\) &   \(1\) &  \(\) &   \(\) & \(\) & \(\) &  \(\) &  \(1\) &   \(\) &   \(\) &   \(\) &    \(\) &    \(1\) \\
\hline
\(3q_1q_2\) &  \(d\equiv 3\,(\mathrm{mod}\,9)\) &     \(\) &   \(1\) & \(1\) &   \(\) & \(\) & \(\) &  \(\) &   \(\) &  \(3\) &   \(\) &   \(\) &    \(\) &    \(3\) \\
\hline
\hline
 \(1\)      &                   \(\varrho_3=1\) &     \(\) &\(3300\) &  \(\) &   \(\) & \(\) & \(\) &  \(\) &   \(\) &   \(\) &\(3300\)&   \(\) &    \(\) & \(3300\) \\
\hline
 \(q\)      &   \(\equiv 2\,(\mathrm{mod}\,3)\) &  \(261\) &    \(\) &  \(\) & \(14\) & \(\) & \(\) & \(4\) &   \(\) &   \(\) & \(36\) &   \(\) &   \(2\) &   \(42\) \\
 \(3\)      &  \(d\equiv 3\,(\mathrm{mod}\,9)\) &   \(27\) &    \(\) &  \(\) &   \(\) & \(\) & \(\) &  \(\) &   \(\) &   \(\) &   \(\) &   \(\) &    \(\) &    \(0\) \\
 \(3\)      &  \(d\equiv 6\,(\mathrm{mod}\,9)\) &   \(34\) &    \(\) &  \(\) &  \(1\) & \(\) & \(\) &  \(\) &   \(\) &   \(\) &  \(3\) &   \(\) &    \(\) &    \(3\) \\
 \(9\)      &  \(d\equiv 6\,(\mathrm{mod}\,9)\) &    \(1\) &    \(\) &  \(\) &   \(\) & \(\) & \(\) &  \(\) &   \(\) &   \(\) &   \(\) &   \(\) &    \(\) &    \(0\) \\
 \(9\)      &  \(d\equiv 2\,(\mathrm{mod}\,3)\) &    \(6\) &    \(\) &  \(\) &   \(\) & \(\) & \(\) &  \(\) &   \(\) &   \(\) &   \(\) &   \(\) &    \(\) &    \(0\) \\
\hline
 \(9\)      &  \(d\equiv 1\,(\mathrm{mod}\,3)\) &   \(10\) &    \(\) &  \(\) &   \(\) & \(\) & \(\) &  \(\) &   \(\) &   \(\) &   \(\) &   \(\) &    \(\) &    \(0\) \\
 \(\ell\)   &   \(\equiv 1\,(\mathrm{mod}\,3)\) &   \(25\) &    \(\) &  \(\) &  \(3\) & \(\) & \(\) & \(3\) &   \(\) &   \(\) &  \(6\) &   \(\) &    \(\) &    \(9\) \\
\hline
 \(q_1q_2\) &                                   &    \(1\) &    \(\) &  \(\) &   \(\) & \(\) & \(\) &  \(\) &   \(\) &   \(\) &   \(\) &   \(\) &    \(\) &    \(0\) \\
 \(3q\)     &  \(d\equiv 3\,(\mathrm{mod}\,9)\) &    \(1\) &    \(\) &  \(\) &   \(\) & \(\) & \(\) &  \(\) &   \(\) &   \(\) &   \(\) &   \(\) &    \(\) &    \(0\) \\
 \(3q\)     &  \(d\equiv 6\,(\mathrm{mod}\,9)\) &    \(1\) &    \(\) &  \(\) &  \(1\) & \(\) & \(\) & \(3\) &   \(\) &   \(\) &   \(\) &   \(\) &    \(\) &    \(3\) \\
\hline
 \(9q\)     &  \(d\equiv 1\,(\mathrm{mod}\,3)\) &    \(1\) &    \(\) &  \(\) &   \(\) & \(\) & \(\) &  \(\) &   \(\) &   \(\) &   \(\) &   \(\) &    \(\) &    \(0\) \\
\hline
\hline
 \(1\)      &                   \(\varrho_3=2\) &     \(\) &    \(\) &  \(\) &   \(\) &\(5\) &\(16\)&  \(\) &   \(\) &   \(\) &  \(4\) &   \(\) &    \(\) &   \(20\) \\
\hline
\hline
            & Summary                           &     \(\) &\(4652\) & \(9\) & \(21\) &\(5\) &\(16\)&\(10\) & \(76\) &\(106\) &\(3349\)& \(79\) &\(1117\) & \(\mathbf{4753}\) \\
\hline
\end{tabular}
\end{center}
\end{table}


According to Table
\ref{tbl:Angell},
the number of non-cyclic totally real cubic fields \(L\) with discriminant \(0<d_L<10^5\) is \(\mathbf{4753}\),
in perfect accordance with the results by Llorente and Oneto
\cite{LlOn1980,LlOn1982},
who discovered the ommission of ten fields in the table by Angell
\cite{An1975,An1976}.
Together with \(51\) cyclic cubic fields in Table
\ref{tbl:AngellCyclic},
the total number is \(\mathbf{4804}\)
(rather than \(4794\), as announced erroneously in
\cite{An1976}).


\renewcommand{\arraystretch}{1.1}

\begin{table}[ht]
\caption{Cyclic cubic discriminants \(d_L=f^2\) in the range \(0<d_L<10^5\)}
\label{tbl:AngellCyclic}
\begin{center}
\begin{tabular}{|rl||rr||r||rr|}
\hline
                  &                          & \multicolumn{2}{c||}{M} & DPF       &  \(f\) &    \(d_L\) \\
 \(f\)            & Condition                        &  \(1\) &  \(2\) & \(\zeta\) &        &            \\
\hline
 \(9\)            & \(d=1\)                          &  \(1\) &   \(\) &     \(1\) &  \(9\) &     \(81\) \\
 \(\ell\)         & \(\equiv +1\,(\mathrm{mod}\,3)\) & \(30\) &   \(\) &    \(30\) &  \(7\) &     \(49\) \\
\hline
 \(9\ell\)        & \(d=1\)                          &   \(\) &  \(4\) &     \(8\) & \(63\) & \(3\,969\) \\
 \(\ell_1\ell_2\) & \(\equiv +1\,(\mathrm{mod}\,3)\) &   \(\) &  \(6\) &    \(12\) & \(91\) & \(8\,281\) \\
\hline
                  & Summary                          & \(31\) & \(10\) &    \(51\) &        &            \\
\hline
\end{tabular}
\end{center}
\end{table}

\noindent
According to Table
\ref{tbl:AngellCyclic},
the number of \textit{cyclic} cubic fields \(L\) with discriminant \(0<d_L<10^5\) is \(\mathbf{51}\),
with \(31\) arising from \textit{singlets} having conductors \(f\) with a single prime divisor,
and \(20\) from \textit{doublets} having two prime divisors of the conductor \(f\).
(M denotes the multiplicity.)


Although we have given a succinct survey of the DPF types
of all \textit{multiplets} in Angell's range \(0<d_L<10^5\)
in the conclusion of
\cite{Ma2021},
we arrange them again in a more ostensive tabular form
with absolute frequency and minimal discriminant \(d_L=f^2\cdot d\).

All \textit{doublets} in Table
\ref{tbl:AngellDoublets}
are \textit{pure}.
In bigger ranges, there will also occur \textit{mixed} doublets,
e.g. in Table
\ref{tbl:LlorenteDoubletsTwo}.
The corresponding \(3\)-class rank is always \(\varrho=0\).

\renewcommand{\arraystretch}{1.1}

\begin{table}[ht]
\caption{Types of doublets in the range \(0<d_L<10^5\)}
\label{tbl:AngellDoublets}
\begin{center}
\begin{tabular}{|c||r||rrr|}
\hline
 DPF Type                         & Frequency &   \(d\) &  \(f\) &     \(d_L\) \\
 \((\tau(L_1),\tau(L_2))\)        &           &         &        &             \\
\hline
 \((\gamma,\gamma)\)              &     \(4\) &  \(33\) & \(45\) & \(66\,825\) \\
 \((\varepsilon,\varepsilon)\)    &     \(5\) & \(373\) & \(10\) & \(37\,300\) \\
\hline
                           Total: &     \(9\) &         &        &             \\
\hline
\end{tabular}
\end{center}
\end{table}


\noindent
The \textit{triplets} with \(\varrho=1\) in Table
\ref{tbl:AngellTriplets}
have been partially classified in a coarse sense by Schmithals in \(1985\)
\cite{Sm1985}.
He merely decided whether capitulation occurs or not,
indicating \(C=1\) by the symbol \lq\lq\(+\)\rq\rq\ and \(C=0\) by \lq\lq\(-\)\rq\rq.
This admits the detection of type \(\varepsilon\)
but fails to distinguish between the types \(\beta_1\) and \(\delta_1\).

\renewcommand{\arraystretch}{1.1}

\begin{table}[ht]
\caption{Types of triplets in the range \(0<d_L<10^5\)}
\label{tbl:AngellTriplets}
\begin{center}
\begin{tabular}{|c|c||r||rrr|}
\hline
             & DPF Type                                     & Frequency &      \(d\) &  \(f\) &     \(d_L\) \\
 \(\varrho\) & \((\tau(L_1),\ldots,\tau(L_3))\)             &           &            &        &             \\
\hline
 \(0\)       & \((\gamma,\gamma,\gamma)\)                   &     \(1\) &     \(69\) & \(18\) & \(22\,356\) \\
 \(0\)       & \((\varepsilon,\varepsilon,\varepsilon)\)    &     \(1\) &    \(717\) &  \(9\) & \(58\,077\) \\
\hline
 \(1\)       & \((\beta_1,\beta_1,\beta_1)\)                &     \(1\) & \(1\,509\) &  \(6\) & \(54\,324\) \\
 \(1\)       & \((\beta_1,\beta_1,\varepsilon)\)            &     \(2\) &\(14\,397\) &  \(2\) & \(57\,588\) \\
 \(1\)       & \((\beta_1,\delta_1,\delta_1)\)              &     \(3\) & \(1\,765\) &  \(7\) & \(86\,485\) \\
 \(1\)       & \((\delta_1,\delta_1,\delta_1)\)             &    \(13\) & \(7\,053\) &  \(2\) & \(28\,212\) \\
\hline
             &                                       Total: &    \(21\) &            &        &             \\
\hline
\end{tabular}
\end{center}
\end{table}


\noindent
The \textit{quartets} in Table
\ref{tbl:AngellQuartets}
belong to unramified cyclic cubic extensions of quadratic fields with \(\varrho=2\).
In fact, they have been classified by Heider and Schmithals in \(1982\)
\cite[p. 24]{HeSm1982}.
In \(2006\), resp. \(2008\), resp. \(2009\), we have detected the remaining capitulation numbers
\(\nu(K)=0\), resp. \(1\), resp. \(2\),
which show up in Table
\ref{tbl:LlorenteQuartetsUnramified}.
See
\cite{Ma2012,Ma2014b}.

\renewcommand{\arraystretch}{1.1}

\begin{table}[ht]
\caption{Types of quartets in the range \(0<d_L<10^5\)}
\label{tbl:AngellQuartets}
\begin{center}
\begin{tabular}{|c||c|r||r|}
\hline
 DPF Type                                  & Capitulation Number \(\nu(K)\)  & Frequency &   \(d_L=d\) \\
 \((\tau(L_1),\ldots,\tau(L_4))\)          & (according to \cite{ChFt1980})  &           &             \\
\hline
 \((\alpha_1,\alpha_1,\alpha_1,\alpha_1)\) & \(4\)                           &     \(1\) & \(62\,501\) \\
 \((\alpha_1,\alpha_1,\alpha_1,\delta_1)\) & \(3\)                           &     \(4\) & \(32\,009\) \\
\hline
                                           &                          Total: &     \(5\) &             \\
\hline
\end{tabular}
\end{center}
\end{table}

\newpage

\subsection{Numerical results by Nicole Moser}
\label{ss:NicoleMoser}

\noindent
In her paper
\cite{Mo1979}
on the units \(U_N\) and class groups \(\mathrm{Cl}_N\)
of dihedral fields \(N\) of degree \(2p\) with an odd prime \(p\),
Nicole Moser has given a small table
\cite[V.4, pp. 72--73]{Mo1979}
of \(34\) totally real cubic fields \(L\)
with discriminants \(0<d_L<1500\)
in order to illustrate her (coarse) classification by
concrete examples for \(p=3\).
She found \(26\) fields of type \(\delta\),
unramified with conductor \(f=1\), without exceptions,
and thus more precisely of our finer type \(\delta_1\).
The frequency \(\frac{26}{34}\approx 76\%\) corresponds to Angell's \(\frac{3349}{4753}\approx 70\%\).
Discriminants \(d_L=1^2\cdot d=d\) are
\[229,257,316,321,469,473,568,697,733,761,785,892,940,985,993,\]
\[1016,1101,\mathbf{1129},1229,1257,1304,1345,1373,1384,1436,1489.\]
In fact, each of the normal closures \(N\)
of degree \(3\) over its quadratic subfield \(K\) 
is precisely the Hilbert \(3\)-class field \(\mathrm{F}_3^1(K)\) of \(K\),
with one exception \(d=\mathbf{1129}\), where \(K\) has class number \(9\).
Here, Moser's table entries \(a=9\), \(h_N=9\) are incorrect
and must be replaced by \(a=3\), \(h_N=3\).
She uses two invariants,
the \textit{index of subfield units} \(a=(U_N:(U_K\cdot U_L\cdot U_{L^\sigma}))\) in Formula
\eqref{eqn:Scholz},
and the \textit{unit norm index} \(b=(U_K:N_{N/K}(U_N))\) in Remark
\ref{rmk:RealDPFTypes}
for her characterization of the types \(\alpha,\beta,\gamma,\delta,\varepsilon\).

Among the remaining \(8\) fields,
one is of type \(\gamma\) with \(d=21\equiv 3\,(\mathrm{mod}\,9)\), \(f=2\cdot 3\), \(d_L=6^2\cdot 21=756\),
and seven are of type \(\varepsilon\).
Among the latter,
five have \(f=2\) and \(d\in\lbrace 37,101,141,197,269\rbrace\), \(d_L\in\lbrace148,404,564,788,1076\rbrace\),
two have \(f=3\) and \(d=69\equiv 6\,(\mathrm{mod}\,9)\), \(d_L=621\),
resp. \(d=93\equiv 3\,(\mathrm{mod}\,9)\), \(d_L=837\).
The frequency \(\frac{7}{34}\approx 21\%\) corresponds to Angell's \(\frac{1117}{4753}\approx 24\%\).

However, it must be pointed out that \(4\) fields are \textit{missing}:
one is of type \(\gamma\) with \(d=13\), \(f=2\cdot 5\), \(d_L=10^2\cdot 13=1300\),
and three are of type \(\varepsilon\)
with \(d=349\), \(f=2\), \(d_L=1396\),
resp. \(d=57\), \(f=5\), \(d_L=1425\),
resp. \(d=373\), \(f=2\), \(d_L=1492\).

On the other hand, it is very instructive that there is also a \textit{superfluous} field:
although \(f=2\) is a \(3\)-admissible conductor for \(d=229\),
since \(2\) remains inert in \(K=\mathbb{Q}(\sqrt{d})\),
\(D=2^2\cdot d=4\cdot 229=916\) is only a \textit{formal} cubic discriminant,
because the defect of \(2\) is \(\delta(2)=1\).
So \(\mathbf{M}_{4d}=\emptyset\) is a \textit{nilet},
and the given polynomial \(X^3-X^2-6X+4\) generates
the cubic field with conductor \(f=1\) and discriminant \(d_L=229\).
In this case, Moser is uncertain whether the type
of the hypothetical cubic field with discriminant \(916\)
is \(\varepsilon\) or \(\gamma\).
Type \(\gamma\), however, is never possible
for a field with prime conductor, such as \(f=2\).
Since \(\varrho=1\) for \(d=229\), the types
\(\beta_1\), \(\delta_1\) and \(\varepsilon\) would be possible,
but Moser's claim \(a=9\) discourages types \(\beta_1\) and \(\delta_1\).
It is mysterious how she determined the invariant \(a\)
for a non-existent field
without knowing the class number \(h_N\).
In view of the errors in Moser's table,
it is worth ones while to state a summarizing theorem
which also pays attention to the modest contribution by cyclic cubic fields.
Cutting off Table
\ref{tbl:Angell}
at \(d_L=1500\) we obtain:


\begin{theorem}
\label{thm:NicoleMoser}
Among the \(44\) totally real cubic fields \(L\)
with discriminants \(0<d_L<1500\), there are
\(2\) (\(4\%\)) of type \(\gamma\),
\(26\) (\(59\%\)) of type \(\delta_1\), and
\(10\) (\(23\%\)) of type \(\varepsilon\).
These \(38\) non-Galois cubic fields are complemented by
\(6\) (\(14\%\)) cyclic cubic fields with conductors
\(f\in\lbrace 7,9,13,19,31,37\rbrace\).
With respect to the multiplicity \(m\),
all \(44\) fields form singlets with \(m=1\).
\end{theorem}


Obviously, Moser was not in possession of Angell's UMT
\cite{An1975},
otherwise she would have been able to detect the gaps in her table.
She rather refers to an unpublished table by Ren\'e Smadja.

Outside of the range \(0<d_L<1500\), Moser gives an example of a field
with type \(\beta\), more precisely our finer type \(\beta_2\),
for \(d=29\), \(f=2\cdot 7\), \(d_L=14^2\cdot 29=5684\).
The conductor is divisible by the prime \(7\) which splits in \(K\),
as required for type \(\beta_2\).

Moser did not know any examples of her type \(\alpha\).
From Table
\ref{tbl:AngellQuartets}
we know that the minimal discriminant of such a field is \(32009\),
discovered by Heider and Schmithals
\cite{HeSm1982}.
Due to \(f=1\)
it is more precisely our finer type \(\alpha_1\).
See Theorem
\ref{thm:Unramified}
and Example
\ref{exm:Unramified}.

\newpage

\section{Update of our \(1991\) classification for \(0<d_L<200\,000\)}
\label{s:GutensteinStreiteben}

\noindent
As announced in
\cite{Ma2021},
Table
\ref{tbl:Angell}
with \(0<d_L<10^5\)
was completed on Tuesday, 29 December 2020.
One week later, on Tuesday, 05 January 2021,
we finished the updated Table
\ref{tbl:GutensteinStreiteben}
containing the revised classification of all totally real cubic fields \(L\)
with discriminants \(0<d_L<2\cdot 10^5\),
which we had investigated in August 1991
\cite{Ma1991c}.

\renewcommand{\arraystretch}{1.1}

\begin{table}[ht]
\caption{Totally real cubic discriminants \(d_L=f^2\cdot d\) in the range \(0<d_L<2\cdot 10^5\)}
\label{tbl:GutensteinStreiteben}
\begin{center}
\begin{tabular}{|rl||r|rrrr||rrrrr|rrr|r|}
\hline
       &           & \multicolumn{5}{c||}{Multiplicity}    & \multicolumn{8}{c|}{Differential Principal Factorization} &       \\
 \(f\) & Condition & \(0\) & \(1\) & \(2\) & \(3\) & \(4\) & \(\alpha_1\) & \(\alpha_3\) & \(\beta_1\) & \(\beta_2\) & \(\gamma\) & \(\delta_1\) & \(\delta_2\) & \(\varepsilon\) & Total \\
\hline
 \(1\)          &                   \(\varrho_3=0\) &\(53848\) &     \(\) &  \(\) &   \(\) &   \(\) &  \(\) &  \(\) &  \(\) &   \(\) &   \(\) &   \(\) &    \(\) &     \(\) &    \(0\) \\
\hline
 \(q\)          &   \(\equiv 2\,(\mathrm{mod}\,3)\) & \(4338\) & \(1656\) &  \(\) &   \(\) &   \(\) &  \(\) &  \(\) &  \(\) &   \(\) &   \(\) &   \(\) &    \(\) & \(1656\) & \(1656\) \\
 \(3\)          &  \(d\equiv 3\,(\mathrm{mod}\,9)\) &  \(551\) &  \(221\) &  \(\) &   \(\) &   \(\) &  \(\) &  \(\) &  \(\) &   \(\) &   \(\) &   \(\) &    \(\) &  \(221\) &  \(221\) \\
 \(3\)          &  \(d\equiv 6\,(\mathrm{mod}\,9)\) &  \(561\) &  \(222\) &  \(\) &   \(\) &   \(\) &  \(\) &  \(\) &  \(\) &   \(\) &   \(\) &   \(\) &    \(\) &  \(222\) &  \(222\) \\
 \(9\)          &  \(d\equiv 6\,(\mathrm{mod}\,9)\) &   \(22\) &   \(68\) &  \(\) &  \(3\) &   \(\) &  \(\) &  \(\) &  \(\) &   \(\) &   \(\) &   \(\) &    \(\) &   \(77\) &   \(77\) \\
 \(9\)          &  \(d\equiv 2\,(\mathrm{mod}\,3)\) &  \(192\) &   \(74\) &  \(\) &   \(\) &   \(\) &  \(\) &  \(\) &  \(\) &   \(\) &   \(\) &   \(\) &    \(\) &   \(74\) &   \(74\) \\
\hline
 \(9\)          &  \(d\equiv 1\,(\mathrm{mod}\,3)\) &  \(184\) &   \(71\) &  \(\) &   \(\) &   \(\) &  \(\) &  \(\) &  \(\) & \(17\) &   \(\) &   \(\) &  \(48\) &    \(6\) &   \(71\) \\
 \(\ell\)       &   \(\equiv 1\,(\mathrm{mod}\,3)\) &  \(624\) &  \(196\) &  \(\) &   \(\) &   \(\) &  \(\) &  \(\) &  \(\) & \(45\) &   \(\) &   \(\) & \(140\) &   \(11\) &  \(196\) \\
\hline
 \(q_1q_2\)     &                                   &   \(65\) &   \(66\) & \(2\) &   \(\) &   \(\) &  \(\) &  \(\) &  \(\) &   \(\) & \(66\) &   \(\) &    \(\) &    \(4\) &   \(70\) \\
 \(3q\)         &  \(d\equiv 3\,(\mathrm{mod}\,9)\) &   \(45\) &   \(46\) & \(1\) &   \(\) &   \(\) &  \(\) &  \(\) &  \(\) &   \(\) & \(46\) &   \(\) &    \(\) &    \(2\) &   \(48\) \\
 \(3q\)         &  \(d\equiv 6\,(\mathrm{mod}\,9)\) &   \(35\) &   \(45\) & \(6\) &   \(\) &   \(\) &  \(\) &  \(\) &  \(\) &   \(\) & \(45\) &   \(\) &    \(\) &   \(12\) &   \(57\) \\
 \(9q\)         &  \(d\equiv 6\,(\mathrm{mod}\,9)\) &     \(\) &     \(\) & \(8\) &  \(1\) &   \(\) &  \(\) &  \(\) &  \(\) &   \(\) & \(15\) &   \(\) &    \(\) &    \(4\) &   \(19\) \\
 \(9q\)         &  \(d\equiv 2\,(\mathrm{mod}\,3)\) &   \(15\) &   \(16\) &  \(\) &   \(\) &   \(\) &  \(\) &  \(\) &  \(\) &   \(\) & \(16\) &   \(\) &    \(\) &     \(\) &   \(16\) \\
\hline
 \(9q\)         &  \(d\equiv 1\,(\mathrm{mod}\,3)\) &   \(10\) &   \(20\) &  \(\) &   \(\) &   \(\) &  \(\) &  \(\) &  \(\) & \(20\) &   \(\) &   \(\) &    \(\) &     \(\) &   \(20\) \\
 \(q\ell\)      &                                   &   \(30\) &   \(57\) & \(3\) &   \(\) &   \(\) &  \(\) &  \(\) &  \(\) & \(56\) &  \(1\) &   \(\) &    \(\) &    \(6\) &   \(63\) \\
 \(3\ell\)      &  \(d\equiv 3\,(\mathrm{mod}\,9)\) &    \(4\) &    \(9\) &  \(\) &   \(\) &   \(\) &  \(\) &  \(\) &  \(\) &  \(9\) &   \(\) &   \(\) &    \(\) &     \(\) &    \(9\) \\
 \(3\ell\)      &  \(d\equiv 6\,(\mathrm{mod}\,9)\) &    \(6\) &    \(5\) &  \(\) &   \(\) &   \(\) &  \(\) &  \(\) &  \(\) &  \(5\) &   \(\) &   \(\) &    \(\) &     \(\) &    \(5\) \\
 \(9\ell\)      &  \(d\equiv 2\,(\mathrm{mod}\,3)\) &    \(1\) &    \(3\) &  \(\) &   \(\) &   \(\) &  \(\) &  \(\) &  \(\) &  \(3\) &   \(\) &   \(\) &    \(\) &     \(\) &    \(3\) \\
\hline
 \(9\ell\)      &  \(d\equiv 1\,(\mathrm{mod}\,3)\) &     \(\) &    \(1\) &  \(\) &   \(\) &   \(\) &  \(\) & \(1\) &  \(\) &   \(\) &   \(\) &   \(\) &    \(\) &     \(\) &    \(1\) \\
\hline
 \(q_1q_2q_3\)  &                                   &     \(\) &     \(\) & \(1\) &   \(\) &   \(\) &  \(\) &  \(\) &  \(\) &   \(\) &  \(2\) &   \(\) &    \(\) &     \(\) &    \(2\) \\
 \(3q_1q_2\)    &  \(d\equiv 3\,(\mathrm{mod}\,9)\) &     \(\) &    \(1\) & \(1\) &   \(\) &   \(\) &  \(\) &  \(\) &  \(\) &   \(\) &  \(3\) &   \(\) &    \(\) &     \(\) &    \(3\) \\
 \(3q_1q_2\)    &  \(d\equiv 6\,(\mathrm{mod}\,9)\) &     \(\) &     \(\) & \(1\) &   \(\) &   \(\) &  \(\) &  \(\) &  \(\) &   \(\) &  \(2\) &   \(\) &    \(\) &     \(\) &    \(2\) \\
\hline
 \(9q_1q_2\)    &  \(d\equiv 1\,(\mathrm{mod}\,3)\) &     \(\) &    \(1\) &  \(\) &   \(\) &   \(\) &  \(\) &  \(\) &  \(\) &   \(\) &  \(1\) &   \(\) &    \(\) &     \(\) &    \(1\) \\
 \(q_1q_2\ell\) &                                   &     \(\) &     \(\) & \(1\) &   \(\) &   \(\) &  \(\) &  \(\) &  \(\) &   \(\) &  \(2\) &   \(\) &    \(\) &     \(\) &    \(2\) \\
 \(3q\ell\)     &  \(d\equiv 3\,(\mathrm{mod}\,9)\) &     \(\) &     \(\) & \(1\) &   \(\) &   \(\) &  \(\) &  \(\) &  \(\) &   \(\) &  \(2\) &   \(\) &    \(\) &     \(\) &    \(2\) \\
\hline
\hline
 \(1\)          &                   \(\varrho_3=1\) &     \(\) & \(6924\) &  \(\) &   \(\) &   \(\) &  \(\) &  \(\) &  \(\) &   \(\) &   \(\) &\(6924\)&    \(\) &     \(\) & \(6924\) \\
\hline
 \(q\)          &   \(\equiv 2\,(\mathrm{mod}\,3)\) &  \(576\) &     \(\) &  \(\) & \(28\) &   \(\) &  \(\) &  \(\) &\(12\) &   \(\) &   \(\) & \(68\) &    \(\) &    \(4\) &   \(84\) \\
 \(3\)          &  \(d\equiv 3\,(\mathrm{mod}\,9)\) &   \(66\) &     \(\) &  \(\) &  \(2\) &   \(\) &  \(\) &  \(\) & \(2\) &   \(\) &   \(\) &  \(3\) &    \(\) &    \(1\) &    \(6\) \\
 \(3\)          &  \(d\equiv 6\,(\mathrm{mod}\,9)\) &   \(63\) &     \(\) &  \(\) &  \(2\) &   \(\) &  \(\) &  \(\) &  \(\) &   \(\) &   \(\) &  \(6\) &    \(\) &     \(\) &    \(6\) \\
 \(9\)          &  \(d\equiv 6\,(\mathrm{mod}\,9)\) &    \(3\) &     \(\) &  \(\) &   \(\) &   \(\) &  \(\) &  \(\) &  \(\) &   \(\) &   \(\) &   \(\) &    \(\) &     \(\) &    \(0\) \\
 \(9\)          &  \(d\equiv 2\,(\mathrm{mod}\,3)\) &   \(16\) &     \(\) &  \(\) &  \(1\) &   \(\) &  \(\) &  \(\) &  \(\) &   \(\) &   \(\) &  \(3\) &    \(\) &     \(\) &    \(3\) \\
\hline
 \(9\)          &  \(d\equiv 1\,(\mathrm{mod}\,3)\) &   \(20\) &     \(\) &  \(\) &   \(\) &   \(\) &  \(\) &  \(\) &  \(\) &   \(\) &   \(\) &   \(\) &    \(\) &     \(\) &    \(0\) \\
 \(\ell\)       &   \(\equiv 1\,(\mathrm{mod}\,3)\) &   \(51\) &     \(\) &  \(\) &  \(5\) &   \(\) &  \(\) &  \(\) & \(4\) &   \(\) &   \(\) & \(10\) &    \(\) &    \(1\) &   \(15\) \\
\hline
 \(q_1q_2\)     &                                   &    \(3\) &     \(\) &  \(\) &   \(\) &   \(\) &  \(\) &  \(\) &  \(\) &   \(\) &   \(\) &   \(\) &    \(\) &     \(\) &    \(0\) \\
 \(3q\)         &  \(d\equiv 3\,(\mathrm{mod}\,9)\) &    \(3\) &     \(\) &  \(\) &   \(\) &   \(\) &  \(\) &  \(\) &  \(\) &   \(\) &   \(\) &   \(\) &    \(\) &     \(\) &    \(0\) \\
 \(3q\)         &  \(d\equiv 6\,(\mathrm{mod}\,9)\) &    \(4\) &     \(\) &  \(\) &  \(1\) &   \(\) &  \(\) &  \(\) & \(3\) &   \(\) &   \(\) &   \(\) &    \(\) &     \(\) &    \(3\) \\
\hline
 \(9q\)         &  \(d\equiv 1\,(\mathrm{mod}\,3)\) &    \(2\) &     \(\) &  \(\) &   \(\) &   \(\) &  \(\) &  \(\) &  \(\) &   \(\) &   \(\) &   \(\) &    \(\) &     \(\) &    \(0\) \\
\hline
\hline
 \(1\)          &                   \(\varrho_3=2\) &     \(\) &     \(\) &  \(\) &   \(\) & \(16\) & \(50\)&  \(\) &  \(\) &   \(\) &   \(\) & \(14\) &    \(\) &     \(\) &   \(64\) \\
\hline
\hline
                & Summary                           &     \(\) & \(9702\) &\(25\) & \(43\) & \(16\) & \(50\)& \(1\) &\(21\) &\(155\) &\(201\) &\(7028\)& \(188\) & \(2301\) &\(\mathbf{9945}\) \\
\hline
\end{tabular}
\end{center}
\end{table}


\noindent
According to Table
\ref{tbl:GutensteinStreiteben},
the total number of all \textit{non-cyclic} totally real cubic fields \(L\) with discriminants \(d<2\cdot 10^5\)
is \(\mathbf{9\,945}\).
Together with \(70\) \textit{cyclic} cubic fields in Table
\ref{tbl:GutensteinCyclic}
the number is \(\mathbf{10\,015}\),

\renewcommand{\arraystretch}{1.1}

\begin{table}[hb]
\caption{Cyclic cubic discriminants \(d_L=f^2\) in the range \(0<d_L<2\cdot 10^5\)}
\label{tbl:GutensteinCyclic}
\begin{center}
\begin{tabular}{|rl||rr||r|}
\hline
                        &            & \multicolumn{2}{c||}{M} & \multicolumn{1}{c|}{DPF} \\
 \(f\)                  & Condition                       &   \(1\) &   \(2\) & \(\zeta\) \\
\hline
 \(9\)                  & \(d=1\)                         &   \(1\) &    \(\) &     \(1\) \\
 \(\ell\)               & \(\equiv 1\,(\mathrm{mod}\,3)\) &  \(41\) &    \(\) &    \(41\) \\
\hline
 \(9\ell\)              & \(d=1\)                         &    \(\) &   \(6\) &    \(12\) \\
 \(\ell_1\ell_2\)       & \(\equiv 1\,(\mathrm{mod}\,3)\) &    \(\) &   \(8\) &    \(16\) \\
\hline
                        & Summary                         &  \(42\) &  \(14\) &    \(70\) \\
\hline
\end{tabular}
\end{center}
\end{table}

\noindent
According to Table
\ref{tbl:GutensteinCyclic},
the number of \textit{cyclic} cubic fields \(L\) with discriminant \(0<d_L<2\cdot 10^5\) is \(\mathbf{70}\),
with \(42\) arising from \textit{singlets} having conductors \(f\) with a single prime divisor,
and \(28\) from \textit{doublets} having two prime divisors of the conductor \(f\).
(M denotes the multiplicity.)

As predicted in the introduction,
several fields of type \(\beta_2\) were unduly classified as type \(\delta_2\),
since Voronoi's algorithm
\cite{Vo1896}
did not find any absolute principal factors
along the chains of lattice minima.
Since the table
\cite[p. 3]{Ma1991c}
accumulates information on conductors with similar behavior
and thus has a format different from Table
\ref{tbl:GutensteinStreiteben},
we compile a translation of critical rows in both tables,
subtracting contributions by cyclic cubic fields in
\cite{Ma1991c}.

\renewcommand{\arraystretch}{1.1}

\begin{table}[ht]
\caption{Accumulation of rows in Table \ref{tbl:GutensteinStreiteben} for comparison with \cite[p. 3]{Ma1991c}}
\label{tbl:Translation}
\begin{center}
\begin{tabular}{|rl||rrrr||rrrrr|rrr|r|}
\hline
       &           & \multicolumn{4}{c||}{Multiplicity}    & \multicolumn{8}{c|}{Differential Principal Factorization} &       \\
 \(f\) & Condition & \(1\) & \(2\) & \(3\) & \(4\) & \(\alpha_1\) & \(\alpha_3\) & \(\beta_1\) & \(\beta_2\) & \(\gamma\) & \(\delta_1\) & \(\delta_2\) & \(\varepsilon\) & Total \\
\hline
 \(9\)          &   \(d\equiv 1\,(\mathrm{mod}\,3)\), \(\varrho=0\) &   \(71\) &   \(\) &   \(\) &   \(\) &   \(\) &  \(\) &            \(\) &           \(17\) &    \(\) &              \(\) &           \(48\) &             \(6\) &   \(71\) \\
 \(\ell\)       &    \(\equiv 1\,(\mathrm{mod}\,3)\), \(\varrho=0\) &  \(196\) &   \(\) &   \(\) &   \(\) &   \(\) &  \(\) &            \(\) &           \(45\) &    \(\) &              \(\) &          \(140\) &            \(11\) &  \(196\) \\
\hline
                & Together                                          &  \(267\) &   \(\) &   \(\) &   \(\) &   \(\) &  \(\) &            \(\) &           \(62\) &    \(\) &              \(\) &          \(188\) &            \(17\) &  \(267\) \\
\hline
 \(\ell\)       &    \(\equiv 1\,(\mathrm{mod}\,3)\), \(\varrho=1\) &     \(\) &   \(\) &  \(5\) &   \(\) &   \(\) &  \(\) &           \(4\) &             \(\) &    \(\) &            \(10\) &             \(\) &             \(1\) &   \(15\) \\
\hline
                & Summary                                           & \(9702\) & \(25\) & \(43\) & \(16\) & \(50\) & \(1\) &          \(21\) &          \(155\) & \(201\) &          \(7028\) &          \(188\) &          \(2301\) & \(9945\) \\
\hline
\hline
 \cite{Ma1991c} & \(m=1\), \(t=s=1\), \(\varrho=0\)                 &  \(267\) &   \(\) &   \(\) &   \(\) &   \(\) &  \(\) &            \(\) &  \(\mathbf{53}\) &    \(\) &              \(\) & \(\mathbf{198}\) &   \(\mathbf{16}\) &  \(267\) \\
\hline
 \cite{Ma1991c} & \(m=3\), \(t=s=1\), \(\varrho=1\)                 &     \(\) &   \(\) &  \(5\) &   \(\) &   \(\) &  \(\) &  \(\mathbf{3}\) &             \(\) &    \(\) &   \(\mathbf{11}\) &             \(\) &             \(1\) &   \(15\) \\
\hline
 \cite{Ma1991c} & total                                             & \(9702\) & \(25\) & \(43\) & \(16\) & \(50\) & \(1\) & \(\mathbf{20}\) & \(\mathbf{146}\) & \(201\) & \(\mathbf{7029}\) & \(\mathbf{198}\) & \(\mathbf{2300}\) & \(9945\) \\
\hline
\end{tabular}
\end{center}
\end{table}

Table
\ref{tbl:Translation}
shows that the failures in
\cite[p. 3]{Ma1991c}
are located in two rows
concerning singlets, resp. triplets, with conductors \(f\) divisible by one prime
which splits in the real quadratic field \(K\) with \(3\)-class rank \(\varrho=0\), resp. \(\varrho=1\),
that is, \(m=1\), resp. \(m=3\), and \(t=s=1\), \(n=0\).
For \(\varrho=0\), there are \(10\) fields of type \(\delta_2\) too much,
\(9\) fields of type \(\beta_2\) too less, and \(1\) field of type \(\varepsilon\) too less.
They concern
\cite[Part I, \(\varrho=0\), Section 3.1--5 and 3.14]{Ma1991c}
and are corrected in Table
\ref{tbl:Corrections}.
For \(\varrho=1\), there is \(1\) field of type \(\delta_1\) too much,
which is of correct type \(\beta_1\).
Consequently, there is a corresponding erroneous impact on the row
\lq\lq total\rq\rq.
Wrong counters in
\cite{Ma1991c}
are printed with \textbf{boldface} font.

\renewcommand{\arraystretch}{1.1}

\begin{table}[ht]
\caption{Corrections in the range \(0<d_L<2\cdot10^5\)}
\label{tbl:Corrections}
\begin{center}
\begin{tabular}{|r|r|r||c|c|}
\hline
 No.    & \(d_L\)      & \(f\)   & Erroneous Type  & Correct Type    \\
        &              &         & \cite{Ma1991c}  & \(\tau(L)\)     \\
\hline
  \(1\) &  \(96\,481\) &   \(7\) & \(\delta_2\)    & \(\beta_2\)     \\
  \(2\) &  \(98\,833\) &   \(7\) & \(\beta_2\)     & \(\varepsilon\) \\
  \(3\) & \(160\,377\) &   \(7\) & \(\varepsilon\) & \(\beta_2\)     \\
  \(4\) & \(179\,144\) &   \(7\) & \(\delta_2\)    & \(\beta_2\)     \\
\hline
  \(5\) & \(130\,329\) &   \(9\) & \(\beta_2\)    & \(\varepsilon\)  \\
\hline
  \(6\) &  \(65\,741\) &  \(13\) & \(\delta_2\)    & \(\beta_2\)     \\
  \(7\) & \(110\,357\) &  \(13\) & \(\delta_2\)    & \(\beta_2\)     \\
  \(8\) & \(114\,413\) &  \(13\) & \(\delta_2\)    & \(\beta_2\)     \\
  \(9\) & \(125\,736\) &  \(13\) & \(\delta_2\)    & \(\beta_2\)     \\
\hline
 \(10\) & \(193\,857\) &  \(19\) & \(\delta_2\)    & \(\beta_2\)     \\
\hline
 \(11\) &  \(93\,217\) &  \(31\) & \(\delta_2\)    & \(\beta_2\)     \\
 \(12\) & \(134\,540\) &  \(31\) & \(\delta_2\)    & \(\beta_2\)     \\
\hline
 \(13\) & \(114\,005\) & \(151\) & \(\delta_2\)    & \(\beta_2\)     \\
\hline
\end{tabular}
\end{center}
\end{table}

We also give arithmetical invariants of the single erroneous triplet explicitly:
For \(d=568\equiv 1\,(\mathrm{mod}\,3)\), \(f=\ell=13\equiv 1\,(\mathrm{mod}\,3)\) and \(d_L=f^2\cdot d=95\,992\),
we have the pure type \((\delta_1,\delta_1,\delta_1)\) in
\cite[Part III, \(\varrho=1\), Section 2.2]{Ma1991c}
instead of the correct type \((\beta_1,\delta_1,\delta_1)\),
since the Voronoi algorithm did not find an absolute principal factor
and we concluded \(A=0\) instead of the correct \(A=1\).
In the coarse classification of Schmithals
\cite{Sm1985},
there is no difference between these two types,
since both are characterized by \((+,+,+)\).

\newpage

\subsection{Gain of arithmetical structure for \(100\,000<d_L<200\,000\)}
\label{ss:SecondGain}

\noindent
Over real quadratic fields \(K\) with \(3\)-class rank \(\varrho=0\),
this range brings some \textit{most remarkable gains}:

\begin{enumerate}
\item
first doublet for \(f=3q\), \(d\equiv 3\,(9)\), and first type \((\varepsilon,\varepsilon)\)
(\(d=5\,277\), \(q=2\), \(d_L=189\,972\)),
\item
first type \(\gamma\) for \(f=q\ell\)
(\(d=21\), \(q=2\), \(\ell=43\), \(f=86\), \(d_L=155\,316\), singlet \((\gamma)\)),
\item
first nilet for \(f=9\ell\), \(d\equiv 2\,(3)\)
(\(d=29\), \(\ell=7\), \(f=63\), \(D=115\,101\)),
\item
first occurrence of a conductor divisible by two splitting primes, \(s=2\),
namely \(f=9\ell\), \(d\equiv 1\,(3)\),
exploited immediately by the first occurrence of type \(\alpha_3\) which requires \(s\ge 2\),
and consequently the \textbf{first verification of the Scholz conjecture}
with unexpected two-dimensional relative principal factorization
instead of two-dimensional capitulation \\
(\(d=37\), \(\ell=7\), \(f=63\), \(d_L=146\,853\), singlet \((\alpha_3)\),
see Theorems \ref{thm:Alpha} and \ref{thm:Alpha3Ramified}),
\item
first occurrence of \(f=q_1q_2q_3\) as a doublet of type \((\gamma,\gamma)\)
(\(d=13\), \(f=110\), \(d_L=157\,300\)),
\item
first occurrence of \(f=3q_1q_2\), \(d\equiv 6\,(9)\), as a doublet of type \((\gamma,\gamma)\) \\
(\(d=213\), \(q_1=2\), \(q_2=5\), \(f=30\), \(d_L=191\,700\)),
\item
first occurrence of \(f=9q_1q_2\), \(d\equiv 1\,(3)\), as a singlet of type \((\gamma)\) \\
(\(d=13\), \(q_1=2\), \(q_2=5\), \(f=90\), \(d_L=105\,300\)),
\item
first occurrence of \(f=q_1q_2\ell\) as a doublet of type \((\gamma,\gamma)\)
(\(d=37\), \(f=70\), \(d_L=181\,300\)),
\item
first occurrence of \(f=3q\ell\), \(d\equiv 3\,(9)\), as a doublet of type \((\gamma,\gamma)\) \\
(\(d=93\), \(q=2\), \(\ell=7\), \(f=42\), \(d_L=164\,052\)).
\end{enumerate}

\noindent
The following phenomena arise within \(3\)-ring class fields \(K_f\), \(f>1\),
over real quadratic fields \(K\) with \(3\)-class rank \(\varrho=1\):

\begin{enumerate}
\item
first occurrence of \(f=3\), \(d\equiv 3\,(9)\), as triplets containing types \(\beta_1,\delta_1,\varepsilon\) \\
(\(d=12\,081\), \(d_L=108\,729\), type \((\delta_1,\delta_1,\delta_1)\), and
\(d=19\,749\), \(d_L=177\,741\), type \((\beta_1,\beta_1,\varepsilon)\)),
\item
first occurrence of \(f=9\), \(d\equiv 2\,(3)\), as a triplet \((\delta_1,\delta_1,\delta_1)\)
(\(d=1\,901\), \(d_L=153\,981\)),
\item
first type \(\varepsilon\) in a triplet with \(f=\ell\)
(\(d=3\,873\), \(\ell=7\), \(d_L=189\,777\), type \((\delta_1,\delta_1,\varepsilon)\)).
\end{enumerate}


\begin{example}
\label{exm:SecondGain}
\textit{Hetero}geneous multiplets arise from \(d=37\) with \(\varrho=0\)
and \(3\)-Selmer space \(V=\langle\eta\rangle\),
where \(\eta\) denotes the fundamental unit of \(K=\mathbb{Q}(\sqrt{37})\).
For (4), we have a singlet
\(\mathrm{Inv}(K_{63})=\lbrack\emptyset;\emptyset,\emptyset;(\alpha_3)\rbrack\),
corresponding to the divisors \((1;7,9;63)\) of \(f=63\),
since \(V(7)=V(9)=0\).
For (8), we have a quartet
\(\mathrm{Inv}(K_{70})=\lbrack\emptyset;(\varepsilon),\emptyset,\emptyset;\emptyset,(\beta_2),\emptyset;(\gamma,\gamma)\rbrack\),
corresponding to \((1;2,5,7;10,35,14;70)\), the divisors of \(f=70\),
since \(V(2)=V\) but \(V(5)=V(7)=0\).
Recall that \(148=2^2\cdot 37\)
is the well-known minimum of all non-cyclic positive cubic discriminants.
Its type is the singlet \((\varepsilon)\).
The type of \(45\,325=35^2\cdot 37\)
is the singlet \((\beta_2)\).
Padding \textit{nilets} \(\emptyset\) illuminate the arithmetical structure.
\end{example}

\newpage

\section{Classifying Ennola and Turunen's range \(0<d_L<500\,000\)}
\label{s:EnnolaTurunen}

\noindent
The increasing contributions by new types of conductors \(f\)
in the range \(0<d_L<5\cdot 10^5\)
enforce a splitting into Table
\ref{tbl:EnnolaTurunen0}
for \(\varrho=0\) and Table
\ref{tbl:EnnolaTurunen1}
for \(\varrho\in\lbrace 1,2\rbrace\).

\renewcommand{\arraystretch}{1.1}

\begin{table}[ht]
\caption{Totally real cubic discriminants \(d_L=f^2\cdot d\) in the range \(0<d_L<5\cdot 10^5\)}
\label{tbl:EnnolaTurunen0}
\begin{center}
\begin{tabular}{|rl||r|rrrr||rrrrr|rrr|r|}
\hline
       &           & \multicolumn{5}{c||}{Multiplicity}    & \multicolumn{8}{c|}{Differential Principal Factorization} &       \\
 \(f\) & Condition & \(0\) & \(1\) & \(2\) & \(3\) & \(4\) & \(\alpha_1\) & \(\alpha_3\) & \(\beta_1\) & \(\beta_2\) & \(\gamma\) & \(\delta_1\) & \(\delta_2\) & \(\varepsilon\) & Total \\
\hline
 \(1\)            &     \(\varrho_3=0\) &\(133534\)&     \(\) &  \(\) &   \(\) &   \(\) &  \(\) &  \(\) &  \(\) &   \(\) &   \(\) &   \(\) &    \(\) &     \(\) &    \(0\) \\
\hline
 \(q\)            &   \(\equiv 2\,(3)\) &\(10515\) & \(4296\) &  \(\) &   \(\) &   \(\) &  \(\) &  \(\) &  \(\) &   \(\) &   \(\) &   \(\) &    \(\) & \(4296\) & \(4296\) \\
 \(3\)            &  \(d\equiv 3\,(9)\) & \(1339\) &  \(573\) &  \(\) &   \(\) &   \(\) &  \(\) &  \(\) &  \(\) &   \(\) &   \(\) &   \(\) &    \(\) &  \(573\) &  \(573\) \\
 \(3\)            &  \(d\equiv 6\,(9)\) & \(1364\) &  \(554\) &  \(\) &   \(\) &   \(\) &  \(\) &  \(\) &  \(\) &   \(\) &   \(\) &   \(\) &    \(\) &  \(554\) &  \(554\) \\
 \(9\)            &  \(d\equiv 6\,(9)\) &   \(46\) &  \(166\) &  \(\) &  \(7\) &   \(\) &  \(\) &  \(\) &  \(\) &   \(\) &   \(\) &   \(\) &    \(\) &  \(187\) &  \(187\) \\
 \(9\)            &  \(d\equiv 2\,(3)\) &  \(453\) &  \(197\) &  \(\) &   \(\) &   \(\) &  \(\) &  \(\) &  \(\) &   \(\) &   \(\) &   \(\) &    \(\) &  \(197\) &  \(197\) \\
\hline
 \(9\)            &  \(d\equiv 1\,(3)\) &  \(468\) &  \(178\) &  \(\) &   \(\) &   \(\) &  \(\) &  \(\) &  \(\) & \(45\) &   \(\) &   \(\) & \(116\) &   \(17\) &  \(178\) \\
 \(\ell\)         &   \(\equiv 1\,(3)\) & \(1530\) &  \(539\) &  \(\) &   \(\) &   \(\) &  \(\) &  \(\) &  \(\) &\(124\) &   \(\) &   \(\) & \(374\) &   \(41\) &  \(539\) \\
\hline
 \(q_1q_2\)       &                     &  \(147\) &  \(159\) &\(10\) &   \(\) &   \(\) &  \(\) &  \(\) &  \(\) &   \(\) &\(161\) &   \(\) &    \(\) &   \(18\) &  \(179\) \\
 \(3q\)           &  \(d\equiv 3\,(9)\) &  \(102\) &  \(107\) & \(8\) &   \(\) &   \(\) &  \(\) &  \(\) &  \(\) &   \(\) &\(107\) &   \(\) &    \(\) &   \(16\) &  \(123\) \\
 \(3q\)           &  \(d\equiv 6\,(9)\) &   \(89\) &  \(109\) &\(14\) &   \(\) &   \(\) &  \(\) &  \(\) &  \(\) &   \(\) &\(109\) &   \(\) &    \(\) &   \(28\) &  \(137\) \\
 \(9q\)           &  \(d\equiv 6\,(9)\) &    \(2\) &     \(\) &\(20\) &  \(3\) &   \(\) &  \(\) &  \(\) &  \(\) &   \(\) & \(39\) &   \(\) &    \(\) &   \(10\) &   \(49\) \\
 \(9q\)           &  \(d\equiv 2\,(3)\) &   \(31\) &   \(41\) & \(1\) &   \(\) &   \(\) &  \(\) &  \(\) &  \(\) &   \(\) & \(41\) &   \(\) &    \(\) &    \(2\) &   \(43\) \\
\hline
 \(9q\)           &  \(d\equiv 1\,(3)\) &   \(30\) &   \(39\) & \(3\) &   \(\) &   \(\) &  \(\) &  \(\) &  \(\) & \(39\) &   \(\) &   \(\) &    \(\) &    \(6\) &   \(45\) \\
 \(q\ell\)        &                     &   \(83\) &  \(129\) & \(6\) &   \(\) &   \(\) &  \(\) &  \(\) &  \(\) &\(125\) &  \(4\) &   \(\) &    \(\) &   \(12\) &  \(141\) \\
 \(3\ell\)        &  \(d\equiv 3\,(9)\) &   \(12\) &   \(20\) & \(1\) &   \(\) &   \(\) &  \(\) &  \(\) &  \(\) & \(20\) &   \(\) &   \(\) &    \(\) &    \(2\) &   \(22\) \\
 \(3\ell\)        &  \(d\equiv 6\,(9)\) &   \(11\) &   \(16\) & \(1\) &   \(\) &   \(\) &  \(\) &  \(\) &  \(\) & \(16\) &   \(\) &   \(\) &    \(\) &    \(2\) &   \(18\) \\
 \(9\ell\)        &  \(d\equiv 6\,(9)\) &     \(\) &     \(\) & \(1\) &   \(\) &   \(\) &  \(\) &  \(\) &  \(\) &  \(2\) &   \(\) &   \(\) &    \(\) &     \(\) &    \(2\) \\
 \(9\ell\)        &  \(d\equiv 2\,(3)\) &    \(4\) &    \(8\) &  \(\) &   \(\) &   \(\) &  \(\) &  \(\) &  \(\) &  \(8\) &   \(\) &   \(\) &    \(\) &     \(\) &    \(8\) \\
\hline
 \(9\ell\)        &  \(d\equiv 1\,(3)\) &    \(3\) &    \(1\) &  \(\) &   \(\) &   \(\) &  \(\) & \(1\) &  \(\) &   \(\) &   \(\) &   \(\) &    \(\) &     \(\) &    \(1\) \\
 \(\ell_1\ell_2\) &                     &    \(1\) &    \(2\) &  \(\) &   \(\) &   \(\) &  \(\) & \(2\) &  \(\) &   \(\) &   \(\) &   \(\) &    \(\) &     \(\) &    \(2\) \\
\hline
 \(q_1q_2q_3\)    &                     &     \(\) &     \(\) & \(1\) &   \(\) &   \(\) &  \(\) &  \(\) &  \(\) &   \(\) &  \(2\) &   \(\) &    \(\) &     \(\) &    \(2\) \\
 \(3q_1q_2\)      &  \(d\equiv 3\,(9)\) &     \(\) &    \(2\) & \(3\) &   \(\) &   \(\) &  \(\) &  \(\) &  \(\) &   \(\) &  \(8\) &   \(\) &    \(\) &     \(\) &    \(8\) \\
 \(3q_1q_2\)      &  \(d\equiv 6\,(9)\) &     \(\) &    \(1\) & \(1\) &   \(\) &   \(\) &  \(\) &  \(\) &  \(\) &   \(\) &  \(3\) &   \(\) &    \(\) &     \(\) &    \(3\) \\
 \(9q_1q_2\)      &  \(d\equiv 2\,(3)\) &     \(\) &    \(2\) &  \(\) &   \(\) &   \(\) &  \(\) &  \(\) &  \(\) &   \(\) &  \(2\) &   \(\) &    \(\) &     \(\) &    \(2\) \\
\hline
 \(9q_1q_2\)      &  \(d\equiv 1\,(3)\) &     \(\) &    \(1\) & \(1\) &   \(\) &   \(\) &  \(\) &  \(\) &  \(\) &   \(\) &  \(3\) &   \(\) &    \(\) &     \(\) &    \(3\) \\
 \(q_1q_2\ell\)   &                     &     \(\) &    \(1\) & \(2\) &   \(\) &   \(\) &  \(\) &  \(\) &  \(\) &  \(1\) &  \(4\) &   \(\) &    \(\) &     \(\) &    \(5\) \\
 \(3q\ell\)       &  \(d\equiv 3\,(9)\) &     \(\) &    \(2\) & \(1\) &   \(\) &   \(\) &  \(\) &  \(\) &  \(\) &   \(\) &  \(4\) &   \(\) &    \(\) &     \(\) &    \(4\) \\
 \(3q\ell\)       &  \(d\equiv 6\,(9)\) &     \(\) &     \(\) & \(2\) &   \(\) &   \(\) &  \(\) &  \(\) &  \(\) &   \(\) &  \(4\) &   \(\) &    \(\) &     \(\) &    \(4\) \\
 \(9q\ell\)       &  \(d\equiv 2\,(3)\) &     \(\) &     \(\) & \(1\) &   \(\) &   \(\) &  \(\) &  \(\) &  \(\) &   \(\) &  \(2\) &   \(\) &    \(\) &     \(\) &    \(2\) \\
\hline
\hline
                  & Subtotal            &     \(\) & \(7143\) &\(77\) & \(10\) &   \(\) &  \(\) & \(3\) &  \(\) &\(380\) &\(493\) &   \(\) & \(490\) & \(5961\) & \(7327\) \\
\hline
\end{tabular}
\end{center}
\end{table}

\noindent
According to Table
\ref{tbl:EnnolaTurunen1},
the total number of all \textit{non-cyclic} totally real cubic fields \(L\) with discriminants \(d<5\cdot 10^5\)
is \(\mathbf{26\,330}\).
Together with \(110\) \textit{cyclic} cubic fields in Table
\ref{tbl:EnnolaCyclic}
the number is \(\mathbf{26\,440}\),

\newpage

\renewcommand{\arraystretch}{1.1}

\begin{table}[ht]
\caption{Table \ref{tbl:EnnolaTurunen0} with \(0<d_L<5\cdot 10^5\) continued for \(\varrho_3\ge 1\)}
\label{tbl:EnnolaTurunen1}
\begin{center}
\begin{tabular}{|rl||r|rrrr||rrrrr|rrr|r|}
\hline
       &           & \multicolumn{5}{c||}{Multiplicity}    & \multicolumn{8}{c|}{Differential Principal Factorization} &       \\
 \(f\) & Condition & \(0\) & \(1\) & \(2\) & \(3\) & \(4\) & \(\alpha_1\) & \(\alpha_3\) & \(\beta_1\) & \(\beta_2\) & \(\gamma\) & \(\delta_1\) & \(\delta_2\) & \(\varepsilon\) & Total \\
\hline
 \(1\)            &     \(\varrho_3=1\) &     \(\) &\(18378\) &  \(\) &   \(\) &   \(\) &  \(\) &  \(\) &  \(\) &   \(\) &   \(\) &\(18378\)&   \(\) &     \(\) &\(18378\) \\
\hline
 \(q\)            &   \(\equiv 2\,(3)\) & \(1603\) &     \(\) &  \(\) & \(92\) &   \(\) &  \(\) &  \(\) &\(48\) &   \(\) &   \(\) & \(213\) &    \(\) &   \(15\) &  \(276\) \\
 \(3\)            &  \(d\equiv 3\,(9)\) &  \(188\) &     \(\) &  \(\) &  \(8\) &   \(\) &  \(\) &  \(\) & \(7\) &   \(\) &   \(\) &  \(15\) &    \(\) &    \(2\) &   \(24\) \\
 \(3\)            &  \(d\equiv 6\,(9)\) &  \(190\) &     \(\) &  \(\) &  \(6\) &   \(\) &  \(\) &  \(\) & \(2\) &   \(\) &   \(\) &  \(15\) &    \(\) &    \(1\) &   \(18\) \\
 \(9\)            &  \(d\equiv 6\,(9)\) &   \(17\) &     \(\) &  \(\) &  \(2\) &   \(\) &  \(\) &  \(\) & \(2\) &   \(\) &   \(\) &   \(3\) &    \(\) &    \(1\) &    \(6\) \\
 \(9\)            &  \(d\equiv 2\,(3)\) &   \(54\) &     \(\) &  \(\) &  \(1\) &   \(\) &  \(\) &  \(\) &  \(\) &   \(\) &   \(\) &   \(3\) &    \(\) &     \(\) &    \(3\) \\
\hline
 \(9\)            &  \(d\equiv 1\,(3)\) &   \(53\) &     \(\) &  \(\) &  \(1\) &   \(\) &  \(\) &  \(\) & \(1\) &   \(\) &   \(\) &   \(2\) &    \(\) &     \(\) &    \(3\) \\
 \(\ell\)         &   \(\equiv 1\,(3)\) &  \(150\) &     \(\) &  \(\) &  \(8\) &   \(\) &  \(\) &  \(\) & \(5\) &   \(\) &   \(\) &  \(16\) &    \(\) &    \(3\) &   \(24\) \\
\hline
 \(q_1q_2\)       &                     &   \(15\) &     \(\) &  \(\) &  \(3\) &   \(\) &  \(\) &  \(\) & \(9\) &   \(\) &   \(\) &    \(\) &    \(\) &     \(\) &    \(9\) \\
 \(3q\)           &  \(d\equiv 3\,(9)\) &   \(11\) &     \(\) &  \(\) &  \(2\) &   \(\) &  \(\) &  \(\) & \(3\) &   \(\) &   \(\) &    \(\) &    \(\) &    \(3\) &    \(6\) \\
 \(3q\)           &  \(d\equiv 6\,(9)\) &   \(18\) &     \(\) &  \(\) &  \(2\) &   \(\) &  \(\) &  \(\) & \(6\) &   \(\) &   \(\) &    \(\) &    \(\) &     \(\) &    \(6\) \\
 \(9q\)           &  \(d\equiv 6\,(9)\) &    \(1\) &     \(\) &  \(\) &   \(\) &   \(\) &  \(\) &  \(\) &  \(\) &   \(\) &   \(\) &    \(\) &    \(\) &     \(\) &    \(0\) \\
 \(9q\)           &  \(d\equiv 2\,(3)\) &    \(1\) &     \(\) &  \(\) &  \(1\) &   \(\) &  \(\) &  \(\) & \(3\) &   \(\) &   \(\) &    \(\) &    \(\) &     \(\) &    \(3\) \\
\hline
 \(9q\)           &  \(d\equiv 1\,(3)\) &    \(3\) &     \(\) &  \(\) &   \(\) &   \(\) &  \(\) &  \(\) &  \(\) &   \(\) &   \(\) &    \(\) &    \(\) &     \(\) &    \(0\) \\
 \(q\ell\)        &                     &   \(10\) &     \(\) &  \(\) &  \(1\) &   \(\) &  \(\) &  \(\) & \(3\) &   \(\) &   \(\) &    \(\) &    \(\) &     \(\) &    \(3\) \\
 \(3\ell\)        &  \(d\equiv 3\,(9)\) &    \(1\) &     \(\) &  \(\) &   \(\) &   \(\) &  \(\) &  \(\) &  \(\) &   \(\) &   \(\) &    \(\) &    \(\) &     \(\) &    \(0\) \\
 \(3\ell\)        &  \(d\equiv 6\,(9)\) &    \(1\) &     \(\) &  \(\) &   \(\) &   \(\) &  \(\) &  \(\) &  \(\) &   \(\) &   \(\) &    \(\) &    \(\) &     \(\) &    \(0\) \\
\hline
\hline
 \(1\)            &     \(\varrho_3=2\) &     \(\) &     \(\) &  \(\) &   \(\) & \(61\) &\(175\)&  \(\) &  \(\) &   \(\) &   \(\) &  \(69\) &    \(\) &     \(\) &  \(244\) \\
\hline
\hline
                  & Subtotal            &     \(\) &\(18378\) &  \(\) &\(127\) & \(61\) &\(175\)&  \(\) &\(89\) &   \(\) &   \(\) &\(18714\)&    \(\) &   \(25\) &\(19003\) \\
\hline
\hline
                  & \textbf{Total}      &     \(\) &\(25521\) &\(77\) &\(137\) & \(61\) &\(175\)& \(3\) &\(89\) &\(380\) &\(493\) &\(18714\)& \(490\) & \(5986\) &\(\mathbf{26330}\)\\
\hline
\end{tabular}
\end{center}
\end{table}


\noindent
According to Table
\ref{tbl:EnnolaCyclic},
the number of \textit{cyclic} cubic fields \(L\) with discriminant \(0<d_L<5\cdot 10^5\) is \(\mathbf{110}\),
with \(60\) arising from \textit{singlets} having conductors \(f\) with a single prime divisor,
and \(50\) from \textit{doublets} having two prime divisors of the conductor \(f\).
(M denotes the multiplicity.)

\renewcommand{\arraystretch}{1.1}

\begin{table}[hb]
\caption{Cyclic cubic discriminants \(d_L=f^2\) in the range \(0<d_L<5\cdot 10^5\)}
\label{tbl:EnnolaCyclic}
\begin{center}
\begin{tabular}{|rl||rr||r|}
\hline
                        &            & \multicolumn{2}{c||}{M} & \multicolumn{1}{c|}{DPF} \\
 \(f\)                  & Condition                       &   \(1\) &   \(2\) & \(\zeta\) \\
\hline
 \(9\)                  & \(d=1\)                         &   \(1\) &    \(\) &     \(1\) \\
 \(\ell\)               & \(\equiv 1\,(\mathrm{mod}\,3)\) &  \(59\) &    \(\) &    \(59\) \\
\hline
 \(9\ell\)              & \(d=1\)                         &    \(\) &   \(9\) &    \(18\) \\
 \(\ell_1\ell_2\)       & \(\equiv 1\,(\mathrm{mod}\,3)\) &    \(\) &  \(16\) &    \(32\) \\
\hline
                        & Summary                         &  \(60\) &  \(25\) &   \(110\) \\
\hline
\end{tabular}
\end{center}
\end{table}


\begin{example}
\label{exm:SecondLine}
In Table
\ref{tbl:EnnolaTurunen0},
the second line with conductor \(f=q\), a prime number \(q\equiv 2\,(\mathrm{mod}\,3)\),
lists \(10515\) \textit{nilets}, starting with the \textit{formal} cubic discriminant \(f^2\cdot d=2^2\cdot 5=20\),
which does \textit{not} belong to an actual cubic field,
and \(4296\) singlets with minimal discriminant \(d_L=f^2\cdot d=2^2\cdot 37=148\) of an \textit{actual} cubic field \(L\).
Theoretical justifications for these facts are given in \cite[Thm. 4.1]{Ma2021}:
the \(3\)-\textit{Selmer space} \(V_3=\langle\eta\rangle\) of the real quadratic field \(K=\mathbb{Q}(\sqrt{d})\)
is generated by the fundamental unit \(\eta\in U_K\).
In the case of a nilet,
the \(3\)-\textit{ring space} mod \(q\), \(V_3(q)\), is the null space of codimension \(\delta_3(q)=1\) in \(V_3\),
since \(\eta\notin\mathcal{O}_q\).
In the case of a singlet, we have \(V_3(q)=V_3\) with \textit{defect} \(\delta_3(q)=0\).
\end{example}

\newpage

\noindent
Since minor counting errors have occurred in the tables by Moser, Angell and Llorente/Quer
(whereas the table by Ennola/Turunen was correct),
we explicitly state the ultimate counters of totally real cubic fields \(L\)
in five ranges of discriminants \(0<d_L<B\)
with various upper bounds \(B\). 

\begin{theorem}
\label{thm:Counters}
The number of cyclic, resp. non-Galois, resp. all,
non-isomorphic totally real cubic fields \(L\)
with discriminants in the range \(0<d_L<B\) is given by
\begin{enumerate}
\item
\(6\), resp. \(38\), resp. \(44\), for \(B=1\,500\),
\item
\(51\), resp. \(4\,753\), resp. \(4\,804\), for \(B=100\,000\),
\item
\(70\), resp. \(9\,945\), resp. \(10\,015\), for \(B=200\,000\),
\item
\(110\), resp. \(26\,330\), resp. \(26\,440\), for \(B=500\,000\),
\item
\(501\), resp. \(592\,421\), resp. \(592\,922\), for \(B=10\,000\,000\).
\end{enumerate}
\end{theorem}

\begin{proof}
See the tables in sections \S\S\
\ref{s:Angell}
--
\ref{s:LlorenteQuer}.
\end{proof}


\noindent
Recall that no examples of the types \(\alpha_2\) and \(\alpha_3\)
occurred in Angell's range \(0<d_L<10^5\),
and type \(\alpha_2\) remained unknown even in Ennola and Turunen's range \(0<d_L<5\cdot 10^5\).
Since this problem is intimately connected with the Scholz Conjecture in \S\
\ref{s:ScholzConjecture},
we now emphasize the following theorem.

\begin{theorem}
\label{thm:Alpha}
The minimal discriminants \(d_L=f^2\cdot d\) of totally real cubic fields \(L\)
with conductor \(f\) and quadratic fundamental discriminant \(d\)
such that \(\tau(L)\) is one of the extremely rare differential principal factorization types
\(\alpha_3\), resp. \(\alpha_2\), are given by
\begin{enumerate}
\item
\(146\,853\) with \(f=63=9\cdot 7\), \(s=2\), and \(d=37\), \(\varrho=0\) (unique field in a singlet, \(m=1\)), resp.
\item
\(966\,397\) with \(f=19\), \(s=1\), and \(d=2\,677\), \(\varrho=1\) (two of the fields in a triplet, \(m=3\)).
\end{enumerate}
\end{theorem}

\begin{proof}
The unique field \(L\) with discriminant \(146\,853\) has been discovered in August \(1991\) already
\cite[Part I, \(\varrho=0\), Section 6.1]{Ma1991c}
and was confirmed in the row with conductor \(f=9\ell\), \(d\equiv 1\,(\mathrm{mod}\,3)\), of Table
\ref{tbl:GutensteinStreiteben}.
According to Theorem
\ref{thm:Alpha3Ramified},
this field forms a singlet with DPF type \(\alpha_3\).

The triplet \((L_1,L_2,L_3)\) with discriminant \(966\,397\) was found by direct search on 19 November \(2017\).
It is now confirmed by gapless construction
in the row with conductor \(f=\ell\equiv 1\,(\mathrm{mod}\,3)\) for \(\varrho_3=1\) in Table
\ref{tbl:LlorenteQuer1}.
According to Theorem
\ref{thm:Alpha2Ramified},
the DPF type of the triplet is \((\alpha_2,\alpha_2,\delta_1)\).
\end{proof}


\subsection{Gain of arithmetical structure for \(200\,000<d_L<500\,000\)}
\label{ss:ThirdGain}

\noindent
The following new features arise within \(3\)-ring class fields \(K_f\), \(f>1\),
over real quadratic fields \(K\) with \(3\)-class rank \(\varrho=0\),

\begin{enumerate}
\item
first doublet of type \((\varepsilon,\varepsilon)\) for \(f=9q\), \(d\equiv 2\,(3)\)
(\(d=1\,157\), \(q=2\), \(d_L=374\,868\)),
\item
first doublet of type \((\varepsilon,\varepsilon)\) for \(f=9q\), \(d\equiv 1\,(3)\)
(\(d=877\), \(q=2\), \(d_L=284\,148\)),
\item
first doublet of type \((\varepsilon,\varepsilon)\) for \(f=3\ell\), \(d\equiv 3\,(9)\)
(\(d=597\), \(\ell=7\), \(d_L=263\,277\)),
\item
first doublet of type \((\varepsilon,\varepsilon)\) for \(f=3\ell\), \(d\equiv 6\,(9)\)
(\(d=1\,068\), \(\ell=7\), \(d_L=470\,988\)),
\item
first occurrence of \(f=9\ell\), \(d\equiv 6\,(9)\), as a doublet \((\beta_2,\beta_2)\)
(\(d=60\), \(\ell=7\), \(d_L=238\,140\)),
\item
first occurrence of \(f=\ell_1\ell_2\) with \(s=2\) as singlets of type \((\alpha_3)\) \\
(\(d=29\), \(f=91\), \(d_L=240\,149\) and \(d=8\), \(f=217\), \(d_L=376\,712\)),
\item
first singlet of type \((\gamma)\) for \(f=3q_1q_2\), \(d\equiv 6\,(9)\)
(\(d=357\), \(f=30\), \(d_L=321\,300\)),
\item
first occurrence of \(f=9q_1q_2\), \(d\equiv 2\,(3)\), as singlet of type \((\gamma)\)
(\(d=53\), \(f=90\), \(d_L=429\,300\)),
\item
first singlet of type \((\beta_2)\) for \(f=q_1q_2\ell\)
(\(d=93\), \(f=70\), \(d_L=455\,700\)),
\item
first singlet of type \((\gamma)\) for \(f=3q\ell\), \(d\equiv 3\,(9)\)
(\(d=165\), \(f=42\), \(d_L=291\,060\)),
\item
first occurrence of \(f=3q\ell\), \(d\equiv 6\,(9)\), as doublet of type \((\gamma,\gamma)\)
(\(f=42\), \(d_L=248\,724\)),
\item
first occurrence of \(f=9q\ell\), \(d\equiv 2\,(3)\), as doublet of type \((\gamma,\gamma)\)
(\(d=29\), \(d_L=460\,404\)).
\end{enumerate}

\noindent
The following phenomena arise within \(3\)-ring class fields \(K_f\), \(f>1\),
over \(K\) with \(\varrho=1\):

\begin{enumerate}
\item
first triplet of type \((\beta_1,\beta_1,\varepsilon)\) for \(f=3\), \(d\equiv 6\,(9)\)
(\(d=52\,197\), \(d_L=469\,773\)),
\item
first occurrence of \(f=9\), \(d\equiv 6\,(9)\), as triplet \((\beta_1,\beta_1,\varepsilon)\)
(\(d=5\,073\), \(d_L=410\,913\)),
\item
first occurrence of \(f=9\), \(d\equiv 1\,(3)\), as triplet \((\beta_1,\delta_1,\delta_1)\)
(\(d=2\,917\), \(d_L=236\,277\)),
\item
first occurrence of \(f=q_1q_2\), as triplet \((\beta_1,\beta_1,\beta_1)\)
(\(d=3\,173\), \(f=10\), \(d_L=317\,300\)),
\item
first occurrence of \(f=3q\), \(d\equiv 3\,(9)\), as triplet \((\beta_1,\beta_1,\beta_1)\)
(\(d=5\,637\), \(f=6\), \(d_L=202\,932\)),
\item
first occurrence of \(f=9q\), \(d\equiv 2\,(3)\), as triplet \((\beta_1,\beta_1,\beta_1)\)
(\(d=1\,373\), \(f=18\), \(d_L=444\,852\)),
\item
first occurrence of \(f=q\ell\), as triplet \((\beta_1,\beta_1,\beta_1)\)
(\(d=1\,101\), \(f=14\), \(d_L=215\,796\)).
\end{enumerate}

\noindent
First \textbf{unramified quartet} of type \((\delta_1,\delta_1,\delta_1,\delta_1)\) for \(\varrho=2\) (\(d=d_L=214\,712\)
\cite{Ma2012,Ma2014b}).

\newpage

\newgeometry{left=0.5cm, right=1cm, top=2cm, bottom=2cm, bindingoffset=5mm}

\section{Classifying Llorente and Quer's range \(0<d_L<10\,000\,000\)}
\label{s:LlorenteQuer}

\noindent
As opposed to the smaller ranges, the extension to Llorente and Quer's upper bound \(10^7\)
caused unexpected complications of two kinds.
Firstly, for ramified extensions with conductor \(f=2\cdot 9=18\), \(d\equiv 1\,(\mathrm{mod}\, 3)\),
at several discriminants \(d_L=f^2\cdot d>4\,941\,972\), \(d>15\,253\), \(\varrho_3=0\),
resp. \(d_L=f^2\cdot d>4\,249\,908\), \(d>13\,117\), \(\varrho_3=1\).
Secondly, for unramified extensions with \(f=1\) and \(\varrho_3=2\),
at several discriminants \(d_L=d>5\,547\,841\).
Thus, we were very releaved,
when a suitable work-around admitted the completion of the following most extensive and expensive Tables
\ref{tbl:LlorenteQuer0}
and
\ref{tbl:LlorenteQuer1}
on Wednesday, 13 January 2021.

\renewcommand{\arraystretch}{1.1}

\begin{table}[ht]
\caption{Totally real cubic discriminants \(d_L=f^2\cdot d\) in the range \(0<d_L<10^7\)}
\label{tbl:LlorenteQuer0}
\begin{center}
\begin{tabular}{|rl||r|rrrrr||rrrrrr|rrr|}
\hline
       &           & \multicolumn{6}{c||}{Multiplicity}    & \multicolumn{9}{c|}{Differential Principal Factorization} \\
 \(f\) & Condition & \(0\) & \(1\) & \(2\) & \(3\) & \(4\) & \(6\) & \(\alpha_1\) & \(\alpha_2\) & \(\alpha_3\) & \(\beta_1\) & \(\beta_2\) & \(\gamma\) & \(\delta_1\) & \(\delta_2\) & \(\varepsilon\) \\
\hline
 \(1\)            &    \(\varrho_3=0\) &\(2623325\)&      \(\) &   \(\) &    \(\) &   \(\) &  \(\) &  \(\) &  \(\) &  \(\) &   \(\) &    \(\) &   \(\) &    \(\)&    \(\) &    \(\) \\
\hline
 \(q\)            &  \(\equiv 2\,(3)\) &\(198952\) & \(88925\) &   \(\) &    \(\) &   \(\) &  \(\) &  \(\) &  \(\) &  \(\) &   \(\) &    \(\) &   \(\) &   \(\) &    \(\) &\(88925\) \\
 \(3\)            & \(d\equiv 3\,(9)\) & \(25596\) & \(11430\) &   \(\) &    \(\) &   \(\) &  \(\) &  \(\) &  \(\) &  \(\) &   \(\) &    \(\) &   \(\) &   \(\) &    \(\) &\(11430\) \\
 \(3\)            & \(d\equiv 6\,(9)\) & \(25563\) & \(11521\) &   \(\) &    \(\) &   \(\) &  \(\) &  \(\) &  \(\) &  \(\) &   \(\) &    \(\) &   \(\) &   \(\) &    \(\) &\(11521\) \\
 \(9\)            & \(d\equiv 6\,(9)\) &   \(947\) &  \(2947\) &   \(\) & \(308\) &   \(\) &  \(\) &  \(\) &  \(\) &  \(\) &   \(\) &    \(\) &   \(\) &   \(\) &    \(\) & \(3871\) \\
 \(9\)            & \(d\equiv 2\,(3)\) &  \(8669\) &  \(3846\) &   \(\) &    \(\) &   \(\) &  \(\) &  \(\) &  \(\) &  \(\) &   \(\) &    \(\) &   \(\) &   \(\) &    \(\) & \(3846\) \\
\hline
 \(9\)            & \(d\equiv 1\,(3)\) &  \(8594\) &  \(3860\) &   \(\) &    \(\) &   \(\) &  \(\) &  \(\) &  \(\) &  \(\) &   \(\) & \(888\) &   \(\) &   \(\) &\(2676\) &  \(296\) \\
 \(\ell\)         &  \(\equiv 1\,(3)\) & \(28591\) & \(11937\) &   \(\) &    \(\) &   \(\) &  \(\) &  \(\) &  \(\) &  \(\) &   \(\) &\(2768\) &   \(\) &   \(\) &\(8252\) &  \(917\) \\
\hline
 \(q_1q_2\)       &                    &  \(2706\) &  \(3092\) &\(429\) &    \(\) &   \(\) &  \(\) &  \(\) &  \(\) &  \(\) &   \(\) &    \(\) &\(3140\)&   \(\) &    \(\) &  \(810\) \\
 \(3q\)           & \(d\equiv 3\,(9)\) &  \(1811\) &  \(2003\) &\(305\) &    \(\) &   \(\) &  \(\) &  \(\) &  \(\) &  \(\) &   \(\) &    \(\) &\(2038\)&   \(\) &    \(\) &  \(575\) \\
 \(3q\)           & \(d\equiv 6\,(9)\) &  \(1826\) &  \(1973\) &\(318\) &    \(\) &   \(\) &  \(\) &  \(\) &  \(\) &  \(\) &   \(\) &    \(\) &\(1997\)&   \(\) &    \(\) &  \(612\) \\
 \(9q\)           & \(d\equiv 6\,(9)\) &    \(58\) &      \(\) &\(340\) &  \(77\) &   \(\) & \(1\) &  \(\) &  \(\) &  \(\) &   \(\) &    \(\) & \(727\)&   \(\) &    \(\) &  \(190\) \\
 \(9q\)           & \(d\equiv 2\,(3)\) &   \(599\) &   \(701\) & \(89\) &    \(\) &   \(\) &  \(\) &  \(\) &  \(\) &  \(\) &   \(\) &    \(\) & \(714\)&   \(\) &    \(\) &  \(165\) \\
\hline
 \(9q\)           & \(d\equiv 1\,(3)\) &   \(610\) &   \(686\) & \(84\) &    \(\) &   \(\) &  \(\) &  \(\) &  \(\) &  \(\) &   \(\) & \(691\) &  \(4\) &   \(\) &  \(34\) &  \(125\) \\
 \(q\ell\)        &                    &  \(1908\) &  \(2308\) &\(280\) &    \(\) &   \(\) &  \(\) &  \(\) &  \(\) &  \(\) &   \(\) &\(2273\) & \(68\) &   \(\) &  \(92\) &  \(435\) \\
 \(3\ell\)        & \(d\equiv 3\,(9)\) &   \(250\) &   \(307\) & \(28\) &    \(\) &   \(\) &  \(\) &  \(\) &  \(\) &  \(\) &   \(\) & \(312\) &  \(1\) &   \(\) &  \(11\) &   \(39\) \\
 \(3\ell\)        & \(d\equiv 6\,(9)\) &   \(254\) &   \(300\) & \(38\) &    \(\) &   \(\) &  \(\) &  \(\) &  \(\) &  \(\) &   \(\) & \(301\) &  \(3\) &   \(\) &   \(6\) &   \(66\) \\
 \(9\ell\)        & \(d\equiv 6\,(9)\) &     \(5\) &      \(\) & \(47\) &  \(14\) &   \(\) &  \(\) &  \(\) &  \(\) &  \(\) &   \(\) & \(110\) &   \(\) &   \(\) &    \(\) &   \(26\) \\
 \(9\ell\)        & \(d\equiv 2\,(3)\) &    \(89\) &   \(105\) &  \(6\) &    \(\) &   \(\) &  \(\) &  \(\) &  \(\) &  \(\) &   \(\) & \(105\) &   \(\) &   \(\) &    \(\) &   \(12\) \\
\hline
 \(9\ell\)        & \(d\equiv 1\,(3)\) &    \(72\) &    \(89\) &  \(3\) &    \(\) &   \(\) &  \(\) &  \(\) &  \(\) &\(72\) &   \(\) &  \(19\) &   \(\) &   \(\) &   \(4\) &     \(\) \\
 \(\ell_1\ell_2\) &                    &    \(60\) &    \(86\) &  \(2\) &    \(\) &   \(\) &  \(\) &  \(\) &  \(\) &\(50\) &   \(\) &  \(38\) &   \(\) &   \(\) &   \(1\) &    \(1\) \\
\hline
 \(q_1q_2q_3\)    &                    &     \(6\) &    \(19\) & \(12\) &    \(\) &   \(\) &  \(\) &  \(\) &  \(\) &  \(\) &   \(\) &    \(\) & \(43\) &   \(\) &    \(\) &     \(\) \\
 \(3q_1q_2\)      & \(d\equiv 3\,(9)\) &    \(14\) &    \(32\) & \(51\) &    \(\) &   \(\) &  \(\) &  \(\) &  \(\) &  \(\) &   \(\) &    \(\) &\(134\) &   \(\) &    \(\) &     \(\) \\
 \(3q_1q_2\)      & \(d\equiv 6\,(9)\) &    \(14\) &    \(25\) & \(40\) &    \(\) &   \(\) &  \(\) &  \(\) &  \(\) &  \(\) &   \(\) &    \(\) &\(105\) &   \(\) &    \(\) &     \(\) \\
 \(9q_1q_2\)      & \(d\equiv 6\,(9)\) &      \(\) &      \(\) &   \(\) &    \(\) &  \(5\) & \(1\) &  \(\) &  \(\) &  \(\) &   \(\) &    \(\) & \(26\) &   \(\) &    \(\) &     \(\) \\
 \(9q_1q_2\)      & \(d\equiv 2\,(3)\) &     \(4\) &    \(13\) & \(16\) &    \(\) &   \(\) &  \(\) &  \(\) &  \(\) &  \(\) &   \(\) &    \(\) & \(45\) &   \(\) &    \(\) &     \(\) \\
\hline
 \(9q_1q_2\)      & \(d\equiv 1\,(3)\) &     \(6\) &    \(11\) & \(14\) &    \(\) &   \(\) &  \(\) &  \(\) &  \(\) &  \(\) &   \(\) &   \(6\) & \(33\) &   \(\) &    \(\) &     \(\) \\
 \(q_1q_2\ell\)   &                    &    \(13\) &    \(35\) & \(44\) &    \(\) &   \(\) &  \(\) &  \(\) &  \(\) &  \(\) &   \(\) &  \(20\) &\(103\) &   \(\) &    \(\) &     \(\) \\
 \(3q\ell\)       & \(d\equiv 3\,(9)\) &     \(8\) &    \(27\) & \(34\) &    \(\) &   \(\) &  \(\) &  \(\) &  \(\) &  \(\) &   \(\) &  \(14\) & \(81\) &   \(\) &    \(\) &     \(\) \\
 \(3q\ell\)       & \(d\equiv 6\,(9)\) &     \(6\) &    \(20\) & \(24\) &    \(\) &  \(1\) &  \(\) &  \(\) &  \(\) &  \(\) &   \(\) &  \(25\) & \(44\) &   \(\) &    \(\) &    \(3\) \\
 \(9q\ell\)       & \(d\equiv 6\,(9)\) &      \(\) &      \(\) &   \(\) &   \(1\) &  \(3\) &  \(\) &  \(\) &  \(\) &  \(\) &   \(\) &   \(1\) & \(14\) &   \(\) &    \(\) &     \(\) \\
 \(9q\ell\)       & \(d\equiv 2\,(3)\) &     \(2\) &    \(11\) & \(12\) &    \(\) &   \(\) &  \(\) &  \(\) &  \(\) &  \(\) &   \(\) &   \(7\) & \(28\) &   \(\) &    \(\) &     \(\) \\
\hline
 \(9q\ell\)       & \(d\equiv 1\,(3)\) &      \(\) &    \(10\) & \(11\) &    \(\) &   \(\) &  \(\) &  \(\) &  \(\) &  \(\) &   \(\) &  \(29\) &  \(3\) &   \(\) &    \(\) &     \(\) \\
\(q\ell_1\ell_2\) &                    &      \(\) &     \(6\) &  \(3\) &    \(\) &   \(\) &  \(\) &  \(\) &  \(\) &  \(\) &   \(\) &  \(12\) &   \(\) &   \(\) &    \(\) &     \(\) \\
\(3\ell_1\ell_2\) & \(d\equiv 3\,(9)\) &      \(\) &     \(1\) &   \(\) &    \(\) &   \(\) &  \(\) &  \(\) &  \(\) &  \(\) &   \(\) &   \(1\) &   \(\) &   \(\) &    \(\) &     \(\) \\
\hline
 \(3q_1q_2\ell\)  &                    &      \(\) &      \(\) &  \(1\) &    \(\) &   \(\) &  \(\) &  \(\) &  \(\) &  \(\) &   \(\) &    \(\) &  \(2\) &   \(\) &    \(\) &     \(\) \\
\hline
\hline
                  & Subtotal           &      \(\) &\(146326\) &\(2231\)& \(400\) &   \(9\)& \(2\) &  \(\) &  \(\) &\(122\)&   \(\) &\(7620\) &\(9353\)&   \(\) &\(11076\)&\(123865\)\\
\hline
\end{tabular}
\end{center}
\end{table}

\newpage

\renewcommand{\arraystretch}{1.1}

\begin{table}[ht]
\caption{Table \ref{tbl:LlorenteQuer0} with \(0<d_L<10^7\) continued for \(\varrho_3\ge 1\)}
\label{tbl:LlorenteQuer1}
\begin{center}
\begin{tabular}{|rl||r|rrrrr||rrrrrr|rrr|}
\hline
       &           & \multicolumn{6}{c||}{Multiplicity}    & \multicolumn{9}{c|}{Differential Principal Factorization} \\
 \(f\) & Condition & \(0\) & \(1\) & \(2\) & \(3\) & \(4\) & \(6\) & \(\alpha_1\) & \(\alpha_2\) & \(\alpha_3\) & \(\beta_1\) & \(\beta_2\) & \(\gamma\) & \(\delta_1\) & \(\delta_2\) & \(\varepsilon\) \\
\hline
 \(1\)            &    \(\varrho_3=1\) &      \(\) &\(413458\) &   \(\) &    \(\) &   \(\) &  \(\) &  \(\) &  \(\) &  \(\) &   \(\) &   \(\) &    \(\) &\(413458\)&  \(\) &     \(\) \\
\hline
 \(q\)            &  \(\equiv 2\,(3)\) & \(38302\) &      \(\) &   \(\) &\(3239\) &   \(\) &  \(\) &  \(\) &  \(\) &  \(\) &\(2022\)&   \(\) &    \(\) &\(6958\)&    \(\) &  \(737\) \\
 \(3\)            & \(d\equiv 3\,(9)\) &  \(4798\) &      \(\) &   \(\) & \(375\) &   \(\) &  \(\) &  \(\) &  \(\) &  \(\) & \(199\)&   \(\) &    \(\) & \(857\)&    \(\) &   \(69\) \\
 \(3\)            & \(d\equiv 6\,(9)\) &  \(4760\) &      \(\) &   \(\) & \(359\) &   \(\) &  \(\) &  \(\) &  \(\) &  \(\) & \(223\)&   \(\) &    \(\) & \(773\)&    \(\) &   \(81\) \\
 \(9\)            & \(d\equiv 6\,(9)\) &   \(393\) &      \(\) &   \(\) &  \(99\) &   \(\) &  \(\) &  \(\) &  \(\) &  \(\) &  \(61\)&   \(\) &    \(\) & \(211\)&    \(\) &   \(25\) \\
 \(9\)            & \(d\equiv 2\,(3)\) &  \(1441\) &      \(\) &   \(\) & \(115\) &   \(\) &  \(\) &  \(\) &  \(\) &  \(\) &  \(83\)&   \(\) &    \(\) & \(241\)&    \(\) &   \(21\) \\
\hline
 \(9\)            & \(d\equiv 1\,(3)\) &  \(1489\) &      \(\) &   \(\) & \(124\) &   \(\) &  \(\) &  \(\) &\(27\) &  \(\) &  \(85\)&   \(\) &    \(\) & \(232\)&   \(7\) &   \(21\) \\
 \(\ell\)         &  \(\equiv 1\,(3)\) &  \(4470\) &      \(\) &   \(\) & \(386\) &   \(\) &  \(\) &  \(\) &\(95\) &  \(\) & \(230\)&  \(8\) &    \(\) & \(706\)&  \(33\) &   \(86\) \\
\hline
 \(q_1q_2\)       &                    &   \(534\) &      \(\) &   \(\) & \(115\) &   \(\) &  \(\) &  \(\) &  \(\) &  \(\) & \(278\)&   \(\) &  \(12\) &   \(\) &    \(\) &   \(55\) \\
 \(3q\)           & \(d\equiv 3\,(9)\) &   \(370\) &      \(\) &   \(\) &  \(78\) &   \(\) & \(1\) &  \(\) &  \(\) &  \(\) & \(187\)&   \(\) &  \(15\) &  \(3\) &    \(\) &   \(35\) \\
 \(3q\)           & \(d\equiv 6\,(9)\) &   \(399\) &      \(\) &   \(\) &  \(84\) &   \(\) & \(1\) &  \(\) &  \(\) &  \(\) & \(174\)&   \(\) &   \(9\) &  \(4\) &    \(\) &   \(71\) \\
 \(9q\)           & \(d\equiv 6\,(9)\) &    \(13\) &      \(\) &   \(\) &  \(25\) &   \(\) & \(1\) &  \(\) &  \(\) &  \(\) &  \(69\)&   \(\) &    \(\) &   \(\) &    \(\) &   \(12\) \\
 \(9q\)           & \(d\equiv 2\,(3)\) &   \(111\) &      \(\) &   \(\) &  \(20\) &   \(\) &  \(\) &  \(\) &  \(\) &  \(\) &  \(43\)&   \(\) &   \(5\) &   \(\) &    \(\) &   \(12\) \\
\hline
 \(9q\)           & \(d\equiv 1\,(3)\) &   \(117\) &      \(\) &   \(\) &  \(24\) &   \(\) &  \(\) &  \(\) & \(6\) &  \(\) &  \(57\)&  \(3\) &    \(\) &   \(\) &   \(3\) &    \(3\) \\
 \(q\ell\)        &                    &   \(364\) &      \(\) &   \(\) &  \(67\) &   \(\) &  \(\) &  \(\) & \(6\) &  \(\) & \(160\)&  \(6\) &   \(2\) &   \(\) &   \(6\) &   \(21\) \\
 \(3\ell\)        & \(d\equiv 3\,(9)\) &    \(44\) &      \(\) &   \(\) &   \(7\) &   \(\) &  \(\) &  \(\) & \(3\) &  \(\) &  \(17\)&  \(1\) &    \(\) &   \(\) &    \(\) &     \(\) \\
 \(3\ell\)        & \(d\equiv 6\,(9)\) &    \(36\) &      \(\) &   \(\) &  \(12\) &   \(\) &  \(\) &  \(\) &  \(\) &  \(\) &  \(24\)&   \(\) &    \(\) &   \(\) &   \(3\) &    \(9\) \\
 \(9\ell\)        & \(d\equiv 6\,(9)\) &      \(\) &      \(\) &   \(\) &   \(4\) &   \(\) &  \(\) &  \(\) &  \(\) &  \(\) &  \(12\)&   \(\) &    \(\) &   \(\) &    \(\) &     \(\) \\
 \(9\ell\)        & \(d\equiv 2\,(3)\) &    \(12\) &      \(\) &   \(\) &    \(\) &   \(\) &  \(\) &  \(\) &  \(\) &  \(\) &   \(\) &   \(\) &    \(\) &   \(\) &    \(\) &     \(\) \\
\hline
 \(9\ell\)        & \(d\equiv 1\,(3)\) &    \(12\) &      \(\) &   \(\) &   \(1\) &   \(\) &  \(\) &  \(\) & \(3\) &  \(\) &   \(\) &   \(\) &    \(\) &   \(\) &    \(\) &     \(\) \\
 \(\ell_1\ell_2\) &                    &     \(4\) &      \(\) &   \(\) &   \(1\) &   \(\) &  \(\) &  \(\) & \(2\) &  \(\) &   \(\) &  \(1\) &    \(\) &   \(\) &    \(\) &     \(\) \\
\hline
 \(q_1q_2q_3\)    &                    &     \(1\) &      \(\) &   \(\) &    \(\) &   \(\) &  \(\) &  \(\) &  \(\) &  \(\) &   \(\) &   \(\) &    \(\) &   \(\) &    \(\) &     \(\) \\
 \(3q_1q_2\)      & \(d\equiv 3\,(9)\) &     \(3\) &      \(\) &   \(\) &   \(1\) &   \(\) &  \(\) &  \(\) &  \(\) &  \(\) &   \(\) &   \(\) &   \(3\) &   \(\) &    \(\) &     \(\) \\
 \(3q_1q_2\)      & \(d\equiv 6\,(9)\) &     \(4\) &      \(\) &   \(\) &   \(4\) &   \(\) &  \(\) &  \(\) &  \(\) &  \(\) &   \(\) &   \(\) &  \(12\) &   \(\) &    \(\) &     \(\) \\
\hline
 \(9q_1q_2\)      & \(d\equiv 1\,(3)\) &      \(\) &      \(\) &   \(\) &   \(1\) &   \(\) &  \(\) &  \(\) &  \(\) &  \(\) &   \(\) &   \(\) &   \(3\) &   \(\) &    \(\) &     \(\) \\
 \(q_1q_2\ell\)   &                    &      \(\) &      \(\) &   \(\) &   \(2\) &   \(\) &  \(\) &  \(\) &  \(\) &  \(\) &   \(\) &   \(\) &   \(6\) &   \(\) &    \(\) &     \(\) \\
 \(3q\ell\)       & \(d\equiv 3\,(9)\) &     \(3\) &      \(\) &   \(\) &    \(\) &   \(\) &  \(\) &  \(\) &  \(\) &  \(\) &   \(\) &   \(\) &    \(\) &   \(\) &    \(\) &     \(\) \\
 \(3q\ell\)       & \(d\equiv 6\,(9)\) &     \(4\) &      \(\) &   \(\) &    \(\) &   \(\) &  \(\) &  \(\) &  \(\) &  \(\) &   \(\) &   \(\) &    \(\) &   \(\) &    \(\) &     \(\) \\
\hline
\hline
 \(1\)            &    \(\varrho_3=2\) &      \(\) &      \(\) &   \(\) &    \(\) &\(2870\)&  \(\) &\(7951\)& \(\) &  \(\) &   \(\) &   \(\) &    \(\) &\(3529\)&    \(\) &     \(\) \\
\hline
 \(q\)            &  \(\equiv 2\,(3)\) &   \(197\) &      \(\) &   \(\) &    \(\) &   \(\) &  \(\) &  \(\) &  \(\) &  \(\) &   \(\) &   \(\) &    \(\) &   \(\) &    \(\) &     \(\) \\
 \(3\)            & \(d\equiv 3\,(9)\) &    \(19\) &      \(\) &   \(\) &    \(\) &   \(\) &  \(\) &  \(\) &  \(\) &  \(\) &   \(\) &   \(\) &    \(\) &   \(\) &    \(\) &     \(\) \\
 \(3\)            & \(d\equiv 6\,(9)\) &    \(18\) &      \(\) &   \(\) &    \(\) &   \(\) &  \(\) &  \(\) &  \(\) &  \(\) &   \(\) &   \(\) &    \(\) &   \(\) &    \(\) &     \(\) \\
 \(9\)            & \(d\equiv 2\,(3)\) &     \(3\) &      \(\) &   \(\) &    \(\) &   \(\) &  \(\) &  \(\) &  \(\) &  \(\) &   \(\) &   \(\) &    \(\) &   \(\) &    \(\) &     \(\) \\
\hline
 \(9\)            & \(d\equiv 1\,(3)\) &     \(3\) &      \(\) &   \(\) &    \(\) &   \(\) &  \(\) &  \(\) &  \(\) &  \(\) &   \(\) &   \(\) &    \(\) &   \(\) &    \(\) &     \(\) \\
 \(\ell\)         &  \(\equiv 1\,(3)\) &     \(6\) &      \(\) &   \(\) &    \(\) &   \(\) &  \(\) &  \(\) &  \(\) &  \(\) &   \(\) &   \(\) &    \(\) &   \(\) &    \(\) &     \(\) \\
\hline
 \(3q\)           & \(d\equiv 3\,(9)\) &     \(3\) &      \(\) &   \(\) &    \(\) &   \(\) &  \(\) &  \(\) &  \(\) &  \(\) &    \(\)&   \(\) &    \(\) &   \(\) &    \(\) &     \(\) \\
 \(3q\)           & \(d\equiv 6\,(9)\) &     \(1\) &      \(\) &   \(\) &    \(\) &   \(\) &  \(\) &  \(\) &  \(\) &  \(\) &    \(\)&   \(\) &    \(\) &   \(\) &    \(\) &     \(\) \\
\hline
\hline
                  & Subtotal           &      \(\) &\(413458\) &  \(0\) &\(5143\) &\(2870\)& \(3\)&\(7951\)&\(142\)& \(0\) &\(3924\)& \(19\) & \(67\)&\(426972\)&  \(52\) & \(1258\) \\
\hline
\hline
                  & \textbf{Total}     &      \(\) &\(559784\) &\(2231\)&\(5543\) &\(2879\)&\(5\)&\(7951\)&\(142\)&\(122\)&\(3924\)&\(7639\)&\(9420\)&\(426972\)&\(11128\)&\(125123\)\\
\hline

\end{tabular}
\end{center}
\end{table}

\noindent
According to Table
\ref{tbl:LlorenteQuer1},
the total number of all \textit{non-cyclic} totally real cubic fields \(L\) with discriminants \(d<10^7\)
is \(\mathbf{592\,421}\).
Together with \(501\) \textit{cyclic} cubic fields in Table
\ref{tbl:LlorenteCyclic}
the number is \(\mathbf{592\,922}\), in perfect accordance with Belabas
\cite[p. 1231 and Tbl. 6.2, p. 1232]{Be1997},
one field less than in the table of Llorente and Quer
\cite{LlQu1988}
(the unknown needle in a gigantic hay stack).

\restoregeometry

\newpage

\noindent
We emphasize the difference between
the \textit{number of discriminants} (without multiplicities),
\[559784+2231+5543+2879+5=\mathbf{570\,442},\]
and the \textit{number of pairwise non-isomorphic fields} (including multiplicities in a weighted sum),
\[1\cdot 559784+2\cdot 2231+3\cdot 5543+4\cdot 2879+6\cdot 5=559784+4462+16629+11516+30=\mathbf{592\,421},\]
which is confirmed by adding the contributions to the \(9\) DPF types,
\(\alpha_1\), \(\alpha_2\), \(\alpha_3\), \(\beta_1\), \(\beta_2\), \(\gamma\), \(\delta_1\), \(\delta_2\), \(\varepsilon\),
\[7951+142+122+3924+7639+9420+426972+11128+125123=\mathbf{592\,421}.\]

\renewcommand{\arraystretch}{1.1}

\begin{table}[hb]
\caption{Cyclic cubic discriminants \(d_L=f^2\) in the range \(0<d_L<10^7\)}
\label{tbl:LlorenteCyclic}
\begin{center}
\begin{tabular}{|rl||rrr||r||rr|}
\hline
                        &                                 & \multicolumn{3}{c||}{M}   & DPF       &    \(f\) &         \(d_L\) \\
 \(f\)                  & Condition                       &   \(1\) &   \(2\) & \(4\) & \(\zeta\) &          &                 \\
\hline
 \(9\)                  & \(d=1\)                         &   \(1\) &    \(\) &  \(\) &     \(1\) &    \(9\) &          \(81\) \\
 \(\ell\)               & \(\equiv 1\,(\mathrm{mod}\,3)\) & \(216\) &    \(\) &  \(\) &   \(216\) &    \(7\) &          \(49\) \\
\hline
 \(9\ell\)              & \(d=1\)                         &    \(\) &  \(33\) &  \(\) &    \(66\) &   \(63\) &      \(3\,969\) \\
 \(\ell_1\ell_2\)       & \(\equiv 1\,(\mathrm{mod}\,3)\) &    \(\) &  \(93\) &  \(\) &   \(186\) &   \(91\) &      \(8\,281\) \\
\hline
 \(9\ell_1\ell_2\)      & \(d=1\)                         &    \(\) &    \(\) & \(6\) &    \(24\) &  \(819\) &    \(670\,761\) \\
 \(\ell_1\ell_2\ell_3\) & \(\equiv 1\,(\mathrm{mod}\,3)\) &    \(\) &    \(\) & \(2\) &     \(8\) & \(1729\) & \(2\,989\,441\) \\
\hline
                        & Summary                         & \(217\) & \(126\) & \(8\) &   \(501\) &          &                 \\
\hline
\end{tabular}
\end{center}
\end{table}

\noindent
According to Table
\ref{tbl:LlorenteCyclic},
the number of \textit{cyclic} cubic fields \(L\) with discriminant \(0<d_L<10^7\) is \(\mathbf{501}\),
with \(217\) arising from \textit{singlets} having conductors \(f\) with a single prime divisor,
\(252\) from \textit{doublets} having two prime divisors of the conductor \(f\),
and \(32\) from \textit{quartets} having three prime divisors of the conductor \(f\).
(M denotes the multiplicity.)

We point out that cyclic cubic fields
are rather contained in \textit{ray class fields} over \(\mathbb{Q}\)
than in ring class fields over real quadratic base fields.
The single possible DPF type \(\zeta\)
has nothing to do with the \(9\) DPF types
\(\alpha_1,\alpha_2,\alpha_3,\beta_1,\beta_2,\gamma,\delta_1,\delta_2,\varepsilon\)
of non-abelian totally real cubic fields in
\cite{Ma2019b}.


\subsection{Unramified Quartets}
\label{ss:Quartets}

\noindent
According to Theorem
\ref{thm:Unramified},
the \(413\,458\) unramified singlets \(N/K\) with conductor \(f=1\)
over quadratic base fields \(K\) with \(3\)-class rank \(\varrho=\varrho_3(K)=1\) 
form an overwhelming crowd of colorless, monotonous, and boring fields which share the common type \(\delta_1\).

In contrast, the \(\mathbf{2\,870}\) \textbf{unramified quartets} \(N/K\)
over quadratic fields \(K\) with \(\varrho=2\) show an interesting statistical distribution of types.
We consider the type \((\tau(L_1),\ldots,\tau(L_4))\) of a quartet \((L_1,\ldots,L_4)\)
as ordered lexicographically, regardless of permutations. Smallest \(d\) see Table
\ref{tbl:LlorenteQuartetsUnramified}.

\renewcommand{\arraystretch}{1.1}

\begin{table}[hb]
\caption{Types of unramified quartets in the range \(0<d_L<10^7\)}
\label{tbl:LlorenteQuartetsUnramified}
\begin{center}
\begin{tabular}{|c||c|r||r|}
\hline
 DPF Type                                  & Capitulation Number \(\nu(K)\)  &   Frequency &    \(d_L=d\) \\
 \((\tau(L_1),\ldots,\tau(L_4))\)          & (according to \cite{ChFt1980})  &             &              \\
\hline
 \((\alpha_1,\alpha_1,\alpha_1,\alpha_1)\) & \(4\)                           &     \(175\) &  \(62\,501\) \\
 \((\alpha_1,\alpha_1,\alpha_1,\delta_1)\) & \(3\)                           &    \(2391\) &  \(32\,009\) \\
 \((\alpha_1,\alpha_1,\delta_1,\delta_1)\) & \(2\)                           &       \(8\) & \(710\,652\) \\
 \((\alpha_1,\delta_1,\delta_1,\delta_1)\) & \(1\)                           &      \(62\) & \(534\,824\) \\
 \((\delta_1,\delta_1,\delta_1,\delta_1)\) & \(0\)                           &     \(234\) & \(214\,712\) \\
\hline
                                           &                          Total: &    \(2870\) &              \\
\hline
\end{tabular}
\end{center}
\end{table}

As known from
\cite{Ma2012}
and
\cite{Ma2014b},
the \(2391\) quartets of \textit{mixed} type \((\alpha_1,\alpha_1,\alpha_1,\delta_1)\) are extremely dominating
with a relative frequency of \(83.31\%\).
Moderate contributions are provided by the \(234\), resp. \(175\), quartets
of \textit{pure} type \((\delta_1,\delta_1,\delta_1,\delta_1)\), resp. \((\alpha_1,\alpha_1,\alpha_1,\alpha_1)\).
Quartets with mixed type \((\alpha_1,\delta_1,\delta_1,\delta_1)\) are rare with \(62\) hits,
and the \(8\) quartets with mixed type \((\alpha_1,\alpha_1,\delta_1,\delta_1)\) are almost negligible.
The reason for this behavior is well understood,
because the corresponding capitulation types \(\varkappa(K)=(\ker(T_{N_1/K}),\ldots,\ker(T_{N_4/K}))\) enforce certain
second \(3\)-class groups \(\mathrm{Gal}(\mathrm{F}_3^2(K)/K)\) of the quadratic base fields \(K\)
which can be realized easily for the quartets with high frequency, due to modest group orders,
but require huge groups in the case of rare quartets (see
\cite{Ma2013}).


\subsection{Other Multiplets}
\label{ss:Multiplets}

\noindent
According to Table
\ref{tbl:LlorenteQuer1},
the number \(2231\) of doublets,
resp. \(5543\) of triplets,
resp. \(2879\) of quartets,
resp. \(5\) of sextets,
agrees with the corresponding counters given in
\cite[Tbl. 2, p. 588]{LlQu1988}, resp.
\cite[Tbl. 3, p. 589]{LlQu1988}, resp.
\cite[Tbl. 4, p. 589]{LlQu1988}, resp.
\cite[p. 588 and Tbl. 5, p. 590]{LlQu1988},
in the paper by Llorente and Quer.
However, there are two misprints in the text below Tbl. 4 on page 589 of
\cite{LlQu1988},
where the authors intended to state that \(2870\) among the \(2879\) quartets
belong to real quadratic fields \(K\) with \(3\)-class rank \(\varrho=2\),
namely the \textit{unramified} quartets in our Table
\ref{tbl:LlorenteQuartetsUnramified}.
But the remaining \(9\) quartets are \textit{ramified}
over real quadratic fields \(K\) with \(3\)-class rank \(\varrho=0\)
and show up in our Table
\ref{tbl:LlorenteQuer0}.
They are analyzed in detail in the following example.


\begin{example}
\label{exm:LlorenteQuartetsRamified}
A common feature of all \(\mathbf{9}\) \textbf{ramified quartets} \((L_1,\ldots,L_4)\)
with discriminants in the range \(0<d_L<10^7\)
is the congruence class of the quadratic fundamental discriminant
\(d\equiv 6\,(\mathrm{mod}\,9)\)
which enables both, conductors with \(3\)-contribution \(v_3(f)=1\) and \(v_3(f)=2\).
The reason for their multiplicity in terms of \(3\)-defects \(\delta_3(f)\)
(co-dimensions of \(3\)-ring spaces \(V_3(f)\)) was discussed in
\cite[Supplements Section, Part 1.a, p. S55, and Part 2.d, pp. S57--S58]{Ma1992}.
Now we are able to present their differential principal factorizations in Table
\ref{tbl:LlorenteQuartetsRamified},
where the type of the conductor establishes the connection with Table
\ref{tbl:LlorenteQuer0}.
A generating polynomial for each member \(L\) of the quartets is given in
\cite[Tbl. 6, p. 591]{LlQu1988},
but we point out that the conductor in the caption of this table should be \(T=3^mT_0>1\) (our \(f\)),
and the discriminant in the table header should be \(D=3^{2m}T_0^2d\) (our \(d_L\)).
\end{example}


\renewcommand{\arraystretch}{1.0}

\begin{table}[ht]
\caption{Nine explicit ramified quartets in the range \(0<d_L<10^7\)}
\label{tbl:LlorenteQuartetsRamified}
\begin{center}
\begin{tabular}{|r|r||r|r|c||c|}
\hline
 No.   & \(d_L\)         & \(d\)    & \(f\)            & Kind of Conductor  & DPF Type                                         \\
       &                 &          &                  &                    & \((\tau(L_1),\ldots,\tau(L_4))\)                 \\
\hline
 \(1\) & \(1\,725\,300\) &  \(213\) &  \(90=9\cdot 2\cdot 5\) & \(9q_1q_2\) & \((\gamma,\gamma,\gamma,\gamma)\)                \\
 \(2\) & \(2\,238\,516\) &  \(141\) & \(126=9\cdot 2\cdot 7\) & \(9q\ell\)  & \((\gamma,\gamma,\gamma,\gamma)\)                \\
 \(3\) & \(2\,891\,700\) &  \(357\) &  \(90=9\cdot 2\cdot 5\) & \(9q_1q_2\) & \((\gamma,\gamma,\gamma,\gamma)\)                \\
 \(4\) & \(4\,641\,300\) &  \(573\) &  \(90=9\cdot 2\cdot 5\) & \(9q_1q_2\) & \((\gamma,\gamma,\gamma,\gamma)\)                \\
 \(5\) & \(6\,810\,804\) &  \(429\) & \(126=9\cdot 2\cdot 7\) & \(9q\ell\)  & \((\gamma,\gamma,\gamma,\gamma)\)                \\
 \(6\) & \(7\,557\,300\) &  \(933\) &  \(90=9\cdot 2\cdot 5\) & \(9q_1q_2\) & \((\gamma,\gamma,\gamma,\gamma)\)                \\
 \(7\) & \(7\,953\,876\) &  \(501\) & \(126=9\cdot 2\cdot 7\) & \(9q\ell\)  & \((\gamma,\gamma,\gamma,\gamma)\)                \\
 \(8\) & \(8\,250\,228\) & \(4677\) &  \(42=3\cdot 2\cdot 7\) & \(3q\ell\)  & \((\gamma,\varepsilon,\varepsilon,\varepsilon)\) \\
 \(9\) & \(8\,723\,700\) & \(1077\) &  \(90=9\cdot 2\cdot 5\) & \(9q_1q_2\) & \((\gamma,\gamma,\gamma,\gamma)\)                \\
\hline
\end{tabular}
\end{center}
\end{table}


\begin{example}
\label{exm:LlorenteSextets}
A particular highlight of the range \(0<d_L<10^7\)
is the occurrence of \(\mathbf{5}\) \textbf{sextets},
which did not show up in smaller tables.
The reason for their multiplicity in terms of \(3\)-defects \(\delta_3(f)\)
(co-dimensions of \(3\)-ring spaces \(V_3(f)\) modulo \(f\) in the \(3\)-Selmer space \(V_3\)) was discussed in
\cite[Supplements Section, Part 1.c, p. S56, Part 2.b, p. S57, Part 2.d, pp. S57--S58, and Part 2.f, p. S58]{Ma1992}.
Now we are able to present their differential principal factorizations in Table
\ref{tbl:LlorenteSextets}.
The leading two sextets are \textit{mixed}, and
the trailing three sextets are \textit{pure}.
The constitution of the sextets is very heterogeneous:
although four of the quadratic fundamental discriminants \(d\)
admit the \textit{irregular} contribution \(9\) to the conductor \(f\)
only three conductors are actually divisible by \(9\),
but they differ either by the \(3\)-rank \(\varrho\)
or by the kind of the conductor.
A generating polynomial for each member \(L\) of the sextets is given in
\cite[Tbl. 5, p. 590]{LlQu1988},
but again we point out that the discriminant in the table header should be \(D=3^{2m}T_0^2d\) (our \(d_L=f^2\cdot d\)).
Types for \(\varrho=0\) are more simple.
\end{example}


\renewcommand{\arraystretch}{1.0}

\begin{table}[hb]
\caption{Five explicit sextets in the range \(0<d_L<10^7\)}
\label{tbl:LlorenteSextets}
\begin{center}
\begin{tabular}{|r|r||r|c|r|c||c|}
\hline
 No.   & \(d_L\)         & \(d\)                     & \(\varrho\) & \(f\)                  & Kind of Conductor  & DPF Type                                                                      \\
       &                 &                           &             &                        &                    & \((\tau(L_1),\ldots,\tau(L_6))\)                                              \\
\hline
 \(1\) & \(3\,054\,132\) &  \(84\,837\equiv 3\,(9)\) & \(1\)       &  \(6=3\cdot 2\)        & \(3q\)             & \((\beta_1,\beta_1,\delta_1,\delta_1,\delta_1,\varepsilon)\)                  \\
 \(2\) & \(4\,735\,467\) & \(131\,541\equiv 6\,(9)\) & \(1\)       &  \(6=3\cdot 2\)        & \(3q\)             & \((\beta_1,\delta_1,\delta_1,\delta_1,\delta_1,\varepsilon)\)                 \\
 \(3\) & \(5\,807\,700\) &      \(717\equiv 6\,(9)\) & \(0\)       & \(90=9\cdot 2\cdot 5\) & \(9q_1q_2\)        & \((\gamma,\gamma,\gamma,\gamma,\gamma,\gamma)\)                               \\
 \(4\) & \(6\,367\,572\) &  \(19\,653\equiv 6\,(9)\) & \(1\)       & \(18=9\cdot 2\)        & \(9q\)             & \((\beta_1,\beta_1,\beta_1,\beta_1,\beta_1,\beta_1)\)                         \\
 \(5\) & \(9\,796\,788\) &  \(30\,237\equiv 6\,(9)\) & \(0\)       & \(18=9\cdot 2\)        & \(9q\)             & \((\varepsilon,\varepsilon,\varepsilon,\varepsilon,\varepsilon,\varepsilon)\) \\
\hline
\end{tabular}
\end{center}
\end{table}


\begin{example}
\label{exm:LlorenteDoublets}
We split the \(2231\) \textbf{doublets}
in the range \(0<d_L<10^7\)
according to the shape of \(f\).
\end{example}


\noindent
\(\bullet\)
In Table
\ref{tbl:LlorenteDoubletsTwo},
we begin with two non-split prime divisors of \(f\),
that is, we accumulate the results for
\(f=q_1q_2\),
\(f=3q\) with \(d\equiv\pm 3\,(\mathrm{mod}\,9)\), and
\(f=9q\) with \(d\equiv 2\,(\mathrm{mod}\,3)\).
The DPF type \((\varepsilon,\varepsilon)\) is highly dominating over
\((\gamma,\varepsilon)\) and \((\gamma,\gamma)\).
Here and in the sequel,
the given paradigms for \(d_L\) are \textit{not} necessarily minimal.
Note the constitution \(1141=429+305+318+89\) of the total frequency.

\renewcommand{\arraystretch}{1.1}

\begin{table}[ht]
\caption{Types of doublets with two non-split prime divisors of \textit{regular} \(f\)}
\label{tbl:LlorenteDoubletsTwo}
\begin{center}
\begin{tabular}{|c||r||rrr|}
\hline
 DPF Type                         & Frequency &    \(d\) &  \(f\) &      \(d_L\) \\
 \((\tau(L_1),\tau(L_2))\)        &           &          &        &              \\
\hline
 \((\gamma,\gamma)\)              &    \(40\) &   \(33\) & \(45\) &  \(66\,825\) \\
 \((\gamma,\varepsilon)\)         &    \(40\) & \(9973\) & \(10\) & \(997\,300\) \\
 \((\varepsilon,\varepsilon)\)    &  \(1061\) &  \(373\) & \(10\) &  \(37\,300\) \\
\hline
                           Total: &  \(1141\) &          &        &              \\
\hline
\end{tabular}
\end{center}
\end{table}


\noindent
\(\bullet\)
The irregular situation \(f=9q\) with \(d\equiv 6\,(\mathrm{mod}\,9)\) in Table
\ref{tbl:LlorenteDoubletsTwoIrr}
shows a reversal of tendencies.
DPF type \((\gamma,\gamma)\) is dominating,
\((\varepsilon,\varepsilon)\) remains moderate,
mixed type
\((\gamma,\varepsilon)\) is almost negligible.

\renewcommand{\arraystretch}{1.1}

\begin{table}[ht]
\caption{Types of doublets with two non-split prime divisors of \textit{irregular} \(f\)}
\label{tbl:LlorenteDoubletsTwoIrr}
\begin{center}
\begin{tabular}{|c||r||rrr|}
\hline
 DPF Type                         & Frequency &    \(d\) &  \(f\) &         \(d_L\) \\
 \((\tau(L_1),\tau(L_2))\)        &           &          &        &                 \\
\hline
 \((\gamma,\gamma)\)              &   \(245\) &  \(213\) & \(18\) &     \(69\,012\) \\
 \((\gamma,\varepsilon)\)         &     \(6\) & \(9213\) & \(18\) & \(2\,985\,012\) \\
 \((\varepsilon,\varepsilon)\)    &    \(89\) &  \(141\) & \(18\) &     \(45\,684\) \\
\hline
                           Total: &   \(340\) &          &        &                 \\
\hline
\end{tabular}
\end{center}
\end{table}


\noindent
\(\bullet\)
Table
\ref{tbl:LlorenteDoubletsSplit}
reveals that, for
\(f=q\ell\),
\(f=3\ell\) with \(d\equiv\pm 3\,(\mathrm{mod}\,9)\),
\(f=9\ell\) with \(d\equiv 2\,(\mathrm{mod}\,3)\), and
\(f=9q\) with \(d\equiv 1\,(\mathrm{mod}\,3)\),
DPF type \((\varepsilon,\varepsilon)\) prevails,
followed by \((\delta_2,\delta_2)\).

\renewcommand{\arraystretch}{1.1}

\begin{table}[ht]
\caption{Types of doublets with a split prime divisor of \textit{regular} \(f\)}
\label{tbl:LlorenteDoubletsSplit}
\begin{center}
\begin{tabular}{|c||r||rrr|}
\hline
 DPF Type                         & Frequency &      \(d\) &  \(f\) &         \(d_L\) \\
 \((\tau(L_1),\tau(L_2))\)        &           &            &        &                 \\
\hline
 \((\beta_2,\beta_2)\)            &    \(14\) &\(23\,717\) & \(14\) & \(4\,648\,532\) \\
 \((\beta_2,\delta_2)\)           &     \(1\) & \(5\,061\) & \(39\) & \(7\,697\,781\) \\
 \((\beta_2,\varepsilon)\)        &    \(14\) & \(7\,589\) & \(14\) & \(1\,487\,444\) \\
 \((\gamma,\varepsilon)\)         &     \(9\) & \(1\,192\) & \(65\) & \(5\,036\,200\) \\
 \((\delta_2,\delta_2)\)          &    \(71\) & \(4\,813\) & \(14\) &    \(943\,348\) \\
 \((\varepsilon,\varepsilon)\)    &   \(327\) &    \(197\) & \(14\) &     \(38\,612\) \\
\hline
                           Total: &   \(436\) &            &        &                 \\
\hline
\end{tabular}
\end{center}
\end{table}


\noindent
\(\bullet\)
Again, the irregular situation \(f=9\ell\) with \(d\equiv 6\,(\mathrm{mod}\,9)\) in Table
\ref{tbl:LlorenteDoubletsSplitIrr}
shows a reversal of tendencies.
DPF type \((\beta_2,\beta_2)\) dominates over \((\varepsilon,\varepsilon)\).

\renewcommand{\arraystretch}{1.1}

\begin{table}[ht]
\caption{Types of doublets with a split prime divisor of \textit{irregular} \(f\)}
\label{tbl:LlorenteDoubletsSplitIrr}
\begin{center}
\begin{tabular}{|c||r||rrr|}
\hline
 DPF Type                         & Frequency &    \(d\) &  \(f\) &       \(d_L\) \\
 \((\tau(L_1),\tau(L_2))\)        &           &          &        &               \\
\hline
 \((\beta_2,\beta_2)\)            &    \(34\) &   \(60\) & \(63\) &  \(238\,140\) \\
 \((\varepsilon,\varepsilon)\)    &    \(13\) &  \(204\) & \(63\) &  \(809\,676\) \\
\hline
                           Total: &    \(47\) &          &        &               \\
\hline
\end{tabular}
\end{center}
\end{table}

\newpage

\noindent
\(\bullet\)
In the case of three non-split prime divisors of \(f\), i.e.,
\(f=q_1q_2q_3\) or
\(f=3q_1q_2\) with \(d\equiv\pm 3\,(\mathrm{mod}\,9)\) or
\(f=9q_1q_2\) with \(d\equiv 2\,(\mathrm{mod}\,3)\),
no table is required,
since all \(119\) \(=12+51+40+16\) occurrences are of type \((\gamma,\gamma)\),
e.g. \(d=93\), \(f=30\), \(d_L=83\,700\).


\begin{example}
\label{exm:LlorenteTriplets}
We split the \(5543\) \textbf{triplets}
in the range \(0<d_L<10^7\)
according to the shape of \(f\).
\end{example}


\noindent
\(\bullet\)
Triplets are usually due to elevated \(3\)-class rank \(\varrho\ge 1\)
of the real quadratic field \(K\).
However, the simplest case of triplets with \(\varrho=0\)
arises for the \textit{irregular} prime power conductor \(f=9\), \(d\equiv 6\,(\mathrm{mod}\,9)\).
The minimal occurrence is \(d=717\), \(d_L=58\,077\).
Each of the \(308\) triplets is embedded in a \textit{hetero}geneous quartet
\(\mathrm{Inv}(K_9)=\lbrack (\varepsilon),(\varepsilon,\varepsilon,\varepsilon)\rbrack\).


\noindent
\(\bullet\)
There are \(77\) cases of triplets with irregular conductors
\(f=9q\), \(d\equiv 6\,(\mathrm{mod}\,9)\), with \(\varrho=0\).
They are all of pure type \((\gamma,\gamma,\gamma)\), for instance
\(d=69\), \(q=2\), \(f=18\), \(d_L=22\,356\).


\noindent
\(\bullet\)
For the irregular case \(f=9\ell\), \(d\equiv 6\,(\mathrm{mod}\,9)\), with \(\varrho=0\),
all \(14\) occurrences are of type \((\beta_2,\beta_2,\beta_2)\),
for instance \(d=177\), \(\ell=7\), \(d_L=702\,513\).
There always exists an associated singlet of type \((\varepsilon)\)
with conductor \(f=3\), that is,
the the triplet and the singlet are embedded in a \textit{hetero}geneous quartet
\(\mathrm{Inv}(K_{9\ell})=\lbrack\emptyset,(\varepsilon),\emptyset,\emptyset,\emptyset,(\beta_2,\beta_2,\beta_2)\rbrack\)
corresponding to the divisors \((1,3,9,\ell,3\ell,9\ell)\) of \(f\).


\noindent
\(\bullet\)
A unique example of \(f=9q\ell\), \(d\equiv 6\,(\mathrm{mod}\,9)\), with \(\varrho=0\),
is given by \(d=69\), \(q=2\), \(\ell=13\), \(f=234\), \(d_L=3\,778\,164\).
It is a triplet of mixed type \((\beta_2,\gamma,\gamma)\).


\noindent
\(\bullet\)
For \(\varrho=1\) and \textit{non-critical split} \(f=\ell\equiv 1\,(\mathrm{mod}\,3)\),
the mixed DPF type \((\beta_1,\delta_1,\delta_1)\) prevails,
followed by mixed type \((\delta_1,\delta_1,\varepsilon)\).
Mixed types have only two distinct components. See Table
\ref{tbl:LlorenteTripletsSplit}.

\renewcommand{\arraystretch}{1.1}

\begin{table}[ht]
\caption{Types of triplets with a \textit{non-critical split} prime conductor \(f=\ell\)}
\label{tbl:LlorenteTripletsSplit}
\begin{center}
\begin{tabular}{|c||r||rrr|}
\hline
 DPF Type                                 & Frequency &       \(d\) &  \(f\) &         \(d_L\) \\
 \((\tau(L_1),\tau(L_2),\tau(L_3))\)      &           &             &        &                 \\
\hline
 \((\alpha_2,\alpha_2,\alpha_2)\)         &    \(10\) & \(32\,204\) &  \(7\) & \(1\,577\,996\) \\
 \((\alpha_2,\alpha_2,\delta_1)\)         &     \(9\) &  \(2\,677\) & \(19\) &    \(966\,397\) \\
 \((\alpha_2,\alpha_2,\delta_2)\)         &    \(23\) &  \(9\,749\) & \(13\) & \(1\,647\,581\) \\
 \((\alpha_2,\delta_2,\delta_2)\)         &     \(1\) &  \(5\,477\) & \(37\) & \(7\,498\,013\) \\
 \((\beta_1,\beta_2,\beta_2)\)            &     \(4\) &  \(7\,244\) & \(19\) & \(2\,615\,084\) \\
 \((\beta_1,\delta_1,\delta_1)\)          &   \(226\) &  \(1\,765\) &  \(7\) &     \(86\,485\) \\
 \((\delta_1,\delta_1,\delta_1)\)         &    \(23\) & \(13\,688\) & \(13\) & \(2\,313\,272\) \\
 \((\delta_1,\delta_1,\delta_2)\)         &     \(1\) & \(30\,553\) & \(13\) & \(5\,163\,457\) \\
 \((\delta_1,\delta_1,\varepsilon)\)      &    \(86\) &  \(3\,873\) &  \(7\) &    \(189\,777\) \\
 \((\delta_1,\delta_2,\delta_2)\)         &     \(2\) & \(44\,641\) &  \(7\) & \(2\,187\,409\) \\
 \((\delta_2,\delta_2,\delta_2)\)         &     \(1\) & \(54\,469\) &  \(7\) & \(2\,668\,981\) \\
\hline
                                   Total: &   \(386\) &             &        &                 \\
\hline
\end{tabular}
\end{center}
\end{table}


\noindent
\(\bullet\)
For \(\varrho=1\) and \textit{critical split} \(f=9\), \(d\equiv 1\,(\mathrm{mod}\,3)\),
again the mixed DPF type \((\beta_1,\delta_1,\delta_1)\) prevails,
followed by mixed type \((\delta_1,\delta_1,\varepsilon)\).
Here, all examples have minimal discriminant \(d_L\). See Table
\ref{tbl:LlorenteTripletsSplitCrit}.

\renewcommand{\arraystretch}{1.1}

\begin{table}[ht]
\caption{Types of triplets with \textit{critical split} prime power conductor \(f=9\)}
\label{tbl:LlorenteTripletsSplitCrit}
\begin{center}
\begin{tabular}{|c||r||rrr|}
\hline
 DPF Type                                 & Frequency &       \(d\) & \(f\) &         \(d_L\) \\
 \((\tau(L_1),\tau(L_2),\tau(L_3))\)      &           &             &       &                 \\
\hline
 \((\alpha_2,\alpha_2,\alpha_2)\)         &     \(4\) & \(14\,197\) & \(9\) & \(1\,149\,957\) \\
 \((\alpha_2,\alpha_2,\delta_1)\)         &     \(2\) & \(31\,069\) & \(9\) & \(2\,516\,589\) \\
 \((\alpha_2,\alpha_2,\delta_2)\)         &     \(5\) & \(15\,529\) & \(9\) & \(1\,257\,849\) \\
 \((\alpha_2,\delta_2,\delta_2)\)         &     \(1\) & \(30\,904\) & \(9\) & \(2\,503\,224\) \\
 \((\beta_1,\delta_1,\delta_1)\)          &    \(85\) &  \(2\,917\) & \(9\) &    \(236\,277\) \\
 \((\delta_1,\delta_1,\delta_1)\)         &     \(6\) & \(13\,861\) & \(9\) & \(1\,122\,741\) \\
 \((\delta_1,\delta_1,\varepsilon)\)      &    \(21\) & \(15\,733\) & \(9\) & \(1\,274\,373\) \\
\hline
                                   Total: &   \(124\) &             &       &                 \\
\hline
\end{tabular}
\end{center}
\end{table}


\noindent
\(\bullet\)
For \(\varrho=1\) and \textit{non-split} \(f=q\equiv 2\,(\mathrm{mod}\,3)\)
or \(f=3\) or \(f=9\), \(d\equiv 2\,(\mathrm{mod}\,3)\),
the pure DPF type \((\delta_1,\delta_1,\delta_1)\) is dominating,
followed by the mixed type \((\beta_1,\beta_1,\varepsilon)\).
Mixed type \((\beta_1,\delta_1,\varepsilon)\) with three distinct components is very rare. See Table
\ref{tbl:LlorenteTripletsNonSplit},
where \(4088=3239+359+375+115\).

\renewcommand{\arraystretch}{1.0}

\begin{table}[ht]
\caption{Types of triplets with a \textit{non-split} prime (power) conductor}
\label{tbl:LlorenteTripletsNonSplit}
\begin{center}
\begin{tabular}{|c||r||rrr|}
\hline
 DPF Type                                  & Frequency &        \(d\) &  \(f\) &         \(d_L\) \\
 \((\tau(L_1),\tau(L_2),\tau(L_3))\)       &           &              &        &                 \\
\hline
 \((\beta_1,\beta_1,\beta_1)\)             &   \(304\) &  \(55\,885\) &  \(2\) &    \(223\,540\) \\
 \((\beta_1,\beta_1,\delta_1)\)            &   \(160\) &  \(30\,965\) &  \(2\) &    \(123\,860\) \\
 \((\beta_1,\beta_1,\varepsilon)\)         &   \(640\) &  \(14\,397\) &  \(2\) &     \(57\,588\) \\
 \((\beta_1,\delta_1,\varepsilon)\)        &     \(5\) & \(417\,077\) &  \(2\) & \(1\,668\,308\) \\
 \((\beta_1,\varepsilon,\varepsilon)\)     &    \(10\) & \(492\,117\) &  \(2\) & \(1\,968\,468\) \\
 \((\delta_1,\delta_1,\delta_1)\)          &  \(2869\) &   \(7\,053\) &  \(2\) &     \(28\,212\) \\
 \((\delta_1,\delta_1,\varepsilon)\)       &    \(11\) & \(486\,461\) &  \(2\) & \(1\,945\,844\) \\
 \((\delta_1,\varepsilon,\varepsilon)\)    &    \(35\) & \(192\,245\) &  \(2\) &    \(768\,980\) \\
 \((\varepsilon,\varepsilon,\varepsilon)\) &    \(54\) & \(197\,445\) &  \(2\) &    \(789\,780\) \\
\hline
                                    Total: &  \(4088\) &              &        &                 \\
\hline
\end{tabular}
\end{center}
\end{table}


\noindent
\(\bullet\)
Table
\ref{tbl:LlorenteTripletsTwoNonSplit}
gives the distribution of DPF types for
\(f=q_1q_2\), \(f=3q\), and \(f=9q\), \(d\equiv 2\,(\mathrm{mod}\,3)\),
where \(297=115+78+84+20\).
Pure DPF type \((\beta_1,\beta_1,\beta_1)\) prevails,
followed by pure type \((\varepsilon,\varepsilon,\varepsilon)\).

\renewcommand{\arraystretch}{1.0}

\begin{table}[ht]
\caption{Types of triplets with \textit{two non-split} ramified primes}
\label{tbl:LlorenteTripletsTwoNonSplit}
\begin{center}
\begin{tabular}{|c||r||rrr|}
\hline
 DPF Type                                  & Frequency &        \(d\) &  \(f\) &         \(d_L\) \\
 \((\tau(L_1),\tau(L_2),\tau(L_3))\)       &           &              &        &                 \\
\hline
 \((\beta_1,\beta_1,\beta_1)\)             &   \(221\) &   \(3\,173\) & \(10\) &    \(317\,300\) \\
 \((\beta_1,\beta_1,\gamma)\)              &     \(5\) &  \(63\,917\) & \(10\) & \(6\,391\,700\) \\
 \((\beta_1,\gamma,\gamma)\)               &     \(6\) &  \(82\,397\) & \(10\) & \(8\,239\,700\) \\
 \((\gamma,\gamma,\gamma)\)                &     \(6\) &   \(9\,413\) & \(22\) & \(4\,555\,892\) \\
 \((\gamma,\gamma,\varepsilon)\)           &     \(3\) &  \(64\,677\) & \(10\) & \(6\,467\,700\) \\
 \((\varepsilon,\varepsilon,\varepsilon)\) &    \(56\) &   \(9\,293\) & \(10\) &    \(929\,300\) \\
\hline
                                    Total: &   \(297\) &              &        &                 \\
\hline
\end{tabular}
\end{center}
\end{table}


\noindent
\(\bullet\)
Table
\ref{tbl:LlorenteTripletsNonSplitAndSplit}
shows the triplets with \(f=q\ell\), \(f=3\ell\), and \(f=9q\), \(d\equiv 1\,(\mathrm{mod}\,3)\).
(There are no hits for \(f=9\ell\), \(d\equiv 2\,(\mathrm{mod}\,3)\).)
The pure DPF type \((\beta_1,\beta_1,\beta_1)\) prevails,
followed by pure type \((\varepsilon,\varepsilon,\varepsilon)\).
Here, \(110=67+7+12+24\).

\renewcommand{\arraystretch}{1.0}

\begin{table}[ht]
\caption{Types of triplets with \textit{non-split and split} ramified primes}
\label{tbl:LlorenteTripletsNonSplitAndSplit}
\begin{center}
\begin{tabular}{|c||r||rrr|}
\hline
 DPF Type                                  & Frequency &        \(d\) &  \(f\) &         \(d_L\) \\
 \((\tau(L_1),\tau(L_2),\tau(L_3))\)       &           &              &        &                 \\
\hline
 \((\alpha_2,\alpha_2,\alpha_2)\)          &     \(5\) &   \(6\,997\) & \(14\) & \(1\,371\,412\) \\
 \((\beta_1,\beta_1,\beta_1)\)             &    \(83\) &   \(1\,101\) & \(14\) &    \(215\,796\) \\
 \((\beta_1,\beta_1,\beta_2)\)             &     \(1\) &  \(21\,324\) & \(21\) & \(9\,403\,884\) \\
 \((\beta_1,\beta_2,\beta_2)\)             &     \(3\) &  \(29\,317\) & \(14\) & \(5\,746\,132\) \\
 \((\beta_2,\beta_2,\beta_2)\)             &     \(1\) &  \(18\,661\) & \(18\) & \(6\,046\,164\) \\
 \((\beta_1,\beta_1,\gamma)\)              &     \(2\) &      \(469\) & \(62\) & \(1\,802\,836\) \\
 \((\delta_2,\delta_2,\delta_2)\)          &     \(4\) &  \(24\,621\) & \(14\) & \(4\,825\,716\) \\
 \((\varepsilon,\varepsilon,\varepsilon)\) &    \(11\) &  \(10\,733\) & \(14\) & \(2\,103\,668\) \\
\hline
                                    Total: &   \(110\) &              &        &                 \\
\hline
\end{tabular}
\end{center}
\end{table}


\noindent
\(\bullet\)
There are only two triplets with \textit{two split} ramified primes: \\
mixed type \((\alpha_2,\alpha_2,\beta_2)\) for \(f=\ell_1\ell_2\)
(\(d=940\), \(f=91\), \(d_L=7\,784\,140\)), and \\
pure type \((\alpha_2,\alpha_2,\alpha_2)\) for \(f=9\ell\), \(d\equiv 1\,(\mathrm{mod}\,3)\)
(\(d=2\,101\), \(f=63\), \(d_L=8\,338\,869\)).

\newpage

\noindent
We conclude this section on multiplets
with information on \textit{singlets}.

\begin{theorem}
\label{thm:Singlets}
(Ramified and unramified \textbf{singlets})
\begin{enumerate}
\item
A \textbf{ramified} singlet (with conductor \(f>1\))
can only be of type
\((\alpha_3)\), \((\beta_2)\), \((\gamma)\), \((\delta_2)\), \((\varepsilon)\).
\item
An \textbf{unramified} singlet (with conductor \(f=1\))
must be of type \((\delta_1)\).
\end{enumerate}
\end{theorem}

\begin{proof}
According to the fundamental inequalities in Corollary
\ref{cor:Estimates}
and the fundamental equation in Corollary
\ref{cor:Trichotomy},
we have:
\begin{enumerate}
\item
For \(f>1\), the multiplicity formula shows that \(3^\varrho\) divides \(m\)
\cite[Thm. 3.2, p. 2215, Thm. 3.3--3.4, p. 2217, and Thm. 4.1--4.2, p. 2224--2225]{Ma2014}.
Thus, a singlet can only occur for \(\varrho=0\),
and this implies \(C=0\), i.e. no capitulation can happen.
By Theorem
\ref{thm:MainReal},
we conclude
\(\tau(L)\notin\lbrace\alpha_1,\alpha_2,\beta_1,\delta_1\rbrace\),
and consequently
\(\tau(L)\in\lbrace\alpha_3,\beta_2,\gamma,\delta_2,\varepsilon\rbrace\).
\item
For \(f=1\), we have the multiplicity formula
\(m=(3^\varrho-1)/2\)
\cite[Thm. 3.1, p. 2214]{Ma2014}.
Therefore, a singlet with \(m=1\) occurs for \(\varrho=1\).
On the other hand, \(f=1\) implies \(t=\) \(s=0\),
and thus \(A=R=0\).
The fundamental equation degenerates to \(U+1=C\),
where \(\varrho=1\) implies the bound \(C\le 1\).
Thus, \(C=1\) and \(U=0\),
that is the unique type \(\delta_1\). \qedhere
\end{enumerate}
\end{proof}

\begin{example}
\label{exm:Singlets}
Indeed, singlets of all the types in Theorem
\ref{thm:Singlets}
actually do occur.
Their minimal discriminants \(d_L\)
are given in Table
\ref{tbl:Singlets}.
\end{example}

\renewcommand{\arraystretch}{1.1}

\begin{table}[ht]
\caption{Smallest occurrences of various singlets}
\label{tbl:Singlets}
\begin{center}
\begin{tabular}{|c||rrr|}
\hline
 DPF Type          &    \(d\) &  \(f\) &      \(d_L\) \\
 \((\tau(L))\)     &          &        &              \\
\hline
 \((\varepsilon)\) &   \(37\) &  \(2\) &      \(148\) \\
 \((\delta_1)\)    &  \(229\) &  \(1\) &      \(229\) \\
 \((\gamma)\)      &   \(21\) &  \(6\) &      \(756\) \\
 \((\delta_2)\)    &   \(53\) &  \(7\) &   \(2\,597\) \\
 \((\beta_2)\)     &   \(29\) & \(14\) &   \(5\,684\) \\
 \((\alpha_3)\)    &   \(37\) & \(63\) & \(146\,853\) \\
\hline
\end{tabular}
\end{center}
\end{table}


\noindent
Concerning the frequency of singlets for \(0<d_L<10^7\),
the first row in Table
\ref{tbl:LlorenteQuer1}
proves that unramified singlets \((\delta_1)\) form
the definite hichamp \(413\,458\) among all contributions.
The last row (Subtotal) in Table 
\ref{tbl:LlorenteQuer0}
illuminates the second extreme contribution \(146\,326\)
by all the other ramified singlets
\((\alpha_3)\), \((\beta_2)\), \((\gamma)\), \((\delta_2)\),
and clearly dominating \((\varepsilon)\).


Another interesting observation is enabled
by the rows with regular conductors \(f\) divisible by exactly two primes,
i.e. \(t=2\), in Table
\ref{tbl:LlorenteQuer0}.
It appears that, under certain conditions,
\textit{non-split extensions} \(N/K\) with \(U_K=N_{N/K}(U_N)\) in the sense of Remark
\ref{rmk:RealDPFTypes}
are forbidden.
Generalizing a proof for singlets of type \((\gamma)\)
by Moser at the top of p. 74 in
\cite{Mo1979},
we partition the case \(t=2\),
according to the number \(0\le s\le 2\) of prime divisors of \(f\)
which split in the real quadratic field \(K\).
The crucial assumption \(\varrho=0\),
that is, the class number of \(K\) is not divisible by \(3\)
(and thus capitulation is discouraged, \(C=0\)),
implies that there are only three possible types of
\textit{split extensions} \(N/K\), namely the singlets
\((\alpha_3)\), \((\beta_2)\), \((\gamma)\).

\begin{theorem}
\label{thm:SplitExtensions}
(Ramified singlets of \textbf{split extensions} \(N/K\)) \\
Suppose that \(K\) is a real quadratic field with \(3\)-class rank \(\varrho=0\),
and let \(f\) be a regular \(3\)-admissible conductor for \(K\)
with exactly two restrictive prime divisors, \(t=2\).
Denote by \(0\le s\le 2\) the number of prime divisors of \(f\)
which split in the real quadratic field \(K\).
Then the following conditions enforce a split extension \(N/K\)
with \(U_N=U_K\cdot E_{N/K}\), where \(E_{N/K}=U_N\cap\ker(N_{N/K})\)
denotes the subgroup of relative units of \(N/K\).
\begin{enumerate}
\item
If \(s=0\) and \(N\) has \(3\)-class number \(1\),
then \(N\) is a singlet of type \((\gamma)\).
\item
If \(s=1\) and \(N\) has \(3\)-class number \(3\),
then \(N\) is a singlet of type \((\beta_2)\).
\item
If \(s=2\) and \(N\) has \(3\)-class number \(9\),
then \(N\) is a singlet of type \((\alpha_3)\) \(\lbrack\)or \((\gamma)\)\(\rbrack\).
\end{enumerate}
\end{theorem}

\begin{proof}
The assumption \(\varrho=0\) implies that the \(3\)-Selmer space \(V\) of \(K\)
is generated by the fundamental unit \(\eta\) of \(K\).
Since we suppose \(t=2\) with \textit{regular restrictive} prime divisors \(q_1,q_2\) of the conductor \(f\),
in the sense of Remark
\ref{rmk:TwoPrimes},
\(\eta\) is not contained in the \(3\)-ring spaces \(V(q_1)\) and \(V(q_2)\),
and the multiplicity of \(f\) is \(m=1\), i.e., we have a \textit{singlet}
\cite[Thm. 3.3, p. 2217]{Ma2014}.

According to
\cite[Thm. A, p. 70]{Mo1979},
the subgroup \(\mathrm{Cl}_N^\sigma\) of weakly ambiguous ideal classes of \(\mathrm{Cl}_N\),
with respect to \(\mathrm{Gal}(N/K)=\langle\sigma\rangle\),
is of order \(\#\mathrm{Cl}_3(K)\cdot 3^{T-1}/3^Q\),
where the norm index is denoted by \(3^Q=(U_K:(U_K\cap N_{N/K}(N^\times)))\)
and \(T\) is the number of prime ideals of \(K\) which ramify in \(N\).
In our situation, we have \(\#\mathrm{Cl}_3(K)=1\), and
\(\#\mathrm{Cl}_3(N)=3^s\) is divisible by \(3^{T-1-Q}\), where \(T=2+s\),
that is, \(s\ge 1+s-Q\), resp. \(Q\ge 1\) and thus \(Q=1\).
A fortiori, the unit norm index is
\(b=(U_K: N_{N/K}(U_N))=(U_K:(U_K\cap N_{N/K}(N^\times)))\cdot ((U_K\cap N_{N/K}(N^\times)):N_{N/K}(U_N))=3\),
because \(\eta\) is not norm of a number in \(N^\times\), let alone of a unit in \(U_N\).
Thus, \(N/K\) is a split extension, in the sense of Remark
\ref{rmk:RealDPFTypes},
and the types \((\delta_1)\) and \((\varepsilon)\) in item (1) of Theorem
\ref{thm:Singlets}
are impossible.

Equation
\eqref{eqn:Scholz}
in additive form,
\(V_N=2\cdot V_L+V_K+E-2\),
where \(V_F:=v_3(\#\mathrm{Cl}_F)\) for a number field \(F\),
degenerates to \(V_N=2\cdot V_L+E-2\) under our assumption \(V_K=0\).
This gives rise to a parity condition for \(E\):
\begin{enumerate}
\item
If \(s=0\) and \(V_N=0\), then \(2\cdot V_L=2-E\) and \(E\) must be even,
\(E=2\), \(N\) of type \((\gamma)\).
\item
If \(s=1\) and \(V_N=1\), then \(2\cdot V_L=3-E\) and \(E\) must be odd,
\(E=1\), \(N\) of type \((\beta_2)\).
\item
If \(s=2\) and \(V_N=2\), then \(2\cdot V_L=4-E\) and \(E\) must be even, \\
\(\lbrack\)either\(\rbrack\) \(E=0\), \(N\) of type \((\alpha_3)\)
\(\lbrack\)or \(E=2\), \(N\) of type \((\gamma)\), conjectured impossible\(\rbrack\). \qedhere
\end{enumerate}
\end{proof}


\section{Statistical evaluation and theoretical interpretation of the tables}
\label{s:Theoretical}

\noindent
Now we illuminate and analyze our extensive numerical (computational, experimental) results
with the aid of statistical evaluations and theoretical statements.

\subsection{Stagnation and evolution of arithmetical structures}
\label{ss:Arithmetic}

\noindent
Some features in the Tables
\ref{tbl:Angell},
\ref{tbl:GutensteinStreiteben},
\ref{tbl:EnnolaTurunen0},
\ref{tbl:EnnolaTurunen1},
\ref{tbl:LlorenteQuer0}, and
\ref{tbl:LlorenteQuer1}
reveal \textit{stagnation},
that is, multiplicities and DPF types remain constant,
and only the statistical counters show monotonic growth,
usually slightly faster than linear.
Other phenomena stick out with conspicuous \textit{evolution},
leading to new multiplicities and new DPF types.
The huge total number \(592\,922\)
of all objects occurring in our investigation
of the extensive range \(0<d_L<10^7\)
admits sound statistical interpretation
and heuristic predictions in unproven conjectures.

Since the conductor \(f=1\) is \(3\)-admissible
for any quadratic fundamental discriminant \(d\),
the quadratic fields \(K=\mathbb{Q}(\sqrt{d})\) with \(3\)-class rank \(\varrho=0\)
must be considered as giving rise to nilets \(\mathbf{M}_{d}=\emptyset\).
Observe that the \(3\)-ring space \(V(1)\) modulo \(1\) coincides with \(3\)-Selmer space \(V\)
and the multiplicity formula for the unramified situation yields
\(m(d)=m(1^2\cdot d)=\frac{1}{2}(3^{\varrho}-1)=\frac{1}{2}(3^0-1)=0\).

In the following conjectures,
of which certain parts are proven theorems,
we always give successive percentages
with respect to the upper bounds \(10^5\), \(2\cdot 10^5\), \(5\cdot 10^5\) and \(10^7\),
in this order.

\begin{conjecture}
\label{cnj:UnramifiedStagnation}
The relative frequency of unramified nilets with \(\varrho=0\)
slightly decreases from \(89.1\%\) over \(88.6\%\) and \(87.9\%\) to \(86.3\%\).
The relative frequency of unramified singlets with \(\varrho=1\)
slightly increases from \(10.9\%\) over \(11.4\%\) and \(12.1\%\) to \(13.6\%\).
All singlets are of permanent type \(\delta_1\), showing \textit{stagnation}.
See Theorems
\ref{thm:Unramified}
and
\ref{thm:Sufficient}.
The relative frequency of \textbf{unramified quartets} with \(\varrho=2\)
is marginal below \(0.1\%\), but they reveal an interesting \textit{evolution} of types:
\begin{enumerate}
\item
Up to \(10^5\),
\(\frac{4}{5}=80\%\) are of mixed type \((\alpha_1,\alpha_1,\alpha_1,\delta_1)\),
\(\frac{1}{5}=20\%\) of pure type \((\alpha_1,\alpha_1,\alpha_1,\alpha_1)\).
\item
Up to \(2\cdot 10^5\),
\(\frac{14}{16}=87.5\%\) are of type \((\alpha_1,\alpha_1,\alpha_1,\delta_1)\),
\(\frac{2}{16}=12.5\%\) of type \((\alpha_1,\alpha_1,\alpha_1,\alpha_1)\).
\item
Up to \(5\cdot 10^5\),
\(\frac{53}{16}=86.9\%\) are of mixed type \((\alpha_1,\alpha_1,\alpha_1,\delta_1)\),
\(\frac{4}{61}=6.6\%\) of pure type \((\alpha_1,\alpha_1,\alpha_1,\alpha_1)\), and also
\(\frac{4}{61}=6.6\%\) of the \textit{new pure} type \((\delta_1,\delta_1,\delta_1,\delta_1)\).
\item
Up to \(10^7\),
\(\frac{2391}{2870}=83.3\%\) are of mixed type \((\alpha_1,\alpha_1,\alpha_1,\delta_1)\),
\(\frac{234}{2870}=8.6\%\) of pure type \((\delta_1,\delta_1,\delta_1,\delta_1)\),
\(\frac{175}{2870}=6.1\%\) of pure type \((\alpha_1,\alpha_1,\alpha_1,\alpha_1)\), 
\(\frac{62}{2870}=2.2\%\) of the \textit{new mixed} type \((\alpha_1,\delta_1,\delta_1,\delta_1)\), and
\(\frac{8}{2870}=0.3\%\) of another \textit{new mixed} type \((\alpha_1,\alpha_1,\delta_1,\delta_1)\).
\end{enumerate}
\end{conjecture}

\begin{conjecture}
\label{cnj:RamifiedStagnation}
For \(3\)-admissible non-split prime(power) conductors
\(f=q\), \(f=3\), and \(f=9\) with \(d\equiv 2\,(\mathrm{mod}\,3)\)
over quadratic fields \(K=\mathbb{Q}(\sqrt{d})\) with \(3\)-class rank \(\varrho=0\),
the relative frequency of nilets
slightly decreases from \(73\%\) over \(72\%\) and \(71\%\) to \(69\%\),
and the relative frequency of singlets
slightly increases from \(27\%\) over \(28\%\) and \(29\%\) to \(31\%\).
All singlets are of permanent type \(\varepsilon\), showing \textit{stagnation}.
See Theorems
\ref{thm:OnePrimeNonSplit}
and
\ref{thm:Sufficient}.
We conjecture the last percentages
for the range \(0<d_L<10^7\)
to be close to their asymptotic limit.
\end{conjecture}


\subsection{New features for \(3\)-class rank \(\varrho=1\)}
\label{ss:NewFeatures}

\noindent
Since ramified extensions \(N/K\) for \(\varrho=2\)
do not occur in the range \(0<d_L<10^7\),
it is sufficient to state the following theorem for \(\varrho\le 1\).

\begin{theorem}
\label{thm:TwoPrimes}
\noindent
Let \(K=\mathbb{Q}(\sqrt{d})\) be a real quadratic base field with
fundamental discriminant \(d\) and \(3\)-class rank \(\varrho\le 1\).
Suppose \(f=q_1\cdot q_2\) is a \textbf{regular} \(3\)-admissible conductor for \(K\)
with \textbf{two prime divisors} \(q_1\) and \(q_2\).
Then the heterogeneous multiplet \(\mathbf{M}(K_f)\) associated with the
\(3\)-ring class field \(K_f\) mod \(f\) of \(K\)
consists of four homogeneous multiplets \(\mathbf{M}_{c^2d}\), \(c\in\lbrace 1,q_1,q_2,f\rbrace\) with multiplicities
\(m(1)\), \(m(q_1)\), \(m(q_2)\) and \(m(f)\).
In this order,
and in dependence on the \(3\)-ring spaces \(V(q_1)\), \(V(q_2)\) and \(V(f)\),
these four multiplicities, forming the \textbf{signature} \(\mathrm{sgn}(\mathbf{M}(K_f))\) of \(\mathbf{M}(K_f)\), are given by
\begin{enumerate}
\item
\((0,\ 1,\ 1,\ 2)\), \quad if \(V(f)=V(q_1)=V(q_2)=V\) (\textbf{doublet}),
\item
\((0,\ 1,\ 0,\ 0)\), \quad if \(0=V(f)=V(q_2)<V(q_1)=V\),
\item
\((0,\ 0,\ 1,\ 0)\), \quad if \(0=V(f)=V(q_1)<V(q_2)=V\),
\item
\((0,\ 0,\ 0,\ 1)\), \quad if \(0=V(f)=V(q_1)=V(q_2)<V\) (\textbf{singlet}),
\end{enumerate}
if \(\varrho=0\), and thus \(3\)-Selmer space \(V\) is \textbf{one-dimensional},
generated by \(\eta\in U_K=\langle -1,\eta\rangle\), and by
\begin{enumerate}
\item
\((1,\ 3,\ 3,\ 6)\), \quad if \(V(f)=V(q_1)=V(q_2)=V\) (\textbf{sextet}),
\item
\((1,\ 3,\ 0,\ 0)\), \quad if \(0<V(f)=V(q_2)<V(q_1)=V\),
\item
\((1,\ 0,\ 3,\ 0)\), \quad if \(0<V(f)=V(q_1)<V(q_2)=V\),
\item
\((1,\ 0,\ 0,\ 3)\), \quad if \(0<V(f)=V(q_1)=V(q_2)<V\) (\textbf{triplet}),
\item
\((1,\ 0,\ 0,\ 0)\), \quad if \(0=V(f)<V(q_1)\ne V(q_2)<V\) (\textbf{nilet} with defect \(\delta=2\)),
\end{enumerate}
if \(\varrho=1\), and thus \(3\)-Selmer space \(V\) is \textbf{two-dimensional},
generated by \(\eta\in U_K\) and \(\theta\in I\setminus U_K\).
\end{theorem}

\begin{proof}
These statements are special cases with \(p=3\) of
\cite[Thm. 5.1]{Ma2021}.
\end{proof}


\begin{remark}
\label{rmk:TwoPrimes}
We emphasize that in the situation with \(\varrho=0\)
a \textit{complete heterogeneous nilet} with signature \((0,0,0,0)\) is impossible,
because there always exists a totally real cubic field \(L\) with discriminant \(d_L\)
equal to either \((q_1q_2)^2d\) or \(q_1^2d\) or \(q_2^2d\).

This is in contrast to the case \(\varrho=1\)
where a \textit{total heterogeneous nilet} with signature \((1,0,0,0)\),
at least with respect to the \textit{ramified} components, can occur.
In this extreme case of a homogeneous nilet \(\mathbf{M}_{f^2d}\) with defect \(\delta(f)=2\),
neither the fundamental unit \(\eta\) nor the other generating \(3\)-virtual unit \(\theta\)
belong to the ring \(R_f\) modulo \(f\) of \(K\),
i.e. both of them are \textit{deficient}.

We also point out that Theorem
\ref{thm:TwoPrimes}
is not only valid for \(f=q_1q_2\) with primes \(q_i\equiv 2\,(\mathrm{mod}\,3)\)
but also for \(f=3q\) with
\(q_1:=3\), \(d\equiv\pm 3\,(\mathrm{mod}\,9)\), \(q_2:=q\equiv\pm 1\,(\mathrm{mod}\,3)\),
for \(f=9q\) with \(q_1:=9\) (the prime power behaves like a prime, formally),
\(d\equiv\pm 1\,(\mathrm{mod}\,3)\), \(q_2:=q\equiv\pm 1\,(\mathrm{mod}\,3)\),
and for \(f=q_1q_2\) with any primes \(q_i\equiv\pm 1\,(\mathrm{mod}\,3)\).
The statement is independent of the decomposition law of the primes \(q_i\) in \(K\),
but it is essential that the conductor is \textit{regular}, that is,
\(9\nmid f\) if \(d\equiv 6\,(\mathrm{mod}\,9)\).
\end{remark}


\begin{example}
\label{exm:TwoPrimes0}
We explicitly consider the statistical results
for \(\varrho=0\), \(f=q_1q_2\) with \(q_1,q_2\equiv 2\,(\mathrm{mod}\,3)\)
in the most extensive range \(0<d_L<10^7\) (Table
\ref{tbl:LlorenteQuer0}).
Since we want to apply \textit{probability theory} to \textit{independent binary properties},
we must start with data concerning prime conductors \(f=q\).
\begin{itemize}
\item
Let \(f=q\equiv 2\,(\mathrm{mod}\,3)\) prime.
Among \(287877\) admissible discriminants \(q^2d\), \\
\phantom{.}\qquad\(198952\) (\(69\%\)) belong to nilets, realizing the \textit{event} \(V(q)=0\), and \\
\phantom{.}\qquad\(88925\) (\(31\%\)) belong to \textbf{singlets}, realizing the \textit{counter event} \(V(q)=V\). 
\item
For \(f=q_1q_2\), the \textit{four field probability table} for \textit{independent events} yields \\
\phantom{.}\qquad\(P=0.69^2\approx 0.476\) for the event \(\lbrack V(q_1)=0\) and \(V(q_2)=0\rbrack\), \\
\phantom{.}\qquad\(P=0.69\cdot 0.31+0.31\cdot 0.69\approx 0.214+0.214=0.428\) for the (symmetric) event \\
\phantom{.}\qquad\qquad\(\lbrack V(q_1)=0\) and \(V(q_2)=V\rbrack\) or \(\lbrack V(q_1)=V\) and \(V(q_2)=0\rbrack\), \\
\phantom{.}\qquad\(P=0.31^2\approx 0.096\) for the event \(\lbrack V(q_1)=V\) and \(V(q_2)=V\rbrack\), \\
and these \textit{theoretical probabilities} are indeed \textit{compatible with} our \textit{experimental result} that
among \(6227\) admissible discriminants \(f^2d\), \\
\(2706\) (\(43\%\approx 42.8\%\)) belong to nilets, \(\lbrack V(q_1)=0 \land V(q_2)=V\rbrack \lor \lbrack V(q_1)=V \land V(q_2)=0\rbrack\), \\
\(3092\) (\(50\%\approx 47.6\%\)) belong to \textbf{singlets}, realizing the event \(\lbrack V(q_1)=0\) and \(V(q_2)=0\rbrack\), \\
\(429\) (\(7\%\approx 9.6\%\)) belong to \textbf{doublets}, realizing the event \(\lbrack V(q_1)=V\) and \(V(q_2)=V\rbrack\).
\end{itemize}
\end{example}


\noindent
Since almost identical probabilities
as for the conductors \(f=q_1q_2\) with \(q_1,q_2\equiv 2\,(\mathrm{mod}\,3)\)
arise for all the other regular conductors with two prime divisors in Theorem
\ref{thm:TwoPrimes},
mentioned explicitly at the end of Remark
\ref{rmk:TwoPrimes},
we are convinced of the following experimental hypothesis.

\begin{conjecture}
\label{cnj:TwoPrimes0}
(\textbf{Probability} for \(m\in\lbrace 0,1,2\rbrace\) when \(\varrho=0\)) \\
The probabilities \(P\) for the occurrence of various multiplets \((L_1,\ldots,L_m)\)
of totally real cubic fields \(L_i\)
among sets of \(3\)-admissible pairs \((f,d)\) of
regular conductors \(f\) and quadratic fundamental discriminants \(d>0\)
with \(\varrho=0\)
are approximately given as follows:
\begin{enumerate}
\item
\(P\approx 31\%\) for a singlet, and \(P\approx 69\%\) for a nilet,
when \(f=q\),
\item
\(P\approx 7\%\) for a doublet, \(P\approx 50\%\) for a singlet, and \(P\approx 43\%\) for a nilet,
when \(f=q_1q_2\).
\end{enumerate} 
\end{conjecture}

\newpage

\begin{example}
\label{exm:TwoPrimes1}
Now we present new features
for \(\varrho=1\), \(f=q_1q_2\) with \(q_1,q_2\equiv 2\,(\mathrm{mod}\,3)\)
in the most extensive range \(0<d_L<10^7\) (Table
\ref{tbl:LlorenteQuer1}).
Again, we must begin with prime conductors \(f=q\).
\begin{itemize}
\item
Let \(f=q\equiv 2\,(\mathrm{mod}\,3)\) prime.
Among \(41541\) admissible discriminants \(q^2d\), \\
\phantom{.}\qquad\(38302\) (\(92.2\%\)) belong to nilets, realizing the event \(V(q)<V\), and \\
\phantom{.}\qquad\(3239\) (\(7.8\%\)) belong to \textbf{triplets}, realizing the counter event \(V(q)=V\). 
\item
For \(f=q_1q_2\), the four field probability table for independent events yields \\
\phantom{.}\qquad\(P=0.922^2\approx 0.850\) for the event \(\lbrack V(q_1)<V\) and \(V(q_2)<V\rbrack\), \\
\phantom{.}\qquad\(P=0.922\cdot 0.078+0.078\cdot 0.922\approx 0.072+0.072=0.144\) for the (symmetric) event \\
\phantom{.}\qquad\qquad\(\lbrack V(q_1)<V\) and \(V(q_2)=V\rbrack\) or \(\lbrack V(q_1)=V\) and \(V(q_2)<V\rbrack\), \\
\phantom{.}\qquad\(P=0.078^2\approx 0.006\) for the event \(\lbrack V(q_1)=V\) and \(V(q_2)=V\rbrack\), \\
but these theoretical probabilities are \textit{not immediately compatible with} our experimental result that
among \(649\) admissible discriminants \(f^2d\), \\
\phantom{.}\qquad\(534\) (\(82.3\%\)) belong to nilets, \\
\phantom{.}\qquad\(115\) (\(17.7\%\)) belong to \textbf{triplets}, \\
\phantom{.}\qquad\(0\) (\(0\%\approx 0.6\%\)) belong to \textbf{sextets}, realizing the event \(\lbrack V(q_1)=V\) and \(V(q_2)=V\rbrack\). \\
Only the case of sextets is compatible, in the sense that it has simply not occurred yet in this range.
At this point, a new phenomenon appears:
the possibility of \textit{elevated defect} \(\delta(f)=2\),
when \(0=V(f)<V(q_1)\ne V(q_2)<V\).
We have to split the event \(\lbrack V(q_1)<V \land V(q_2)<V\rbrack\),
with theoretical probability \(85.0\%\),
into two cases,
a \textbf{triplet} for \(0<V(f)=V(q_1)=V(q_2)<V\) with \textit{experimental} probability \(17.7\%\),
and a \textbf{nilet} for \(0=V(f)<V(q_1)\ne V(q_2)<V\) with unknown probability,
which can now be calculated as \(85.0\%-17.7\%=67.3\%\), an astonishingly high value.
Eventually, the sum of the probabilities for nilets with \(\delta=1\) and nilets with \(\delta=2\),
that is, \(14.4\%+67.3\%=81.7\%\approx 82.3\%\) agrees with the experimental probability for all nilets, indeed.
\end{itemize}
\end{example}

\begin{conjecture}
\label{cnj:TwoPrimes1}
(\textbf{Probability} for \(m\in\lbrace 0,3,6\rbrace\) when \(\varrho=1\)) \\
The probabilities \(P\) for the occurrence of various multiplets \((L_1,\ldots,L_m)\)
of totally real cubic fields \(L_i\)
among sets of \(3\)-admissible pairs \((f,d)\) of
regular conductors \(f\) and quadratic fundamental discriminants \(d>0\)
with \(\varrho=1\)
are approximately given as follows:
\begin{enumerate}
\item
\(P\approx 8\%\) for a triplet, and \(P\approx 92\%\) for a nilet,
when \(f=q\),
\item
\(P\approx 1\%\) for a sextet, \(P\approx 17\%\) for a triplet, and \(P\approx 82\%\) for a nilet,
when \(f=q_1q_2\).
Among the \(82\%\) for a nilet, there are
\(18\%\) nilets with \(\delta=1\) and \(82\%\) nilets with \(\delta=2\).
\end{enumerate} 
\end{conjecture}

\begin{example}
\label{exm:2And5}
It is illuminating to give particular realizations
of the various multiplets in Theorem
\ref{thm:TwoPrimes}.
Let \(q_1=2\) and \(q_2=5\) and consider the composite conductor \(f=q_1q_2=10\).
\end{example}

\noindent
\(\bullet\)
Among quadratic fundamental discriminants \(d\) with \(\varrho=0\),
there are four \(d\in\lbrace 5,\mathbf{13},21,29\rbrace\)
which give rise to nilets \(\mathbf{M}_{4d}=\emptyset\)
before we find a singlet with conductor \(2\) for \(d=37\), \(d_L=148\), and
there are eight \(d\in\lbrace 8,12,\mathbf{13},17,28,33,37,53\rbrace\)
giving rise to nilets \(\mathbf{M}_{25d}=\emptyset\)
until a singlet with conductor \(5\) occurs for \(d=57\), \(d_L=1425\).
The consequence of the simultaneous nilets
\(\mathbf{M}_{4d}=\mathbf{M}_{25d}=\emptyset\) for \(d=\mathbf{13}\)
is the existence of a \textit{singlet} with conductor \(f=10\) and \(d_L=1300\)
in spite of positive defect \(\delta(10)=1\).
A \textit{nilet} \(\mathbf{M}_{100d}=\emptyset\) with conductor \(f=10\) arises for \(d=37\),
because \(\mathbf{M}_{4d}\) is a singlet and \(\mathbf{M}_{25d}=\emptyset\) is a nilet.
We have to wait for the sixteenth discriminant \(d\) for which \(f=10\) is admissible
in order to encounter the first \textit{doublet} \(\mathbf{M}_{100d}\) for \(d=373\), \(d_L=37300\)
with vanishing defect \(\delta(10)=0\).

\noindent
\(\bullet\)
Among quadratic fundamental discriminants \(d\) with \(\varrho=1\),
the probability \(P\approx 92\%>69\%\) for a nilet with prime conductor is higher, and thus
we have to skip \(56\) discriminants, commencing with \(d\in\lbrace 229,469,\mathbf{733},\ldots\rbrace\)
until the first triplet \(\mathbf{M}_{4d}\) with conductor \(2\) occurs for \(d=7053\), \(d_L=28212\).
Similarly, we must overleap \(7\) discriminants, beginning with \(d\in\lbrace 257,473,568,697,\mathbf{733},\ldots\rbrace\)
before we find a triplet \(\mathbf{M}_{25d}\) with conductor \(5\) for \(d=1257\), \(d_L=31425\).
Now the new feature of elevated defect \(\delta=2\) for positive \(3\)-class rank sets in:
The consequence of the simultaneous nilets
\(\mathbf{M}_{4d}=\mathbf{M}_{25d}=\emptyset\) for \(d=\mathbf{733}\)
is not at all a triplet, but rather a \textit{nilet} \(\mathbf{M}_{100d}=\emptyset\) with \(f=10\),
because the ring spaces \(V(2)\) and \(V(5)\) have trivial meet, whence \(\delta(10)=2\).
This phenomenon continues for further six discriminants starting with \(d\in\lbrace 1373,1957,2213\rbrace\)
until \(V(10)=V(2)=V(5)\) coincide for \(d=3173\), \(d_L=317300\),
giving rise to the first \textit{triplet} \(\mathbf{M}_{100d}\).
Even later, the first nilet with moderate defect \(\delta(10)=1\)
(it is the \(24\)th in the series of nilets)
occurs for \(d=7053\), since
\(\mathbf{M}_{4d}\) is a triplet and \(\mathbf{M}_{25d}=\emptyset\) is a nilet.
This ostensively shows the dominant role of
the \(82\%\) nilets \(\mathbf{M}_{100d}=\emptyset\) with \(\delta=2\)
as opposed to the  \(18\%\) with \(\delta=1\),
according to Conjecture
\ref{cnj:TwoPrimes1}.

\newpage

\subsection{Unramified extensions}
\label{ss:Unramified}

\noindent
The unique conductor without prime divisors is \(f = 1\).
It is \(3\)-admissible for \textit{any} quadratic fundamental discriminant \(d\).

Among the \(3\,039\,653\) quadratic fundamental discriminants in the range \(0 < d < 10^7\),
there are \(2\,623\,325\), resp. \(413\,458\), resp. \(2\,870\),
which give rise to real quadratic number fields \(K = \mathbb{Q}(\sqrt{d})\)
with \(3\)-class rank \(\varrho = \varrho_3(K) = 0\), resp. \(1\), resp. \(2\).
According to the multiplicity formula
\(m = m_3(K,1) = \frac{3^\varrho-1}{3-1}\),
there exist \(0\), resp. \(413\,458\), resp. \(11\,480\),
totally real cubic fields \(L\) with discriminant \(d_L = f^2\cdot d = 1^2\cdot d = d\),
occurring in \textit{nilets}, resp. \textit{singlets}, resp. \textit{quartets}.
The associated normal closure \(N\) of each of these non-Galois cubic fields \(L\)
is \textit{unramified} over its unique quadratic subfield \(K\).

\begin{example}
\label{exm:Unramified}
The smallest discriminant 
with \(\varrho = 0\) is
\(d = 5\).
Although it is an actual quadratic fundamental discriminant,
it is only a \textit{formal} cubic discriminant belonging to a nilet.
The minimal discriminants
\(d = 229\),
resp.
\(d = 32\,009\),
with \(\varrho = 1\), resp. \(\varrho = 2\),
are both, fundamental discriminants of real quadratic fields
and \textit{actual} discriminants of totally real cubic fields
belonging to a singlet, resp. quartet.
The latter two discriminants are contained in the table of Angell
with \(0 < d_L < 10^5\) already.
\end{example}


In the sequel,
we briefly speak about the \textit{type}
\(\tau(N) = \tau(L) \in \lbrace\alpha_1,\alpha_2,\alpha_3,\beta_1,\beta_2,\gamma,\delta_1,\delta_2,\varepsilon\rbrace\)
of a totally real \(S_3\)-field \(N\), resp. its three conjugate cubic subfields \(L\),
when we specify the \textit{differential principal factorization type} of \(N\), resp. \(L\).

\begin{theorem}
\label{thm:Unramified}
Let \(L\) be a non-Galois totally real cubic field
whose normal closure \(N\) is unramified over its quadratic subfield \(K\),
with conductor \(f = 1\).
\begin{enumerate}
\item
If the \(3\)-class group \(\mathrm{Cl}_3(K)\) is non-trivial cyclic,
then \(L\) must be of type \(\tau(L) = \delta_1\).
\item
If \(K\) has \(3\)-class rank \(\varrho \ge 2\),
then two types \(\tau(L) \in \lbrace\alpha_1,\delta_1\rbrace\) are possible for \(L\).
\end{enumerate}
\end{theorem}

\begin{proof}
See Theorem
\ref{thm:Sufficient} (1)
for item (1),
and Theorem
\ref{thm:MainReal}
with \(t=s=0\) and thus \(A=R=0\)
for item (2).
\end{proof}


\subsection{Conductors with a single prime divisor}
\label{ss:OnePrime}

\begin{example}
\label{exm:OnePrimeRankTwo}
It is conspicuous,
that the range \(0 < d < 10^7\) contains an abundance of \(197\) \textit{nilets}
with \textit{formal} cubic discriminants \(f^2\cdot d\)
such that the conductor \(f=q\) is a prime \(q\equiv 2\,(\mathrm{mod}\,3)\)
and the fundamental discriminant \(d\) belongs to a real quadratic field \(K\)
with \(3\)-class rank \(\varrho = 2\).
The smallest values of \(d\) occurring among these \(197\) cases are
\(32\,009\), \(42\,817\), \(62\,501\).
However, the associated formal cubic discriminants appear in reverse order
\(250\,004 = 2^2\cdot 62\,501\), \(1\,070\,425 = 5^2\cdot 42\,817\), \(3\,873\,089 = 11^2\cdot 32\,009\),
due to the conductors which increase in the opposite direction.
In particular, the smallest formal cubic discriminant \(250\,004\)
lies in the range \(0 < d < 5\cdot 10^5\) of Ennola and Turunen already.
Actual \textit{nonets} (\(m=9\)) of cubic fields with these discriminants \textit{do not exist}.
According to a private communication by Karim Belabas on 31 January 2002,
the discriminant \(\mathbf{18\,251\,060}=2^2\cdot 4\,562\,765\) in Theorem
\ref{thm:Alpha1Ramified}
is not only minimal with a ramified (\(f=2\)) component of type \(\alpha_1\),
as required for the proof of the \textbf{Scholz Conjecture},
but even the minimal discriminant of totally real cubic \textbf{nonets} at all
(see \texttt{http://www.algebra.at/KarimDan5.htm}).
\end{example}


\begin{theorem}
\label{thm:OnePrimeNonSplit}
Let \(L\) be a totally real cubic field
whose normal closure \(N\) is ramified over its quadratic subfield \(K\)
with \(\varrho=0\) and conductor \(f\) divisible by a single non-split prime,
\begin{enumerate}
\item
either \(f = q\) a prime \(q\equiv 2\,(\mathrm{mod}\,3)\), inert in \(K\),
\item
or \(f=3\) with \(d\equiv 3\,(\mathrm{mod}\,9)\) or \(d\equiv 6\,(\mathrm{mod}\,9)\)
\item
or \(f=9\) with \(d\equiv 6\,(\mathrm{mod}\,9)\)
\item
or \(f=9\) with \(d\equiv 2\,(\mathrm{mod}\,3)\).
\end{enumerate}
In the second and third case, \(3\) ramifies in \(K\),
in the fourth case, \(3\) remains inert in \(K\). \\
Then \(L\) must necessarily be of type \(\tau(L) = \varepsilon\).
\end{theorem}

\begin{proof}
See Theorem
\ref{thm:Sufficient} (2).
\end{proof}



\newpage

\subsection{General conditions for differential principal factorizations}
\label{ss:Conditions}

\noindent
The nine possible types
\(\tau(L)=\tau(N)\in\lbrace\alpha_1,\alpha_2,\alpha_3,\beta_1,\beta_2,\gamma,\delta_1,\delta_1,\varepsilon\rbrace\)
of differential principal factorizations of a non-cyclic totally real cubic field \(L\),
more precisely of the totally real Galois closure \(N\) of \(L\),
are defined with the aid of three invariants \(A\), \(R\) and \(C\)
which are \(\mathbb{F}_3\)-dimensions of canonical subspaces
of the vector space \(\mathcal{P}_{N/K}/\mathcal{P}_K\) of primitive ambiguous principal ideals
of \(N\) over its quadratic subfield \(K\).


The most restrictive necessary conditions are imposed by the three types \(\alpha_1,\alpha_3,\gamma\)
which are characterized by two-dimensional subspaces.

\begin{theorem}
\label{thm:TwoDimensions}
(Necessary conditions for \textbf{two-dimensional} subspaces)

\begin{enumerate}
\item
For type \(\gamma\) with two-dimensional \textbf{absolute} principal factorization \(A=2\),
the conductor \(f\) must have at least two prime divisors, \(t\ge 2\).
\item
For type \(\alpha_3\) with two-dimensional \textbf{relative} principal factorization \(R=2\),
the conductor \(f\) must have at least two prime divisors which split in \(K\), \(s\ge 2\)
(and a fortiori \(t\ge 2\)).
\item
For type \(\alpha_1\) with two-dimensional \textbf{capitulation} \(C=2\),
the \(3\)-class rank \(\varrho\) of \(K\) must be at least two
(independently of the conductor \(f\ge 1\)).
\end{enumerate}
\end{theorem}

\begin{proof}
We make use of the fundamental inequalities in Corollary
\ref{cor:Estimates}:
\[
0\le A\le\min(n+s,2), \quad
0\le R\le\min(s,2), \quad
0\le C\le\min(\varrho,2).
\]

\begin{enumerate}
\item
Type \(\gamma\) \(\Longleftrightarrow\) \(A=2\) \(\Longrightarrow\) \(\min(n+s,2)=2\), i.e. \(t=n+s\ge 2\).
\item
Type \(\alpha_3\) \(\Longleftrightarrow\) \(R=2\) \(\Longrightarrow\) \(\min(s,2)=2\), i.e. \(s\ge 2\), and thus \(t=n+s\ge s\ge 2\).
\item
Type \(\alpha_1\) \(\Longleftrightarrow\) \(C=2\) \(\Longrightarrow\) \(\min(\varrho,2)=2\), i.e. \(\varrho\ge 2\). \qedhere
\end{enumerate}
\end{proof}


\noindent
Looser necessary conditions are required for non-trivial subspaces.

\begin{theorem}
\label{thm:NonTrivial}
(Necessary conditions for \textbf{one-dimensional} subspaces)

\begin{enumerate}
\item
For the types \(\beta_1,\beta_2,\varepsilon\)
with one-dimensional \textbf{absolute} principal factorization \(A=1\),
the conductor \(f\) must have at least one prime divisor, \(t\ge 1\).
\item
For the types \(\alpha_2,\beta_2,\delta_2\)
with one-dimensional \textbf{relative} principal factorization \(R=1\),
the conductor \(f\) must have at least one prime divisor which splits in \(K\), \(s\ge 1\)
(thus \(t\ge 1\)).
\item
For the types \(\alpha_2,\beta_1,\delta_1\)
with one-dimensional \textbf{capitulation} \(C=1\),
the \(3\)-class rank \(\varrho\) of \(K\) must be at least one.
\end{enumerate}
For each of the types \(\alpha_2,\beta_1,\beta_2\), two suitable among these conditions may be combined.
\end{theorem}

\begin{proof}
According to the definitions of DPF types and the fundamental inequalities in Corollary
\ref{cor:Estimates}:

\begin{enumerate}
\item
Type \(\beta_1,\beta_2,\varepsilon\) \(\Longleftrightarrow\) \(A=1\) \(\Longrightarrow\) \(\min(n+s,2)\ge 1\), i.e. \(t=n+s\ge 1\).
\item
Type \(\alpha_2,\beta_2,\delta_2\) \(\Longleftrightarrow\) \(R=1\) \(\Longrightarrow\) \(\min(s,2)\ge 1\), i.e. \(s\ge 1\), and thus \(t=n+s\ge s\ge 1\).
\item
Type \(\alpha_2,\beta_1,\delta_1\) \(\Longleftrightarrow\) \(C=1\) \(\Longrightarrow\) \(\min(\varrho,2)\ge 1\), i.e. \(\varrho\ge 1\). \qedhere
\end{enumerate}
\end{proof}


\noindent
Due to the fact that the occurrence of absolute principal factorizations is usually unpredictable
as soon as the conductor \(f>1\) has at least one prime divisor, \(t\ge 1\),
there a only very few sufficient conditions for DPF types.
Only two types can be enforced unambiguously.

\begin{theorem}
\label{thm:Sufficient}
(\textbf{Sufficient} conditions for types \(\delta_1\) and \(\varepsilon\))
\begin{enumerate}
\item
If \(N/K\) is unramified with conductor \(f=1\) and \(K\) has \(3\)-class rank \(\varrho=1\),
then \(\tau(N)=\delta_1\).
\item
If the conductor \(f\) of \(N/K\) has precisely one prime divisor which does not split in \(K\)
and the class number of \(K\) is not divisible by \(3\),
then \(\tau(N)=\varepsilon\).
\end{enumerate}
In both cases, there exists a unit \(H\in U_N\) such that \(\eta=N_{N/K}(H)\) is a fundamental unit of \(K\).
\end{theorem}

\begin{proof}
According to the fundamental inequalities in Corollary
\ref{cor:Estimates}
and the fundamental equation in Corollary
\ref{cor:Trichotomy},
we have:

\begin{enumerate}
\item
\(t=0\), \(\varrho=1\)
\(\Longrightarrow\) \(A\le\min(n+s,2)=\min(t,2)=0\), \(s\le t=0\), \(R\le\min(s,2)=0\), \(C\le\min(\varrho,2)=1\),
but on the other hand \(C=0+0+C=A+R+C=U+1\ge 1\)
\(\Longrightarrow\) \(A=R=0\), \(C=1\) \(\Longleftrightarrow\) Type \(\delta_1\).
\item
\(t=1\), \(s=0\), \(\varrho=0\)
\(\Longrightarrow\) \(A\le\min(n+s,2)=\min(t,2)=1\), \(R\le\min(s,2)=0\), and \(C\le\min(\varrho,2)=0\),
but on the other hand \(A=A+0+0=A+R+C=U+1\ge 1\)
\(\Longrightarrow\) \(A=1\), \(R=C=0\) \(\Longleftrightarrow\) Type \(\varepsilon\).
\end{enumerate}
In both cases, we obtain \(U=0\) as a byproduct, i.e. \(N_{N/K}(U_N)=U_K\).
\end{proof}


\section{Complete verification of the Scholz conjecture}
\label{s:ScholzConjecture}

\begin{figure}[ht]
\caption{Hasse subfield diagram of the normal closure \(N/\mathbb{Q}\) of \(L\)}
\label{fig:AbsoluteSubfields}

{\small

\setlength{\unitlength}{1.0cm}
\begin{picture}(5,5)(-7,-9.4)



\put(-6,-9){\circle*{0.2}}
\put(-6,-9.2){\makebox(0,0)[ct]{\(\mathbb{Q}=L\cap K\)}}
\put(-7,-9){\makebox(0,0)[rc]{rational number field}}

\put(-6,-9){\line(2,1){2}}
\put(-5,-8.7){\makebox(0,0)[lt]{\(\lbrack K:\mathbb{Q}\rbrack=2\)}}

\put(-4,-8){\circle*{0.2}}
\put(-4,-8.2){\makebox(0,0)[ct]{\(K\)}}
\put(-3,-8){\makebox(0,0)[lc]{quadratic field}}


\put(-6.2,-7.5){\makebox(0,0)[rc]{\(\lbrack L:\mathbb{Q}\rbrack=3\)}}
\put(-6,-9){\line(0,1){3}}
\put(-4,-8){\line(0,1){3}}



\put(-6,-6){\circle{0.2}}
\put(-6,-5.8){\makebox(0,0)[cb]{\(L\)}}
\put(-5.8,-6){\makebox(0,0)[lt]{\(L^\prime,L^{\prime\prime}\)}}
\put(-7,-6){\makebox(0,0)[rc]{three conjugate cubic fields}}

\put(-6,-6){\line(2,1){2}}

\put(-4,-5){\circle*{0.2}}
\put(-4,-4.8){\makebox(0,0)[cb]{\(N=L\cdot K\)}}
\put(-3,-5){\makebox(0,0)[lc]{\(S_3\)-field (dihedral field of degree \(6\))}}


\end{picture}

}

\end{figure}
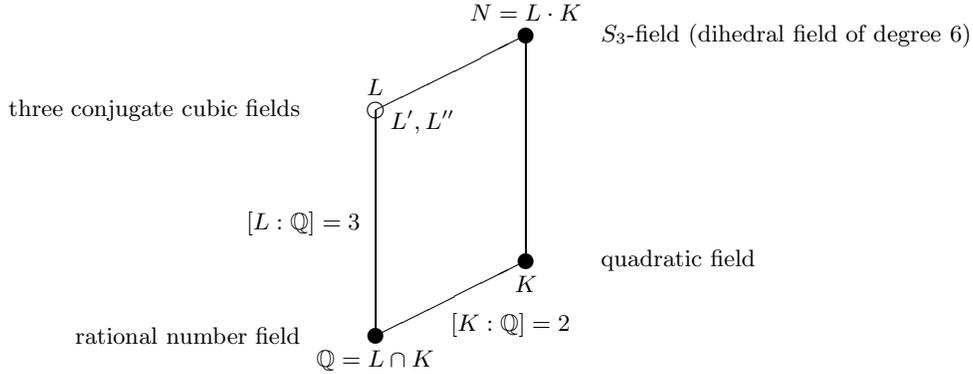

\noindent
Let \(L\) be a \textit{non-cyclic totally real} cubic field.
Then \(L\) is non-Galois over the rational number field \(\mathbb{Q}\)
with two conjugate fields \(L^\prime\) and \(L^{\prime\prime}\).
The Galois closure \(N\) of \(L\) is a totally real dihedral field of degree \(6\),
i.e. an \(S_3\)-field,
which contains a unique real quadratic field \(K\),
as illustrated in Figure
\ref{fig:AbsoluteSubfields}.


In \(1930\), Hasse
\cite{Ha1930}
determined the \textit{discriminants} \(d_L\) of \(L\)
\cite[pp. 567 (1) and 575]{Ha1930}
and \(d_N\) of \(N\)
\cite[p. 566 (2)]{Ha1930},
in dependence on the discriminant \(d=d_K\) of \(K\)
and on the class field theoretic \textit{conductor} \(f=f_{N/K}\)
of the cyclic cubic, and thus abelian, relative extension \(N/K\):

\begin{equation}
\label{eqn:Hasse}
d_L=f^2\cdot d, \quad \text{ and } \quad d_N=f^4\cdot d^3.
\end{equation}

Three years later, Scholz
\cite[p. 216]{So1933}
determined the \textit{relation}

\begin{equation}
\label{eqn:Scholz}
h_N=\frac{a}{9}\cdot h_L^2\cdot h_K
\end{equation}

\noindent
\textit{between the class numbers} of the fields \(N\), \(L\) and \(K\),
in dependence on the \textit{index of subfield units}, \(a=(U_N:U_0)=3^E\),
where \(U_0=\langle U_K,U_L,U_{L^\prime},U_{L^{\prime\prime}}\rangle\)
and \(E\in\lbrace 0,1,2\rbrace\).

Note that \(E=0\), respectively \(a=1\), is the \textit{distinguished situation}
where the unit group \(U_N\) of the normal field \(N\)
is entirely generated by all proper subfield units, that is, \(U_N=U_0\).


Scholz was able to give explicit numerical examples
\cite[p. 216]{So1933}
for \(E=1\) (e.g. \(d_L=229\)),
and \(E=2\) (e.g. \(d_L=148\)),
but not for \(E=0\),
and he formulated the following hypothesis.

\begin{figure}[ht]
\caption{Hilbert and ring class fields over \(K\)}
\label{fig:RingClassFields}

{\small

\setlength{\unitlength}{1.0cm}
\begin{picture}(10,6)(-8,-10)



\put(-8,-9){\circle*{0.2}}
\put(-8,-9.2){\makebox(0,0)[ct]{\(K\)}}

\put(-8,-9){\line(-2,1){2}}
\put(-8,-9){\line(2,1){2}}
\put(-6.8,-8.5){\makebox(0,0)[lt]{\(\lbrack N:K\rbrack=3\)}}

\put(-10,-8){\line(2,1){2}}
\put(-6,-8){\line(-2,1){2}}

\put(-8,-7.9){\makebox(0,0)[cb]{unramified quartet}}
\put(-10,-8){\circle*{0.2}}
\put(-8.7,-8){\circle*{0.2}}
\put(-7.3,-8){\circle*{0.2}}
\put(-6,-8){\circle*{0.2}}
\put(-5.8,-8){\makebox(0,0)[lc]{\(N\)}}

\put(-8,-7){\circle*{0.2}}

\put(-8,-7){\line(0,1){2}}

\put(-9.2,-5){\makebox(0,0)[rc]{(1)}}
\put(-8,-5){\circle*{0.2}}
\put(-7.8,-5){\makebox(0,0)[lc]{\(\mathrm{F}_{3,1}(K)=\mathrm{F}_3^1(K)\)}}
\put(-7.8,-5.5){\makebox(0,0)[lc]{Hilbert \(3\)-class field of \(K\)}}



\put(0,-9){\circle*{0.2}}
\put(0,-9.2){\makebox(0,0)[ct]{\(K\)}}

\put(0,-9){\line(0,1){1}}
\put(0.2,-8.5){\makebox(0,0)[lc]{\(\lbrack N:K\rbrack=3\)}}
\put(0,-8){\line(0,1){3}}

\put(-0.2,-7.9){\makebox(0,0)[rb]{? ramified singlet ?}}
\put(0,-8){\circle{0.2}}
\put(0.2,-8){\makebox(0,0)[lc]{\(N\)}}

\put(-1.2,-5){\makebox(0,0)[rc]{(2)}}
\put(0,-5){\circle*{0.2}}
\put(0.2,-5){\makebox(0,0)[lc]{\(\mathrm{F}_{3,f}(K)=K_f\)}}
\put(0.2,-5.5){\makebox(0,0)[lc]{\(3\)-ring class field mod \(f\) of \(K\)}}


\end{picture}

}

\end{figure}
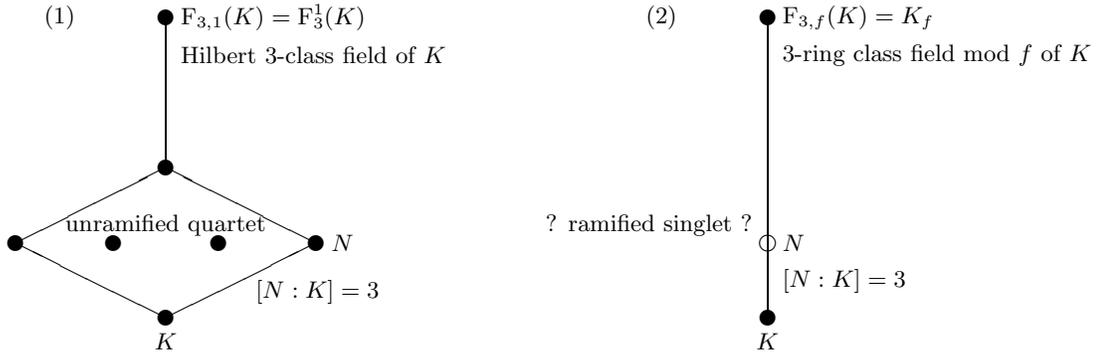

\begin{conjecture}
\label{cnj:Scholz}
(The \textbf{Scholz Conjecture, 1933}, illustrated in Figure
\ref{fig:RingClassFields}) \\
There should exist non-Galois totally real cubic fields \(L\)
whose Galois closure \(N\) is either 
\begin{enumerate}
\item
\textit{unramified}, with conductor \(f=1\),
over some real quadratic field \(K\) with \(3\)-class rank \(\varrho_3(K)=2\)
whose complete \(3\)-elementary class group capitulates in \(N\)
such that \(U_N=U_0\)
(in the terminology of Scholz, \(N\) is an \textit{absolute class field} over \(K\))
\cite[p. 219]{So1933},
or
\item
\textit{ramified}, with conductor \(f>1\),
over some real quadratic field \(K\)
such that \(U_N=U_0\)
(here, Scholz calls \(N\) a \textit{ring class field} over \(K\), by abuse of language)
\cite[p. 221]{So1933}.
\end{enumerate}
\end{conjecture}


\noindent
We point out that, in the \textit{unramified} situation \(f=1\),
\(d_L=d\) is a quadratic fundamental discriminant, and
\(d_N=d^3\) is a perfect cube,
according to Formula
\eqref{eqn:Hasse}.
In this unramified case,
the verification of Conjecture
\ref{cnj:Scholz}
can be obtained from a more general theorem, since
any real quadratic field \(K\) with \(3\)-class rank \(\varrho_3(K)=2\)
possesses a multiplet of four unramified cyclic cubic extensions
\(N_1,\ldots,N_4\),
that is a \textit{quartet} of absolutely dihedral fields of degree \(6\)
\cite{Ma2012}
with non-Galois totally real cubic subfields \(L_1,\ldots,L_4\),
each of them selected among three conjugate fields.


For such a quartet,
Chang and Foote
\cite{ChFt1980}
introduced the concept of the \textit{capitulation number} \(0\le\nu(K)\le 4\),
defined as the number of those members of the quartet
in which the complete \(3\)-elementary class group of \(K\) capitulates.
For this number \(\nu(K)\), the following theorem holds.

\begin{theorem}
\label{thm:Alpha1Unramified}
For each value \(0\le\nu\le 4\),
there exists a real quadratic field \(K\) with \(3\)-class rank \(\varrho_3(K)=2\)
and capitulation number \(\nu(K)=\nu\).
It is even possible to restrict the claim to fields with
elementary \(3\)-class group of type \(\mathrm{Cl}_3(K)\simeq C_3\times C_3\).
\end{theorem}

\newpage

\begin{proof}
From the viewpoint of finite \(p\)-group theory,
this theorem is a proven statement about the possible \textit{transfer kernel types}
of finite metabelian \(3\)-groups \(G\) with
abelianization \(G/G^\prime\simeq (3,3)\)
applied to the second \(3\)-class group \(G:=\mathrm{Gal}(F_3^2(K)/K)\) of \(K\)
\cite{Ma2012}.
However, it is easier to give explicit numerical paradigms for each value of \(\nu(K)\).
We have the following minimal occurrences: \\
\(\nu(K)=4\) for \(d_K=62\,501\), \\
\(\nu(K)=3\) for \(d_K=32\,009\), \\
\(\nu(K)=2\) for \(d_K=710\,652\), \\
\(\nu(K)=1\) for \(d_K=534\,824\), \\
\(\nu(K)=0\) for \(d_K=214\,712\), \\
which have been computed by ourselves in
\cite{Ma2012}.
The existence of these cases completes the proof.
\end{proof}


\begin{remark}
\label{rmk:Alpha1Unramified}
We have the priority of discovering the first examples of
real quadratic fields \(K\) with \(\nu(K)\in\lbrace 0,1,2\rbrace\) in
\cite{Ma2012}.
However, the first examples of
real quadratic fields \(K\) with \(\nu(K)\in\lbrace 3,4\rbrace\)
are due to Heider and Schmithals
\cite{HeSm1982},
who performed a mainframe computation on the CDC Cyber of the University at Cologne,
and thus the following corollary is proven since \(1982\) already.
\end{remark}

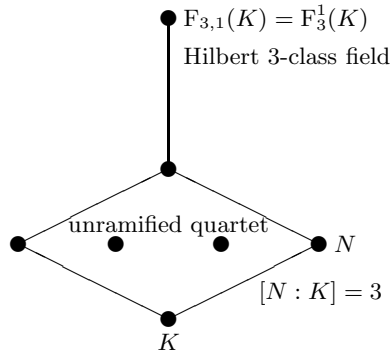
\begin{figure}[h]
\caption{Hilbert class field over \(K\)}
\label{fig:HilbertClassFieldQuartet}

{\small

\setlength{\unitlength}{1.0cm}
\begin{picture}(10,5)(-12,-9.6)



\put(-8,-9){\circle*{0.2}}
\put(-8,-9.2){\makebox(0,0)[ct]{\(K\)}}

\put(-8,-9){\line(-2,1){2}}
\put(-8,-9){\line(2,1){2}}
\put(-6.8,-8.5){\makebox(0,0)[lt]{\(\lbrack N:K\rbrack=3\)}}

\put(-10,-8){\line(2,1){2}}
\put(-6,-8){\line(-2,1){2}}

\put(-8,-7.9){\makebox(0,0)[cb]{unramified quartet}}
\put(-10,-8){\circle*{0.2}}
\put(-8.7,-8){\circle*{0.2}}
\put(-7.3,-8){\circle*{0.2}}
\put(-6,-8){\circle*{0.2}}
\put(-5.8,-8){\makebox(0,0)[lc]{\(N\)}}

\put(-8,-7){\circle*{0.2}}

\put(-8,-7){\line(0,1){2}}

\put(-8,-5){\circle*{0.2}}
\put(-7.8,-5){\makebox(0,0)[lc]{\(\mathrm{F}_{3,1}(K)=\mathrm{F}_3^1(K)\)}}
\put(-7.8,-5.5){\makebox(0,0)[lc]{Hilbert \(3\)-class field of \(K\)}}


\end{picture}

}

\end{figure}

\begin{corollary}
\label{cor:Alpha1Unramified}
(Verification of Conjecture \ref{cnj:Scholz}, (1) for \textbf{unramified} extensions; see Figure
\ref{fig:HilbertClassFieldQuartet}) \\
There exist non-Galois totally real cubic fields \(L\)
whose Galois closure \(N\) is unramified, with conductor \(f=1\),
over a real quadratic field \(K\) with \(3\)-class rank \(\varrho_3(K)=2\)
whose complete \(3\)-elementary class group capitulates in \(N\),
and which therefore has \(U_N=U_0\).
The minimal discriminant of such a field \(L\) is \(d_L=32\,009\)
(discovered in \cite{HeSm1982},
actually, the first three members of this \textbf{quartet} with DPF type
\((\alpha_1,\alpha_1,\alpha_1,\delta_1)\)
in Table
\ref{tbl:AngellQuartets}
satisfy the relation \(U_N=U_0\)).
\end{corollary}


\begin{proof}
It suffices to take a real quadratic field \(K\) with \(1\le \nu(K)\le 4\)
in Theorem
\ref{thm:Alpha1Unramified}.
In view of the minimal discriminant, we select \(\nu(K)=3\)
and obtain \(U_N=U_0\) for \(d_L=d_K=32\,009\).
\end{proof}


\noindent
Concerning the \textit{ramified} situation  \(f>1\) in Conjecture
\ref{cnj:Scholz} (2),
Scholz does not explicitly impose any conditions
on the underlying real quadratic field \(K\).
We suppose that he also tacitly assumed
a real quadratic field \(K\) with \(3\)-class rank \(\varrho_3(K)=2\).
However, more recent extensions of the theory of dihedral fields
by means of \textit{differential principal factorizations} and \textit{Galois cohomology},
two concepts which we have expanded thoroughly in the preparatory sections \S\S\
\ref{ss:NormKernel},
\ref{ss:CapitulationKernel}, and
\ref{ss:GaloisCohomology},
revealed that for \(U_N=U_0\)
no constraints on the \(p\)-class rank \(\varrho_p(K)\) are required.
In \(1975\), Nicole Moser
\cite{Mo1979}
used the \textit{Galois cohomology}
\(\hat{\mathrm{H}}^0(G,U_N)\simeq U_K/\mathrm{N}_{N/K}(U_N)\)
of the unit group \(U_N\) of the normal closure \(N\) as a module over \(G=\mathrm{Gal}(N/K)\)
to establish a \textit{fine structure} with five possible types \(\alpha,\beta,\gamma,\delta,\varepsilon\)
on the \textit{coarse} classification by three possible values of the index of subfield units: \\
\((U_N:U_0)=1\) \(\Longleftrightarrow\) type \(\alpha\) with \((U_K:\mathrm{N}_{N/K}(U_N))=3\), \\
\((U_N:U_0)=3\) \(\Longleftrightarrow\) type \(\beta\) with \((U_K:\mathrm{N}_{N/K}(U_N))=3\) or type \(\delta\) with \((U_K:\mathrm{N}_{N/K}(U_N))=1\), \\
\((U_N:U_0)=9\) \(\Longleftrightarrow\) type \(\gamma\) with \((U_K:\mathrm{N}_{N/K}(U_N))=3\) or type \(\varepsilon\) with \((U_K:\mathrm{N}_{N/K}(U_N))=1\). \\
Thus, Moser's refinement does not illuminate the situation \(U_N=U_0\) (\(\Longleftrightarrow\) type \(\alpha\)) of Scholz's conjecture more closely.
Meanwhile, Barrucand and Cohn
\cite{BaCo1971}
had coined the concept of \textit{(differential) principal factorization} ((D)PF) for pure cubic fields.
In \(1991\), we generalized the theory of DPFs for dihedral fields of both signatures
\cite{Ma1991b},
and we obtained a \textit{hyperfine structure} by splitting Moser's types further
according to the \(\mathbb{F}_p\)-dimensions
\(C\) of the capitulation kernel \(\ker(T_{K,N})\) and
\(R\) of the space of relative DPFs of \(N/K\),
which we recalled in the preparatory section \S\
\ref{ss:RealDihedralAndQuinticTypes}.
In particular, type \(\alpha\) with \(U_N=U_0\) splits into three subtypes: \\
type \(\alpha_1\) \(\Longleftrightarrow\) \(C=2\), \(R=0\), which implies \(\varrho_p(K)\ge 2\), \\
type \(\alpha_2\) \(\Longleftrightarrow\) \(C=1\), \(R=1\), which implies \(\varrho_p(K)\ge 1\) and a split prime divisor of \(f\) \((s\ge 1)\), \\
type \(\alpha_3\) \(\Longleftrightarrow\) \(C=0\), \(R=2\), which is compatible with any \(\varrho_p(K)\ge 0\), but requires \(s\ge 2\).


Consequently, we were led to the following refinement of Conjecture
\ref{cnj:Scholz}, (2).

\begin{conjecture}
\label{cnj:Mayer}
(Conjecture of D. C. Mayer, \(1991\)) \\
Non-Galois totally real cubic fields \(L\)
whose Galois closure \(N\) is \textit{ramified}, with conductor \(f>1\),
over some real quadratic field \(K\),
and is of type \(\alpha\), with \(U_N=U_0\),
should exist for each of the following three situations:
\begin{enumerate}
\item[(2.1)]
type \(\alpha_1\) with
\(\dim_{\mathbb{F}_3}(\ker(T_{K,N}))=2\) and \(\varrho_3(K)=2\), \(s=0\),
\item[(2.2)]
type \(\alpha_2\) with
\(\dim_{\mathbb{F}_3}(\ker(T_{K,N}))=1\) and \(\varrho_3(K)=1\), \(s=1\),
\item[(2.3)]
type \(\alpha_3\) with
\(\dim_{\mathbb{F}_3}(\ker(T_{K,N}))=0\) and \(\varrho_3(K)=0\), \(s=2\),
\end{enumerate}
where \(T_{K,N}:\,\mathrm{Cl}_3(K)\to\mathrm{Cl}_3(N)\), \(\mathfrak{a}\cdot\mathcal{P}_K\mapsto(\mathfrak{a}\mathcal{O}_N)\cdot\mathcal{P}_N\),
denotes the \textit{transfer homomorphism} of \(3\)-classes from \(K\) to \(N\),
and \(s\) counts the prime divisors of the conductor \(f\) which \textit{split} in \(K\).
\end{conjecture}


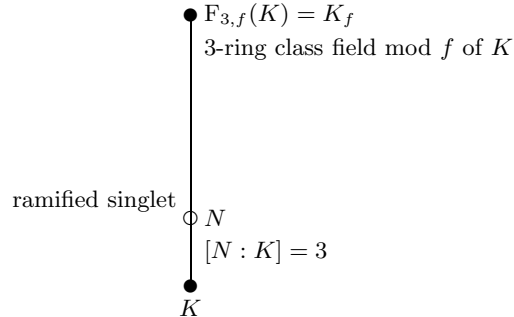
\begin{figure}[ht]
\caption{Ring class field modulo \(f=63=3^2\cdot 7\) over \(K\)}
\label{fig:RingClassFieldSingulet}

{\small

\setlength{\unitlength}{0.9cm}
\begin{picture}(10,4.5)(-4,-9.2)



\put(0,-9){\circle*{0.2}}
\put(0,-9.2){\makebox(0,0)[ct]{\(K\)}}

\put(0,-9){\line(0,1){1}}
\put(0.2,-8.5){\makebox(0,0)[lc]{\(\lbrack N:K\rbrack=3\)}}
\put(0,-8){\line(0,1){3}}

\put(-0.2,-7.9){\makebox(0,0)[rb]{ramified singlet}}
\put(0,-8){\circle{0.2}}
\put(0.2,-8){\makebox(0,0)[lc]{\(N\)}}

\put(0,-5){\circle*{0.2}}
\put(0.2,-5){\makebox(0,0)[lc]{\(\mathrm{F}_{3,f}(K)=K_f\)}}
\put(0.2,-5.5){\makebox(0,0)[lc]{\(3\)-ring class field mod \(f\) of \(K\)}}


\end{picture}

}

\end{figure}


\begin{theorem}
\label{thm:Alpha3Ramified}
(Verification of Conjecture \ref{cnj:Mayer}, (2.3), and Conjecture \ref{cnj:Scholz}, (2); see Figure
\ref{fig:RingClassFieldSingulet}) \\
There exist non-Galois totally real cubic fields \(L\)
whose Galois closure \(N\) is ramified,
with conductor \(f>1\) divisible by two prime divisors which split in \(K\), i.e. \(s=2\),
over a real quadratic field \(K\) with \(3\)-class rank \(\varrho_3(K)=0\),
without capitulation in \(N\),
and such that \(U_N=U_0\).
The minimal discriminant of such a field \(L\) is \(d_L=\mathbf{146\,853}=(7\cdot 9)^2\cdot 37\)
(which forms a \textbf{singlet} \cite{Ma1991c}).
\end{theorem}

\begin{proof}
This was proved in the numerical supplement
\cite{Ma1991c}
of our paper
\cite{Ma1991b}
by computing a gapless list of all \(10\,015\) totally real cubic fields \(L\)
with discriminants \(d_L<200\,000\)
on the AMDAHL mainframe of the University of Manitoba.
There occurred the minimal discriminant \(d_L=146\,853=f^2\cdot d_K\)
with \(d_K=37\) and conductor \(f=63=3^2\cdot 7\)
divisible by two primes which both split in \(K\), i.e. \(s=2\).
This is a necessary requirement for a two-dimensional relative principal factorization with \(R=2\)
and is unique up to \(d_L<200\,000\).
(The next is \(d_L=240\,149\) with \(f=7\cdot 13\).)
There is only a single field \(L\) with this discriminant \(d_L=146\,853\) (forming a singlet).
\end{proof}


Our discovery of the truth of Theorem
\ref{thm:Alpha3Ramified}
with the aid of the list
\cite{Ma1991c}
was a random hit without explicit intention to find a verification of Scholz's conjecture.
Unfortunately,
\cite{Ma1991c}
does not contain examples of the unique missing DPF type \(\alpha_2\).
It required more than \(25\) years until we focused on an attack against this lack of information.
In contrast to the techniques of
\cite{Ma1991c},
we did not use the Voronoi algorithm
\cite{Vo1896}
after cumbersome preparation of generating polynomials for totally real cubic fields,
but rather Fieker's class field theory routines of Magma
\cite{BCP1997,BCFS2020,Fi2001,MAGMA2020}
for a direct generation of the fields as subfields of \(3\)-ray class fields modulo conductors \(f>1\).

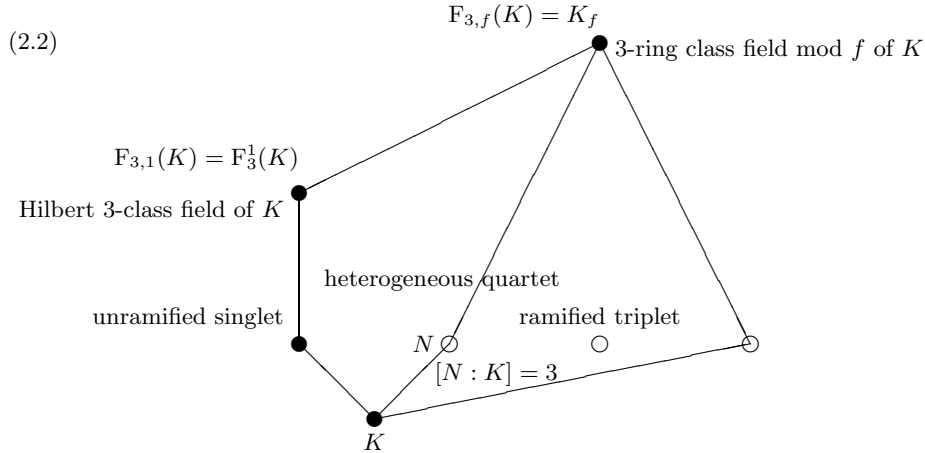
\begin{figure}[ht]
\caption{Heterogeneous \textbf{quartet} modulo \(f=19\) over \(K\)}
\label{fig:RingClassFieldQuartet}

{\small

\setlength{\unitlength}{1.0cm}
\begin{picture}(10,6)(-8,-9)



\put(-5,-9){\circle*{0.2}}
\put(-5,-9.2){\makebox(0,0)[ct]{\(K\)}}

\put(-5,-9){\line(-1,1){1}}
\put(-5,-9){\line(1,1){1}}
\put(-5,-9){\line(5,1){5}}

\put(-4.1,-7.3){\makebox(0,0)[cb]{heterogeneous quartet}}
\put(-6.2,-7.8){\makebox(0,0)[rb]{unramified singlet}}
\put(-6,-8){\circle*{0.2}}

\put(-6,-8){\line(0,1){2}}

\put(-9.2,-4){\makebox(0,0)[rc]{(2.2)}}
\put(-6,-6){\circle*{0.2}}
\put(-6,-5.7){\makebox(0,0)[rb]{\(\mathrm{F}_{3,1}(K)=\mathrm{F}_3^1(K)\)}}
\put(-6.2,-6.2){\makebox(0,0)[rc]{Hilbert \(3\)-class field of \(K\)}}



\put(-4.2,-8.4){\makebox(0,0)[lc]{\(\lbrack N:K\rbrack=3\)}}
\put(-2,-4){\line(-2,-1){4}}
\put(-2,-4){\line(-1,-2){2}}
\put(-2,-4){\line(1,-2){2}}

\put(-2,-7.8){\makebox(0,0)[cb]{ramified triplet}}
\put(-4,-8){\circle{0.2}}
\put(-2,-8){\circle{0.2}}
\put(0,-8){\circle{0.2}}
\put(-4.2,-8){\makebox(0,0)[rc]{\(N\)}}

\put(-2,-4){\circle*{0.2}}
\put(-2,-3.8){\makebox(0,0)[rb]{\(\mathrm{F}_{3,f}(K)=K_f\)}}
\put(-1.8,-4.1){\makebox(0,0)[lc]{\(3\)-ring class field mod \(f\) of \(K\)}}


\end{picture}

}

\end{figure}

\begin{theorem}
\label{thm:Alpha2Ramified}
(Verification of Conjecture \ref{cnj:Mayer}, (2.2), and Conjecture \ref{cnj:Scholz}, (2); see Figure
\ref{fig:RingClassFieldQuartet}) \\
There exist non-Galois totally real cubic fields \(L\)
whose Galois closure \(N\) is ramified,
with conductor \(f>1\) divisible by a single prime divisor that splits in \(K\), i.e. \(s=1\),
over a real quadratic field \(K\) with \(3\)-class rank \(\varrho_3(K)=1\),
with one-dimensional capitulation of the elementary \(3\)-class group in \(N\),
and such that \(U_N=U_0\).
The minimal discriminant of such a field \(L\) is \(d_L=\mathbf{966\,397}=19^2\cdot 2\,677\)
(the first two fields of a \textbf{triplet} \((\alpha_2,\alpha_2,\delta_1)\), discovered \(19\) November \(2017\)).
\end{theorem}

\begin{proof}
The proof is conducted in the following section \S\
\ref{ss:Rank1}.
\end{proof}

\begin{figure}[ht]
\caption{Heterogeneous \textbf{tridecuplet} modulo \(f=2\) over \(K\)}
\label{fig:RingClassFieldTridecuplet}

{\small

\setlength{\unitlength}{1.0cm}
\begin{picture}(10,7)(-8,-9)



\put(-5,-9){\circle*{0.2}}
\put(-5,-9.2){\makebox(0,0)[ct]{\(K\)}}

\put(-5,-9){\line(-5,1){5}}
\put(-5,-9){\line(-1,1){1}}
\put(-5,-9){\line(1,1){1}}
\put(-5,-9){\line(5,1){5}}

\put(-10,-8){\line(2,1){2}}
\put(-6,-8){\line(-2,1){2}}

\put(-5.5,-7.3){\makebox(0,0)[cb]{heterogeneous tridecuplet}}
\put(-8,-7.9){\makebox(0,0)[cb]{unramified quartet}}
\put(-10,-8){\circle*{0.2}}
\put(-8.7,-8){\circle*{0.2}}
\put(-7.3,-8){\circle*{0.2}}
\put(-6,-8){\circle*{0.2}}

\put(-8,-7){\circle*{0.2}}

\put(-8,-7){\line(0,1){2}}

\put(-9.2,-3){\makebox(0,0)[rc]{(2.1)}}
\put(-8,-5){\circle*{0.2}}
\put(-8,-4.8){\makebox(0,0)[rb]{\(\mathrm{F}_{3,1}(K)=\mathrm{F}_3^1(K)\)}}
\put(-7.8,-5.2){\makebox(0,0)[lc]{Hilbert \(3\)-class field of \(K\)}}



\put(-4.2,-8.4){\makebox(0,0)[lc]{\(\lbrack N:K\rbrack=3\)}}
\put(0,-3){\line(-4,-1){8}}
\put(0,-3){\line(-4,-5){4}}
\put(0,-3){\line(0,-1){5}}

\put(-2,-7.8){\makebox(0,0)[cb]{ramified nonet}}
\put(-4,-8){\circle{0.2}}
\put(-3.5,-8){\circle{0.2}}
\put(-3,-8){\circle{0.2}}
\put(-2.5,-8){\circle{0.2}}
\put(-2,-8){\circle{0.2}}
\put(-1.5,-8){\circle{0.2}}
\put(-1,-8){\circle{0.2}}
\put(-0.5,-8){\circle{0.2}}
\put(0,-8){\circle{0.2}}
\put(-4.2,-8){\makebox(0,0)[rc]{\(N\)}}

\put(0,-3){\circle*{0.2}}
\put(0,-2.8){\makebox(0,0)[rb]{\(\mathrm{F}_{3,f}(K)=K_f\)}}
\put(0.2,-3.1){\makebox(0,0)[lc]{\(3\)-ring class field mod \(f\) of \(K\)}}


\end{picture}

}

\end{figure}
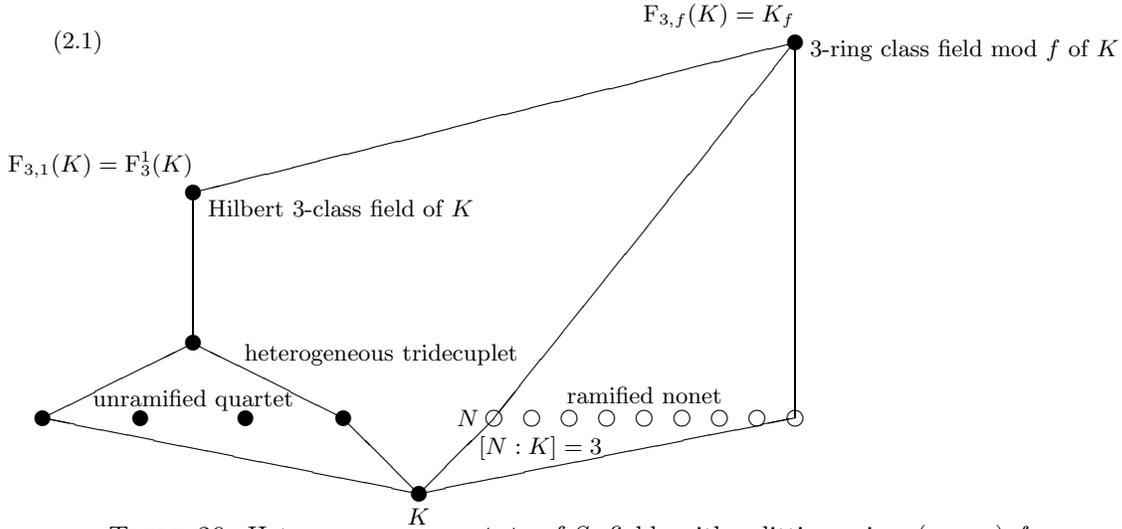

\begin{theorem}
\label{thm:Alpha1Ramified}
(Verification of Conjecture \ref{cnj:Mayer}, (2.1), and Conjecture \ref{cnj:Scholz}, (2); see Figure
\ref{fig:RingClassFieldTridecuplet}) \\
There exist non-Galois totally real cubic fields \(L\)
whose Galois closure \(N\) is ramified,
with conductor \(f>1\) divisible only by prime divisors which do not split in \(K\), i.e. \(s=0\),
over a real quadratic field \(K\) with \(3\)-class rank \(\varrho_3(K)=2\),
with two-dimensional capitulation of the elementary \(3\)-class group in \(N\),
and such that \(U_N=U_0\).
The minimal discriminant of such a field \(L\) is
\[
d_L=\mathbf{18\,251\,060}=2^2\cdot 4\,562\,765
\]
(the first five fields of a \textbf{nonet} \((\alpha_1,\alpha_1,\alpha_1,\alpha_1,\alpha_1,\beta_1,\delta_1,\delta_1,\delta_1)\), discovered \(23\) November \(2017\)).
\end{theorem}

\begin{proof}
The proof is conducted in the following section \S\
\ref{ss:Rank2}.
\end{proof}

The proof of Theorem
\ref{thm:Alpha2Ramified}
and Theorem
\ref{thm:Alpha1Ramified}
is conducted in the following sections
on real quadratic base fields with \(3\)-class rank \(1\) and \(2\).

\subsection{Real quadratic base fields with \(3\)-class rank \(1\)}
\label{ss:Rank1}

\noindent
In Table
\ref{tbl:Rank1},
we present the results of our search
for the \textit{minimal discriminant} \(d_L\), resp. \(d_N\),
of a non-Galois totally real cubic field \(L\), resp. its normal closure \(N\),
with \textit{differential principal factorization type} \(\alpha_2\).
Since \(\varrho=1\),
the unramified component is a \textit{singlet}, which must be of DPF type \(\delta_1\).
Since \(t=s=1\),
the DPF types \(\alpha_2,\beta_1,\beta_2,\delta_1,\delta_2,\varepsilon\) would be possible,
for each member of the ramified \textit{triplet},
but only the types \(\alpha_2,\delta_1,\delta_2\) occur usually.

The desired minimum is clearly given by \(d_L=19^2\cdot 2\,677=966\,397\)
with two occurrences of ramified extensions having DPF type \(\alpha_2\).
For \(f=3^2\), the condition \(d\equiv 1\,(\mathrm{mod}\,3)\) is required.

\renewcommand{\arraystretch}{1.1}

\begin{table}[ht]
\caption{Heterogeneous \textbf{quartets} of \(S_3\)-fields with splitting prime (power) \(f\)}
\label{tbl:Rank1}
\begin{center}
\begin{tabular}{|r|r|r||c||ccc|}
\hline
        &             &                      & unramified component & \multicolumn{3}{c|}{ramified components}   \\
 \(f\)  & \(d\)       & \(d_L=f^2\cdot d\)   & \(\delta_1\)         & \(\alpha_2\) & \(\delta_1\) & \(\delta_2\) \\
\hline
\(3^2\) & \(14\,197\) &  \(1\,149\,957\)     & \(1\)                & \(3\)        & \(0\)        & \(0\)        \\
 \(7\)  & \(21\,781\) &  \(1\,067\,269\)     & \(1\)                & \(2\)        & \(1\)        & \(0\)        \\
 \(13\) &  \(9\,749\) &  \(1\,647\,581\)     & \(1\)                & \(2\)        & \(0\)        & \(1\)        \\
 \(19\) &  \(2\,677\) &\(\mathbf{966\,397}\) & \(1\)                & \(2\)        & \(1\)        & \(0\)        \\
 \(31\) &  \(3\,877\) &  \(3\,725\,797\)     & \(1\)                & \(2\)        & \(0\)        & \(1\)        \\
 \(37\) &  \(5\,477\) &  \(7\,498\,013\)     & \(1\)                & \(1\)        & \(0\)        & \(2\)        \\
 \(43\) &  \(4\,933\) &  \(9\,121\,117\)     & \(1\)                & \(3\)        & \(0\)        & \(0\)        \\
 \(61\) &  \(3\,981\) & \(14\,813\,301\)     & \(1\)                & \(3\)        & \(0\)        & \(0\)        \\
 \(67\) &  \(4\,493\) & \(20\,169\,077\)     & \(1\)                & \(2\)        & \(0\)        & \(1\)        \\
 \(73\) & \(10\,733\) & \(57\,196\,157\)     & \(1\)                & \(3\)        & \(0\)        & \(0\)        \\
\hline
\end{tabular}
\end{center}
\end{table}


\noindent
Since we know a small candidate \(d_L=966\,397\) for the minimal discriminant,
and since the smallest quadratic fundamental discriminant with \(\varrho=1\) is \(d=229\),
we only have to investigate prime and composite conductors \(f=\sqrt{\frac{d_L}{d_K}}\) with \(s\ge 1\) and
\[f\le\sqrt{\frac{966\,397}{229}}\approx\sqrt{4220}\approx 64.9,\]
which are divisible by a split prime, that is,
\[f\in\lbrace 7,9=3^2,13,14=2\cdot 7,18=2\cdot 3^2,19,21=3\cdot 7,26=2\cdot 13,31,35=5\cdot 7,37,\]
\[38=2\cdot 19,39=3\cdot 13,42=2\cdot 3\cdot 7,43,45=5\cdot 3^2,57=3\cdot 19,61,62=2\cdot 31,63=7\cdot 3^2\rbrace.\]

{\normalsize

\renewcommand{\arraystretch}{1.1}

\begin{table}[ht]
\caption{Heterogeneous \textbf{quartets} of \(S_3\)-fields with conductor \(f\), where \(s=1\)}
\label{tbl:Rank1Coarse}
\begin{center}
\begin{tabular}{|r|c||r|r||c||ccccc|}
\hline
                    &                    &             &                      & unramified component & \multicolumn{5}{c|}{ramified components}                               \\
 \(f\)              & condition          & \(d\)       & \(d_L=f^2\cdot d\)   & \(\delta_1\)         & \(\alpha_2\) & \(\beta_1\) & \(\beta_2\) & \(\delta_1\) & \(\delta_2\) \\
\hline
              \(7\) &                    & \(21\,781\) &      \(1\,067\,269\) & \(1\)                & \(2\)        & \(0\)       & \(0\)       & \(1\)        & \(0\)        \\
            \(3^2\) & \(d\equiv 1\,(3)\) & \(14\,197\) &      \(1\,149\,957\) & \(1\)                & \(3\)        & \(0\)       & \(0\)       & \(0\)        & \(0\)        \\
             \(13\) &                    &  \(9\,749\) &      \(1\,647\,581\) & \(1\)                & \(2\)        & \(0\)       & \(0\)       & \(0\)        & \(1\)        \\
             \(19\) &                    &  \(2\,677\) &\(\mathbf{966\,397}\) & \(1\)                & \(2\)        & \(0\)       & \(0\)       & \(1\)        & \(0\)        \\
             \(31\) &                    &  \(3\,877\) &      \(3\,725\,797\) & \(1\)                & \(2\)        & \(0\)       & \(0\)       & \(0\)        & \(1\)        \\
             \(37\) &                    &  \(5\,477\) &      \(7\,498\,013\) & \(1\)                & \(1\)        & \(0\)       & \(0\)       & \(0\)        & \(2\)        \\
             \(43\) &                    &  \(4\,933\) &      \(9\,121\,117\) & \(1\)                & \(3\)        & \(0\)       & \(0\)       & \(0\)        & \(0\)        \\
             \(61\) &                    &  \(3\,981\) &     \(14\,813\,301\) & \(1\)                & \(3\)        & \(0\)       & \(0\)       & \(0\)        & \(0\)        \\
\hline
       \(2\cdot 7\) &                    &  \(6\,997\) &      \(1\,371\,412\) & \(1\)                & \(3\)        & \(0\)       & \(0\)       & \(0\)        & \(0\)        \\
     \(2\cdot 3^2\) & \(d\equiv 1\,(3)\) & \(16\,141\) &      \(5\,229\,684\) & \(1\)                & \(3\)        & \(0\)       & \(0\)       & \(0\)        & \(0\)        \\
       \(3\cdot 7\) & \(d\equiv 3\,(9)\) & \(28\,137\) &     \(12\,408\,417\) & \(1\)                & \(3\)        & \(0\)       & \(0\)       & \(0\)        & \(0\)        \\
       \(3\cdot 7\) & \(d\equiv 6\,(9)\) & \(57\,516\) &     \(25\,364\,556\) & \(1\)                & \(3\)        & \(0\)       & \(0\)       & \(0\)        & \(0\)        \\
      \(2\cdot 13\) &                    & \(21\,557\) &     \(14\,572\,532\) & \(1\)                & \(3\)        & \(0\)       & \(0\)       & \(0\)        & \(0\)        \\
       \(5\cdot 7\) &                    & \(14\,457\) &     \(17\,709\,825\) & \(1\)                & \(3\)        & \(0\)       & \(0\)       & \(0\)        & \(0\)        \\
      \(2\cdot 19\) &                    & \(13\,765\) &     \(19\,876\,660\) & \(1\)                & \(3\)        & \(0\)       & \(0\)       & \(0\)        & \(0\)        \\
      \(3\cdot 13\) & \(d\equiv 3\,(9)\) & \(51\,528\) &     \(78\,374\,088\) & \(1\)                & \(3\)        & \(0\)       & \(0\)       & \(0\)        & \(0\)        \\
      \(3\cdot 13\) & \(d\equiv 6\,(9)\) & \(37\,176\) &     \(56\,544\,696\) & \(1\)                & \(3\)        & \(0\)       & \(0\)       & \(0\)        & \(0\)        \\
\(2\cdot 3\cdot 7\) & \(d\equiv 3\,(9)\) &\(891\,237\) & \(1\,572\,142\,068\) & \(1\)                & \(4\)        & \(1\)       & \(1\)       & \(0\)        & \(0\)        \\
\(2\cdot 3\cdot 7\) & \(d\equiv 6\,(9)\) &\(474\,261\) &    \(836\,596\,404\) & \(1\)                & \(2\)        & \(0\)       & \(1\)       & \(0\)        & \(0\)        \\
     \(5\cdot 3^2\) & \(d\equiv 1\,(3)\) & \(24\,952\) &     \(50\,527\,800\) & \(1\)                & \(1\)        & \(0\)       & \(0\)       & \(1\)        & \(1\)        \\
      \(3\cdot 19\) & \(d\equiv 3\,(9)\) & \(24\,393\) &     \(79\,252\,857\) & \(1\)                & \(3\)        & \(0\)       & \(0\)       & \(0\)        & \(0\)        \\
      \(3\cdot 19\) & \(d\equiv 6\,(9)\) & \(39\,417\) &    \(128\,065\,833\) & \(1\)                & \(3\)        & \(0\)       & \(0\)       & \(0\)        & \(0\)        \\
      \(2\cdot 31\) &                    &  \(7\,573\) &     \(29\,110\,612\) & \(1\)                & \(3\)        & \(0\)       & \(0\)       & \(0\)        & \(0\)        \\
     \(7\cdot 3^2\) & \(d\equiv 1\,(3)\) &  \(2\,941\) &     \(11\,672\,829\) & \(1\)                & \(3\)        & \(0\)       & \(0\)       & \(0\)        & \(0\)        \\
     \(7\cdot 3^2\) & \(d\equiv 2\,(3)\) & \(23\,993\) &     \(95\,228\,217\) & \(1\)                & \(3\)        & \(0\)       & \(0\)       & \(0\)        & \(0\)        \\
\hline
\end{tabular}
\end{center}
\end{table}

}

\noindent
The result of the investigations is summarized in Table
\ref{tbl:Rank1Coarse},
which clearly shows that \(d_L=\mathbf{966\,397}\), for \(d=2\,677\) and splitting prime conductor \(f=19\)
bigger than the conductor \(f=1\) of unramified extensions \(N/K\),
is the desired \textbf{minimal discriminant} of a totally real cubic field with
ramified extension \(N/K\), DPF type \(\alpha_2\) and \(U_N=U_0\).
The information has been computed with Fieker's class field theoretic routines of Magma
\cite{Fi2001,MAGMA2020}.


\subsection{Real quadratic base fields with \(3\)-class rank \(2\)}
\label{ss:Rank2}
\noindent
In this situation, the unramified \textit{quartet} is non-trivial,
since two DPF types \(\alpha_1\) and \(\delta_1\) are possible.
These quartets have been thoroughly studied in
\cite{Ma2012},
and in Table
\ref{tbl:Rank2Cond2}
and
\ref{tbl:Rank2Cond5},
we use the corresponding notation for \textit{capitulation types}.

In Table
\ref{tbl:Rank2Cond2},
we present the results of the crucial search
for the \textit{minimal discriminant} \(d_L\), resp. \(d_N\),
of a non-Galois totally real cubic field \(L\), resp. its normal closure \(N\),
with \textit{differential principal factorization type} \(\alpha_1\)
such that \(N/K\) is a \textit{ramified} extension of a
real quadratic field \(K\) with \(3\)-class rank \(\varrho=2\).
We tried to fix the minimal possible conductor \(f>1\), namely \(f=2\).
This experiment was motivated by the fact that
the conductor \(f\) enters the expression \(d_L=f^2\cdot d\) in its second power,
whereas the quadratic discriminant \(d\) enters linearly.
Consequently, the probability to find the minimum of \(d_L\)
is higher for small \(f\) than for small \(d\).

The table is ordered by increasing quadratic fundamental discriminants \(d\)
and gives \(d_L=2^2\cdot d\) and the \textit{Artin pattern} \((\varkappa,\tau)\)
of the \textit{heterogeneous tridecuplet} of cyclic cubic relative extensions \(N/K\)
consisting of an \textit{unramified quartet} \((N_{1,1},\ldots,N_{1,4})\) with conductor \(f^\prime=1\)
and a \textit{ramified nonet} \((N_{2,1},\ldots,N_{2,9})\) with conductor \(f=2\),
grouped by the possible two, resp. four, DPF types \(\alpha_1,\delta_1\), resp. \(\alpha_1,\beta_1,\delta_1,\varepsilon\).
Transfer kernels \(\varkappa\) are abbreviated by digits,
\(0\) for two-dimensional and
\(1,\ldots,4\) for one-dimensional principalization,
and an asterisk \(\ast\) for a trivial kernel.
Transfer targets \(\tau\) are abbreviated
by logarithmic abelian type invariants of \(3\)-class groups.
Symbolic exponents always denote iteration.

The desired minimum is given by \(d_L=4\cdot 4\,562\,765=\mathbf{18\,251\,060}\)
with five occurrences of ramified extensions with DPF type \(\alpha_1\).
Generally,
there is an abundance of ramified extensions with two-dimensional capitulation kernel:
at least three and at most all nine of a nonet.


\renewcommand{\arraystretch}{1.1}

\begin{table}[ht]
\caption{Artin pattern \((\varkappa,\tau)\) of heterogeneous multiplets modulo \(f=2\)}
\label{tbl:Rank2Cond2}
\begin{center}
\begin{tabular}{|r||l|cc|cc||cc|cc|cc|cc|}
\hline
                 & \multicolumn{5}{c||}{unramified components} & \multicolumn{8}{c|}{ramified components} \\
                 &      & \multicolumn{2}{c|}{\(\alpha_1\)} & \multicolumn{2}{c||}{\(\delta_1\)} & \multicolumn{2}{c|}{\(\alpha_1\)} & \multicolumn{2}{c|}{\(\beta_1\)} & \multicolumn{2}{c|}{\(\delta_1\)} & \multicolumn{2}{c|}{\(\varepsilon\)} \\
         \(d_K\) & Type & \(\varkappa\) & \(\tau\)          & \(\varkappa\) & \(\tau\)           & \(\varkappa\) & \(\tau\)          & \(\varkappa\) & \(\tau\)         & \(\varkappa\) & \(\tau\)          & \(\varkappa\) & \(\tau\)             \\
\hline
 \(\mathbf{4\,562\,765}\) & a.\(3^\ast\) & \(0^3\) & \((1^2)^3\) & \(1\)    & \(1^3\)        & \(0^5\) & \(2^21^2,(1^4)^4\) & \(1\) & \(1^5\)  & \(14^2\)    & \((21^3)^3\)     &       &          \\
 \(7\,339\,397\) & a.\(3^\ast\) & \(0^3\) & \((1^2)^3\) & \(1\)    & \(1^3\)        & \(0^7\) & \((1^4)^7\)        & \(2\) & \(21^3\) & \(1\)       & \(21^3\)         &       &          \\
 \(7\,601\,461\) & a.\(3\)      & \(0^3\) & \((1^2)^3\) & \(1\)    & \(21\)         & \(0^6\) & \(2^21^2,(1^4)^5\) &       &          & \(234\)     & \((21^3)^3\)     &       &          \\
 \(7\,657\,037\) & a.\(3\)      & \(0^3\) & \((1^2)^3\) & \(1\)    & \(21\)         & \(0^6\) & \((1^4)^6\)        & \(1\) & \(21^3\) & \(12\)      & \(1^5,21^3\)     &       &          \\
 \(7\,736\,749\) & a.\(3^\ast\) & \(0^3\) & \((1^2)^3\) & \(1\)    & \(1^3\)        & \(0^7\) & \((1^4)^7\)        &       &          & \(4^2\)     & \((21^3)^2\)     &       &          \\
 \(8\,102\,053\) & a.\(3^\ast\) & \(0^3\) & \((1^2)^3\) & \(1\)    & \(1^3\)        & \(0^7\) & \((1^4)^7\)        &       &          & \(23\)      & \(1^5,21^3\)     &       &          \\
 \(9\,182\,229\) & a.\(2\)      & \(0^3\) & \((1^2)^3\) & \(4\)    & \(21\)         & \(0^8\) & \(2^21^2,(1^4)^7\) &       &          & \(2\)       & \(21^3\)         &       &          \\
 \(9\,500\,453\) & a.\(3\)      & \(0^3\) & \((1^2)^3\) & \(1\)    & \(21\)         & \(0^8\) & \(2^21^2,(1^4)^7\) &       &          & \(3\)       & \(21^3\)         &       &          \\
 \(9\,533\,357\) & a.\(3\)      & \(0^3\) & \((1^2)^3\) & \(1\)    & \(21\)         & \(0^6\) & \((1^4)^6\)        & \(1\) & \(21^3\) & \(23\)      & \((21^3)^2\)     &       &          \\
\(11\,003\,845\) & a.\(3\)      & \(0^3\) & \((1^2)^3\) & \(1\)    & \(21\)         & \(0^4\) & \((1^4)^4\)        &       &          & \(12^24^2\) & \(1^5,(21^3)^4\) &       &          \\
\(12\,071\,253\) & a.\(3\)      & \(0^3\) & \((1^2)^3\) & \(1\)    & \(21\)         & \(0^7\) & \((1^4)^7\)        & \(3\) & \(21^3\) & \(2\)       & \(21^3\)         &       &          \\
\(14\,266\,853\) & a.\(3\)      & \(0^3\) & \((1^2)^3\) & \(1\)    & \(21\)         & \(0^8\) & \(2^21^2,(1^4)^7\) &       &          & \(4\)       & \(21^3\)         &       &          \\
\(14\,308\,421\) & a.\(3^\ast\) & \(0^3\) & \((1^2)^3\) & \(1\)    & \(1^3\)        & \(0^4\) & \((1^4)^4\)        &       &          & \(1^2234\)  & \(2^31,(21^3)^3,1^5\) &  &          \\
\(14\,315\,765\) & a.\(3\)      & \(0^3\) & \((1^2)^3\) & \(1\)    & \(21\)         & \(0^7\) & \((1^4)^7\)        &       &          & \(23\)      & \((21^3)^2\)     &       &          \\
\(14\,395\,013\) & a.\(3^\ast\) & \(0^3\) & \((1^2)^3\) & \(1\)    & \(1^3\)        & \(0^6\) & \((1^4)^6\)        & \(1\) & \(21^3\) & \(23\)      & \((21^3)^2\)     &       &          \\
\(15\,131\,149\) & D.\(10\)     &         &             & \(2414\) & \((21)^3,1^3\) & \(0^7\) & \((1^4)^7\)        & \(1\) & \(21^3\) & \(1\)       & \(21^3\)         &       &          \\
\(16\,385\,741\) & a.\(3^\ast\) & \(0^3\) & \((1^2)^3\) & \(1\)    & \(1^3\)        & \(0^4\) & \((1^4)^4\)        &       &          & \(23^24\)   & \((21^3)^4\)     & \(\ast\) & \(32^21\) \\
\hline
\end{tabular}
\end{center}
\end{table}


Table
\ref{tbl:Rank2Cond5}
shows analogous results for the conductor \(f=5\),
that is, \(d_L=5^2\cdot d\).
The minimum \(d_L=25\cdot 1\,049\,512=26\,237\,800\)
is clearly beaten by the minimum \(4\cdot 4\,562\,765=18\,251\,060\) in Table
\ref{tbl:Rank2Cond2}.


\renewcommand{\arraystretch}{1.1}

\begin{table}[ht]
\caption{Artin pattern \((\varkappa,\tau)\) of heterogeneous multiplets modulo \(f=5\)}
\label{tbl:Rank2Cond5}
\begin{center}
\begin{tabular}{|r||l|cc|cc||cc|cc|cc|cc|}
\hline
                 & \multicolumn{5}{c||}{unramified components} & \multicolumn{8}{c|}{ramified components} \\
                 &      & \multicolumn{2}{c|}{\(\alpha_1\)} & \multicolumn{2}{c||}{\(\delta_1\)} & \multicolumn{2}{c|}{\(\alpha_1\)} & \multicolumn{2}{c|}{\(\beta_1\)} & \multicolumn{2}{c|}{\(\delta_1\)} & \multicolumn{2}{c|}{\(\varepsilon\)} \\
         \(d_K\) & Type & \(\varkappa\) & \(\tau\)          & \(\varkappa\) & \(\tau\)           & \(\varkappa\) & \(\tau\)          & \(\varkappa\) & \(\tau\)         & \(\varkappa\) & \(\tau\)          & \(\varkappa\) & \(\tau\)             \\
\hline
 \(1\,049\,512\) & a.\(3\)      & \(0^3\) & \((1^2)^3\) & \(1\)    & \(21\)         & \(0^4\) & \((1^4)^4\)            &       &          & \(234^3\) & \((21^3)^5\) &       &          \\
 \(2\,461\,537\) & a.\(2\)      & \(0^3\) & \((1^2)^3\) & \(4\)    & \(21\)         & \(0^7\) & \((1^4)^7\)            &       &          & \(12\)    & \((21^3)^2\) &       &          \\
 \(2\,811\,613\) & a.\(3^\ast\) & \(0^3\) & \((1^2)^3\) & \(1\)    & \(1^3\)        & \(0^5\) & \(2^21^2,(1^4)^4\)     & \(2\) & \(21^3\) & \(123\)   & \((21^3)^3\) &       &          \\
 \(3\,091\,133\) & a.\(3\)      & \(0^3\) & \((1^2)^3\) & \(1\)    & \(21\)         & \(0^4\) & \((1^4)^4\)            & \(4\) & \(21^3\) & \(1^32\)  & \((21^3)^4\) &       &          \\
 \(5\,858\,753\) & G.\(19\)     &         &             & \(2143\) & \((21)^4\)     & \(0^7\) & \((2^21^2)^3,(1^4)^4\) &       &          & \(3\)     & \(21^3\)     & \(\ast\) & \(21^4\)  \\
 \(6\,036\,188\) & D.\(10\)     &         &             & \(3431\) & \(1^3,(21)^3\) & \(0^8\) & \((1^4)^8\)            &       &          &           &              & \(\ast\) & \(2^21^2\) \\
\hline
\end{tabular}
\end{center}
\end{table}


Since we know a small candidate \(d_L=18\,251\,060\) for the minimal discriminant,
and since the smallest quadratic discriminant with \(\varrho=2\) is \(d=32\,009\),
we only have to investigate prime and composite conductors \(f=\sqrt{\frac{d_L}{d_K}}\) with
\[f\le\sqrt{\frac{18\,251\,060}{32\,009}}\approx\sqrt{570.2}\approx 23.9,\]
that is,
\[f\in\lbrace 2,3,5,6=2\cdot 3,7,9=3^2,10=2\cdot 5,11,13,14=2\cdot 7,\]
\[15=3\cdot 5,17,18=2\cdot 3^2,19,21=3\cdot 7,22=2\cdot 11,23\rbrace.\]

\noindent
The result of the investigations is summarized in Table
\ref{tbl:Rank2Coarse},
which clearly shows that \(d_L=\mathbf{18\,251\,060}\), for \(d=4\,562\,765\) and the smallest possible conductor \(f=2\)
bigger than the conductor \(f=1\) of unramified extensions \(N/K\),
is the desired \textbf{minimal discriminant} of a totally real cubic field with
ramified extension \(N/K\), DPF type \(\alpha_1\) and \(U_N=U_0\).
The information has been computed with Fieker's class field theoretic routines of Magma
\cite{Fi2001,MAGMA2020}.

\newpage

{\normalsize

\renewcommand{\arraystretch}{1.1}

\begin{table}[ht]
\caption{Heterogeneous \textbf{tridecuplets} of \(S_3\)-fields with conductor \(f\)}
\label{tbl:Rank2Coarse}
\begin{center}
\begin{tabular}{|r|c||r|r||cc||cccc|}
\hline
               &                    &                 &                      & \multicolumn{2}{c||}{unramified components} & \multicolumn{4}{c|}{ramified components} \\
 \(f\)         & condition          & \(d\)           & \(d_L=f^2\cdot d\)   & \(\alpha_1\) & \(\delta_1\) & \(\alpha_1\) & \(\beta_1\) & \(\delta_1\) & \(\varepsilon\) \\
\hline
         \(2\) &            &\(\mathbf{4\,562\,765}\)&\(\mathbf{18\,251\,060}\)& \(3\)        & \(1\)        & \(5\)        & \(1\)       & \(3\)        & \(0\)           \\
         \(3\) & \(d\equiv 3\,(9)\) & \(9\,964\,821\) &     \(89\,683\,389\) & \(3\)        & \(1\)        & \(4\)        & \(0\)       & \(4\)        & \(1\)           \\
         \(5\) &                    & \(1\,049\,512\) &     \(26\,237\,800\) & \(3\)        & \(1\)        & \(4\)        & \(0\)       & \(5\)        & \(0\)           \\
         \(7\) &                    &    \(966\,053\) &     \(47\,336\,597\) & \(3\)        & \(1\)        & \(4\)        & \(0\)       & \(4\)        & \(1\)           \\
       \(3^2\) & \(d\equiv 1\,(3)\) & \(1\,482\,568\) &    \(120\,088\,008\) & \(3\)        & \(1\)        & \(5\)        & \(1\)       & \(2\)        & \(1\)           \\
       \(3^2\) & \(d\equiv 2\,(3)\) & \(2\,515\,388\) &    \(203\,746\,428\) & \(3\)        & \(1\)        & \(6\)        & \(1\)       & \(2\)        & \(0\)           \\
       \(3^2\) & \(d\equiv 6\,(9)\) &    \(621\,429\) &     \(50\,335\,749\) & \(3\)        & \(1\)        & \(6\)        & \(0\)       & \(3\)        & \(0\)           \\
        \(11\) &                    &    \(476\,152\) &     \(57\,614\,392\) & \(3\)        & \(1\)        & \(7\)        & \(0\)       & \(2\)        & \(0\)           \\
        \(13\) &                    & \(1\,122\,573\) &    \(189\,714\,837\) & \(3\)        & \(1\)        & \(7\)        & \(0\)       & \(2\)        & \(0\)           \\
        \(17\) &                    &    \(665\,832\) &    \(192\,425\,848\) & \(3\)        & \(1\)        & \(7\)        & \(0\)       & \(2\)        & \(0\)           \\
        \(19\) &                    &    \(635\,909\) &    \(229\,563\,149\) & \(3\)        & \(1\)        & \(5\)        & \(3\)       & \(1\)        & \(0\)           \\
        \(23\) &                    &    \(390\,876\) &    \(206\,773\,404\) & \(3\)        & \(1\)        & \(7\)        & \(1\)       & \(1\)        & \(0\)           \\
\hline
  \(2\cdot 3\) & \(d\equiv 3\,(9)\) & \(5\,963\,493\) &    \(214\,685\,748\) & \(3\)        & \(1\)        & \(7\)        & \(2\)       & \(0\)        & \(0\)           \\
  \(2\cdot 3\) & \(d\equiv 6\,(9)\) & \(4\,305\,957\) &    \(155\,014\,452\) & \(0\)        & \(4\)        & \(6\)        & \(3\)       & \(0\)        & \(0\)           \\
  \(2\cdot 5\) &                    &    \(363\,397\) &     \(36\,339\,700\) & \(3\)        & \(1\)        & \(6\)        & \(3\)       & \(0\)        & \(0\)           \\
  \(2\cdot 7\) &                    &    \(358\,285\) &     \(70\,223\,860\) & \(4\)        & \(0\)        & \(7\)        & \(2\)       & \(0\)        & \(0\)           \\
  \(3\cdot 5\) & \(d\equiv 3\,(9)\) & \(4\,845\,432\) & \(1\,090\,222\,200\) & \(3\)        & \(1\)        & \(6\)        & \(3\)       & \(0\)        & \(0\)           \\
  \(3\cdot 5\) & \(d\equiv 6\,(9)\) & \(1\,646\,817\) &    \(370\,533\,825\) & \(3\)        & \(1\)        & \(6\)        & \(3\)       & \(0\)        & \(0\)           \\
\(2\cdot 3^2\) & \(d\equiv 1\,(3)\) & \(2\,142\,445\) &    \(694\,152\,180\) & \(3\)        & \(1\)        & \(6\)        & \(3\)       & \(0\)        & \(0\)           \\
\(2\cdot 3^2\) & \(d\equiv 2\,(3)\) &    \(635\,909\) &    \(206\,034\,516\) & \(3\)        & \(1\)        & \(6\)        & \(3\)       & \(0\)        & \(0\)           \\
\(2\cdot 3^2\) & \(d\equiv 6\,(9)\) & \(2\,538\,285\) &    \(822\,404\,340\) & \(3\)        & \(1\)        & \(6\)        & \(3\)       & \(0\)        & \(0\)           \\
  \(3\cdot 7\) & \(d\equiv 3\,(9)\) & \(3\,597\,960\) & \(1\,586\,700\,360\) & \(3\)        & \(1\)        & \(6\)        & \(3\)       & \(0\)        & \(0\)           \\
  \(3\cdot 7\) & \(d\equiv 6\,(9)\) & \(3\,122\,232\) & \(1\,376\,904\,312\) & \(0\)        & \(4\)        & \(6\)        & \(3\)       & \(0\)        & \(0\)           \\
 \(2\cdot 11\) &                    & \(2\,706\,373\) & \(1\,309\,884\,532\) & \(3\)        & \(1\)        & \(6\)        & \(3\)       & \(0\)        & \(0\)           \\
\hline
\end{tabular}
\end{center}
\end{table}

}


\subsection{Scholz conjecture for \(p\ge 5\)}
\label{ss:DihedralScholz}

\noindent
We have been curious if the conjecture of Scholz
can also be verified for dihedral fields \(N/\mathbb{Q}\) of degrees \(10\) and \(14\).
This is indeed the case,
and the root discriminants \(f^2\cdot d\) in the following theorem are probably minimal.

\begin{theorem}
\label{thm:ScholzDihedral}
Let \(p\) be an odd prime number.
Suppose \(L\) is a non-Galois number field of degree \(p\)
with totally real absolutely dihedral Galois closure \(N\) of degree \(2p\),
and let \(K\) be the unique real quadratic subfield of \(N\).
Then \(N\) satisfies the condition \(U_N=U_0:=\langle U_K,U_L,U_{L^{(1)}},\ldots,U_{L^{(p-1)}}\rangle\),
\begin{enumerate}
\item
if \(d_L=(f^2\cdot d)^2\) with \(f=11\cdot 31\), \(d=5\), \(f^2\cdot d=581\,405\), when \(p=5\),
\item
if \(d_L=(f^2\cdot d)^3\) with \(f=29\cdot 43\), \(d=13\), \(f^2\cdot d=20\,215\,117\), when \(p=7\).
\end{enumerate}
In both cases, the conductor is of the form \(f=\ell_1\cdot\ell_2\)
with prime numbers \(\ell_i\equiv +1\,(\mathrm{mod}\,p)\) which split in \(K\),
and \(L\) is a singlet with differential principal factorization type \(\tau(L)=\alpha_3\).
\end{theorem}

\begin{proof}
By immediate inspection of real quadratic fields \(K\)
with \(p\)-class rank \(\varrho_p=0\)
and \(p\)-admissible conductors \(f\),
divisible by two primes which split in \(K\),
with the aid of Magma.
\end{proof}


\section{Conclusion}
\label{Conclusion}

\noindent
In this paper, we have given the \textit{complete classification} of all multiplets
of totally real cubic fields \(L\) in the range \(0<d_L<10^7\) of Llorente and Quer
\cite{LlQu1988}
according to their \textit{differential principal factorizations}
(Tables \ref{tbl:LlorenteQuer0} and \ref{tbl:LlorenteQuer1}).
Inspired by discussions after our two presentations at the
West Coast Number Theory Conference in Asilomar, December \(1990\),
we had attempted this classification in August \(1991\) already,
but we were forced to restrict the range to the upper bound \(2\cdot 10^5\) in
\cite{Ma1991c}.
In spite of the required correction of \(14\) errors
(Tables \ref{tbl:Translation} and \ref{tbl:Corrections}),
the table
\cite{Ma1991c}
and the associated theory
\cite{Ma1991b}
were a \textit{masterpiece of outstanding innovations}
concerning DPF types of multiplets of dihedral fields
and a \textit{role model} for the present paper
and its predecessor
\cite{Ma2021}.

We have also given the complete verification of the \textit{Conjecture of Arnold Scholz}
(Conjecture \ref{cnj:Scholz}).
It was necessary to develop the new concept of \textit{relative} principal factorizations
in order to illuminate the full reach of this conjecture,
which we have reformulated more ostensively in Conjecture
\ref{cnj:Mayer}.
Due to the computational challenges,
the proof of each of the different perspectives of the conjecture
was established many years after Scholz's paper in \(1933\)
\cite{So1933}:
Corollary \ref{cor:Alpha1Unramified} on \(f=1\) was proved
\(49\) years later in \(1982\)
\cite{HeSm1982},
Theorem \ref{thm:Alpha3Ramified} concerning the type \(\alpha_3\) singulet
\(58\) years later in \(1991\)
\cite{Ma1991c}, 
Theorem \ref{thm:Alpha2Ramified} on the triplet containing type \(\alpha_2\)
even \(84\) years later on 19 November \(2017\), and
Theorem \ref{thm:Alpha1Ramified} on the nonet containing type \(\alpha_1\) with \(f>1\)
also \(84\) years later on 23 November \(2017\).

In our ultimate Table
\ref{tbl:Tendencies},
we emphasize the apparent \textit{asymptotic tendencies} of DPF types,
based on five ranges of discriminants \(0<d_L<B\)
with increasing upper bounds \(B\).
Relative frequencies are rounded to integer percentages.
It is striking that the normal closures \(N\) of an overwhelming proportion
with \(93\%\) of all totally real cubic fields \(L\)
have a unit group \(U_N\) which is a
\textit{non-split extension} of \(U_K=\langle -1,\eta\rangle\)
if considered as a module over the integral group algebra \(\mathbb{Z}\lbrack S_3\rbrack\),
since it contains a unit \(H\) such that \(N_{N/K}(H)=H\cdot H^\sigma\cdot H^{\sigma^2}=\eta\),
according to Remark
\ref{rmk:RealDPFTypes}.

\noindent
This phenomenon is due to \textit{extremely dominating}
unramified extensions \(N/K\) with conductor \(f=1\), \(\varrho=1\),
and mandatory type \(\delta_1\) (\(72\%\)), and
ramified extensions \(N/K\) with regular prime(power) conductor \(f\), \(\varrho=0\),
and mandatory type \(\varepsilon\) (\(21\%\)).
In contrast,
the contributions by the rare types \(\alpha_2\) and \(\alpha_3\)
and by the cyclic cubic fields \(\zeta\) are in fact \textit{negligible}.
In spite of its distinctive dominance for conductors with two or more prime divisors,
type \(\gamma\) remains \textit{marginal} with \(2\%\).
Other marginal (but not negligible) contributions arise from type \(\alpha_1\),
due to increasing occurrences of \(\varrho=2\),
from type \(\beta_1\), due to capitulation in ramified extensions with \(\varrho=1\),
and from the types \(\beta_2\) and \(\delta_2\),
due to conductors \(f\) with a prime divisor which splits in the quadratic subfield \(K<N\).

\renewcommand{\arraystretch}{1.1}

\begin{table}[ht]
\caption{Tendencies of the statistical distribution of DPF types}
\label{tbl:Tendencies}
\begin{center}
\begin{tabular}{|c||rr||rr||rr||rr||rr|}
\hline
 \(B\) & \multicolumn{2}{c||}{\(1500\)}& \multicolumn{2}{c||}{\(10^5\)}&\multicolumn{2}{c||}{\(2\cdot 10^5\)}&\multicolumn{2}{c||}{\(5\cdot 10^5\)}& \multicolumn{2}{c|}{\(10^7\)} \\
 Type           & \(\#\) & \(\%\) & \(\#\)   & \(\%\) & \(\#\)   & \(\%\) & \(\#\)    & \(\%\) & \(\#\)     & \(\%\) \\
\hline
\(\alpha_1\)    &  \(0\) &  \(0\) &   \(16\) &  \(0\) &   \(50\) &  \(1\) &   \(175\) &  \(1\) &   \(7951\) &  \(1\) \\
\(\alpha_2\)    &  \(0\) &  \(0\) &    \(0\) &  \(0\) &    \(0\) &  \(0\) &     \(0\) &  \(0\) &    \(142\) &  \(0\) \\
\(\alpha_3\)    &  \(0\) &  \(0\) &    \(0\) &  \(0\) &    \(1\) &  \(0\) &     \(3\) &  \(0\) &    \(122\) &  \(0\) \\
\(\beta_1\)     &  \(0\) &  \(0\) &   \(10\) &  \(0\) &   \(21\) &  \(0\) &    \(89\) &  \(0\) &   \(3924\) &  \(1\) \\
\(\beta_2\)     &  \(0\) &  \(0\) &   \(76\) &  \(2\) &  \(155\) &  \(2\) &   \(380\) &  \(1\) &   \(7639\) &  \(1\) \\
\(\gamma\)      &  \(2\) &  \(4\) &  \(106\) &  \(2\) &  \(201\) &  \(2\) &   \(493\) &  \(2\) &   \(9420\) &  \(2\) \\
\hline
\(\delta_1\)    & \(26\) & \(59\) & \(3349\) & \(70\) & \(7028\) & \(70\) & \(18714\) & \(71\) & \(426972\) & \(72\) \\
\(\delta_2\)    &  \(0\) &  \(0\) &   \(79\) &  \(2\) &  \(188\) &  \(2\) &   \(490\) &  \(2\) &  \(11128\) &  \(2\) \\
\(\varepsilon\) & \(10\) & \(23\) & \(1117\) & \(23\) & \(2301\) & \(23\) &  \(5986\) & \(23\) & \(125123\) & \(21\) \\
\hline
\(\zeta\)       &  \(6\) & \(14\) &   \(51\) &  \(1\) &   \(70\) &  \(1\) &   \(110\) &  \(0\) &    \(501\) &  \(0\) \\
\hline
\end{tabular}
\end{center}
\end{table}


\section{Acknowledgements}
\label{s:Thanks}

\noindent
The author gratefully acknowledges
that his research was supported by the Austrian Science Fund (FWF):
projects J0497-PHY and P26008-N25.



\end{document}